\documentclass[10pt,a4paper, BCOR15mm]{scrartcl}
\usepackage[utf8]{inputenc}
\usepackage[T1]{fontenc}
\usepackage[english]{babel}
\usepackage{csquotes}
\usepackage{amssymb}
\usepackage{epsfig}
\usepackage{amsthm}
\usepackage{enumerate}
\usepackage{color}
\usepackage{url}
\usepackage{here}
\usepackage{caption}
\usepackage{subcaption}
\usepackage{rotating}
\usepackage{tabularx}
\usepackage{appendix}
\usepackage{mathtools}
\usepackage{tabu}
\usepackage[capitalize]{cleveref}
\addtokomafont{sectioning}{\rmfamily} 
\newtheorem{theorem}{Theorem}[section]

\newtheorem{remark}[theorem]{Remark}
\newtheorem{algorithm}[theorem]{Algorithm}

\newcommand{\R}{\mathbb{R}}
\newcommand{\N}{\mathbb{N}}
\newcommand{\eps}{\varepsilon}

\DeclareMathOperator{\Span}{span}

\newcommand{\h}{\mathcal{H}}

\newcommand{\LO}{\mathrm{L}^{2}(\Omega)}
\DeclareMathOperator*{\argmin}{argmin} 
\DeclareMathOperator*{\argmax}{argmax}

\begin{document}

%Title
\begin{center}
\textbf{\Large Study on parameter choice methods for the RFMP with respect to downward continuation}
\end{center}
\begin{center}
M. Gutting\footnote[1]{Geomathematics Group, Department of Mathematics, University of Siegen, Emmy-Noether-Campus,
Walter-Flex-Str. 3, 57068 Siegen, Germany, E-Mail addresses: gutting@mathematik.uni-siegen.de, kretz@mathematik.uni-siegen.de, michel@mathematik.uni-siegen.de}, B. Kretz\footnotemark[1], V. Michel\footnotemark[1], R. Telschow\footnote[2]{Computational Science Center, University of Vienna, Oskar Morgenstern-Platz 1, 1090 Vienna, Austria, E-Mail address: roger.telschow@univie.ac.at}
\end{center}

%Front
\begin{center}
\textbf{Abstract}
\end{center}
\begin{abstract}
Recently, the regularized functional matching pursuit (RFMP) was introduced as a greedy algorithm for linear ill-posed inverse problems.
This algorithm incorporates the Tikhonov-Phillips regularization which implies the necessity of a parameter choice. In this paper, some known parameter choice methods are evaluated with respect to their performance in the RFMP and its enhancement, the regularized orthogonal functional matching pursuit (ROFMP). As an example of a linear inverse problem, the downward continuation of gravitational field data from the satellite orbit to the Earth's surface is chosen, because it is exponentially ill-posed. For the test scenarios, different satellite heights with several noise-to-signal ratios and kinds of noise are combined. The performances of the parameter choice strategies in these scenarios are analyzed. For example, it is shown that a strongly scattered set of data points is an essentially harder challenge for the regularization than a regular grid. The obtained results yield a first orientation which parameter choice methods are feasible for the RFMP and the ROFMP.
\end{abstract}



%Intro 
\section{Introduction}

The gravitational field of the Earth is an important reference in the geosciences. It is an indicator for mass transports and mass reallocation on the Earth's surface.
These displacements of masses can be caused by ocean currents, evaporation, changes of the groundwater level, ablating of continental ice sheets, changes in the mean sea level or climate change (see e.g.\ \cite{IlkFlury, kusche2014}).\\
However, it is difficult to model the gravitational field, because terrestrial measurements are not globally available. In addition, the points of measurement on the sea are more scattered than those on the continents. This has motivated the launch of satellite missions with a focus on the gravitational field (see e.g.\ \cite{goce, goce_grace, champ, gravity}). Naturally, those data are given at a satellite orbit, not on the Earth's surface.  Additionally, the measurements are only given pointwise and are afflicted with noise. The problem of getting the potential from the satellite orbit onto the Earth's surface is the so-called downward continuation problem, which is a severely unstable and, therefore, ill-posed inverse problem (see e.g.\ \cite{freeden_michel, schneider}). \\
Traditionally, the gravitational potential of the Earth has been represented in terms of orthogonal spherical polynomials (i.e.\ spherical harmonics $Y_{n,j}$, see e.g.  \cite{freeden_gutting, michel, mueller}) as in the case of the Earth Gravitational Model 2008 (EGM2008, see \cite{egm}).
An advantage of this representation is that the upward continuation operator $\Psi$ which maps a potential $F$ from the Earth's surface (which we assume here to be the unit sphere $\Omega$) to the orbit $r\Omega$ with $r>1$ has the singular value decomposition
\begin{align}
\left(\Psi F\right)\left(x\right) = \sum\limits_{n=0}^{\infty} \sum\limits_{j=-n}^{n} \left\langle F,Y_{n,j} \right\rangle_{\LO} r^{-n} \, Y_{n,j}^r \left(x \right),\label{psi}
\end{align}
where $Y_{n,j}^r(x)\coloneqq \frac{1}{r} Y_{n,j} \left( \frac{x}{r} \right)$, $x \in r\Omega$. 
Its inverse is, therefore, given by 
\begin{align}
\Psi^{+}G = \sum\limits_{n=0}^{\infty} \sum\limits_{j=-n}^{n} \left\langle G,Y_{n,j}^{r} \right\rangle_{\mathrm{L}^{2}(r\Omega)} r^n Y_{n,j} = \sum\limits_{n=0}^{\infty} \sum\limits_{j=-n}^{n} \left\langle G,Y_{n,j}^{r} \right\rangle_{\mathrm{L}^{2}(r\Omega)} \sigma_{n}^{-1} Y_{n,j}
\end{align}
in the sense of $\LO$ and for all $G \in \Psi(\LO) \subset \mathrm{L}^{2}(r\Omega)$. Note that the singular values of $\Psi^+$, which are given by $(\sigma_n^{-1})_n= (r^n)_n$, increase exponentially. For details, see \cite{schneider, telschow}. 

The outline of this paper is as follows. \cref{rfmp} deals with the RFMP and its enhancement, the ROFMP, which are used here for the regularization of $\Psi^{+}$.  For both algorithms, the essential theoretical results are recapitulated.
In \cref{methods}, the parameter choice methods under consideration for the RFMP and ROFMP are summarized and details of their implementation for the test cases are explained.
In \cref{spec_methods}, the relevant details of the considered test scenarios are outlined. \cref{comparison} analyzes and compares the results for the various parameter choice strategies.

%RFMP Intro
\section{RFMP} \label{rfmp}
In this section, we briefly resume the regularized functional matching pursuit (RFMP), which was introduced in \cite{fischer, fischer_michel, michel_handbook, michel_telschow}, and an orthogonalized modification of it (see \cite{michel_telschow_pp, telschow}). It is an algorithm for the regularization of linear inverse problems.

According to \cite{michel_orzlowski, michel_telschow, telschow}, we use an arbitrary Hilbert space $\h \subset \mathrm{L}^2(\Omega)$.\\
Let an operator $\mathcal{F}:\h \rightarrow \R^{l}$ be given which is continuous and linear. Concerning the downward continuation, we have a vector $\boldsymbol{y}\in \R^l$ of measurements at a satellite orbit, that means our data are given pointwise. The inverse problem consists of the determination of a function $F \in \h$ such that 
\begin{align}
\mathcal{F}F=\boldsymbol{y}=((\Psi F)(x_j))_{j=1,\dots,l}, \label{task_points}
\end{align}
 where $(x_j)_{j=1,\dots,l}$ is a set of points at satellite height. 
In the following, we use bold letters for vectors in $\R^l$.

To find an approximation for our function F, we need to have a set of trial functions $\mathcal{D} \subset \h \setminus \{0\}$, which we call the \emph{dictionary}.
Our unknown function $F$ is expanded in terms of dictionary elements, that means we can represent it as
$F=\sum_{k=1}^{\infty} \alpha_{k} d_{k} \mathrm{~with~} \alpha_{k} \in \R \mathrm{~and~} d_{k} \in \mathcal{D} \mathrm{~for~all~} k \in \N$.

%RFMP algorithm
\subsection{The algorithm} \label{rfmp_algorithm}
The idea of the RFMP is the iterative construction of a sequence of approximations $(F_{n})_{n}$. This means that we add a basis function $d_k$ of the dictionary to the approximation in each step. This basis function is furthermore equipped with a coefficient $\alpha_k$.

Since the considered inverse problem is ill-posed, we use the Tikhonov-Phillips regularization, that is, our task is to find a function $F$ which minimizes
\begin{align}
\left\|\boldsymbol{y}-\mathcal{F}F \right\|_{\R^l}^2 + \lambda \left\|F\right\|_{\h}^2.
\end{align}
That means, if we have the approximation $F_n$ up to step $n$, our greedy algorithm chooses $\alpha_{n+1} \in \R$ and $d_{n+1} \in \mathcal{D}$ such that
\begin{align}
\left\|\boldsymbol{y}-\mathcal{F}\left(F_{n}+\alpha_{n+1}d_{n+1}\right)\right\|_{\R^l}^{2} + \lambda \left\|F_{n}+\alpha_{n+1}d_{n+1}\right\|_{\h}^{2}
\end{align}
is minimized. Here, $\lambda>0$ is the regularization parameter.

We can state the following algorithm for the RFMP.
\begin{algorithm} \label{rfmp_algo}
Let $\boldsymbol{y} \in \R^l$ and an operator $\mathcal{F}:\h \rightarrow \R^l$ (linear and continuous) be given.
\begin{description}
\item[(1) Initialization] 
Set $n \coloneqq 0$, $F_{0} \coloneqq 0$ and $\boldsymbol{R}_{0} \coloneqq \boldsymbol{y}-\mathcal{F}F_{0}=\boldsymbol{y}$, choose a stopping criterion (we stop, if $\left\|\boldsymbol{R}_{n+1}\right\|<\varrho$ for a given $\varrho > 0$ or $\alpha_{n+1} < \delta$ for a given $\delta>0$ or $n+1 > N$ for a given $N \in \N$, see also \cref{sec41}), and choose a regularization parameter $\lambda \in \R^{+}$.
\item[(2) Iteration]
Build $F_{n+1} \coloneqq F_{n}+\alpha_{n+1}d_{n+1}$ such that the following is fulfilled:
\begin{align}
d_{n+1}& \coloneqq \argmax_{d \in \mathcal{D}} \frac{\left( \left\langle \boldsymbol{R}_{n},\mathcal{F}d\right\rangle_{\R^l} - \lambda \left\langle F_{n},d\right\rangle_{\h} \right)^{2}}{\left\|\mathcal{F}d\right\|_{\R^l}^{2} + \lambda \left\| d\right\|_{\h}^{2}}, \label{dn+1} \\
\alpha_{n+1}& \coloneqq \frac{\left\langle \boldsymbol{R}_{n},\mathcal{F}d_{n+1}\right\rangle_{\R^l}-\lambda \left\langle F_{n},d_{n+1}\right\rangle_{\h}}{\left\|\mathcal{F}d_{n+1}\right\|_{\R^l}^2+\lambda \left\|d_{n+1}\right\|_{\h}^{2}}.
\end{align}
Set $\boldsymbol{R}_{n+1} \coloneqq \boldsymbol{R}_{n}-\alpha_{n+1}\mathcal{F}d_{n+1}$.
\item[(3) Stopping criterion]
$F_{n+1}$ is the output, if the stopping criterion is fulfilled. Otherwise, increase $n$ and go to step 2.
\end{description}
\end{algorithm}

The maximization, which is necessary to get $d_{n+1}$, is implemented by evaluating the fraction for all $d \in \mathcal{D}$ in each iteration and picking a maximizer. Since many involved terms can be calculated in a preprocessing, the numerical expenses can be kept low (see \cite{michel_handbook}).
For a convergence proof of the RFMP, we refer to \cite{michel_orzlowski}. Briefly, under certain conditions, one can show that the sequence $(F_n)_n$ converges to the solution $F_\infty$ of the Tikhonov-regularized normal equation
\begin{align}
(\mathcal{F}^* \mathcal{F}+\lambda \mathcal{I}) F_\infty=\mathcal{F}^* y, \label{normalequ}
\end{align}
where $\mathcal{I}$ is the identity operator and $\mathcal{F}^*$ is the adjoint operator to $\mathcal{F}$.

%ROFMP
\subsection{ROFMP} \label{rofmp_algorithm}
The regularized orthogonal functional matching pursuit (ROFMP) is an advancement of the RFMP from the previous section.

The basic idea is to project the residual onto the span of the chosen vectors, i.e.,
\begin{align}
\mathcal{V}_n  \coloneqq  \Span\{ \mathcal{F} d_1,\dots,\mathcal{F} d_n \} \subset \R^l \,,
\end{align}
and then adjust the previously chosen coefficients in such a way that the residual is afterwards contained in the orthogonal complement of the span. Since this so-called backfitting (cf. \cite{mallat_zhang,pati_et_al}) might not be optimal, we implement the so-called prefitting (cf. \cite{vincent_bengio}), where the next function and all coefficients are chosen simultaneously to guarantee optimality at every single stage of the algorithm.
Moreover, let $\mathcal{W}_n \coloneqq \mathcal{V}_n^{\bot}$ and the orthogonal projections on $\mathcal{V}_{n}$ and $\mathcal{W}_{n}$ are denoted by $\mathcal{P}_{\mathcal{V}_{n}}$ and $\mathcal{P}_{\mathcal{W}_{n}}$, respectively.
All in all, our aim is to find
\begin{align}
\left(\alpha_{n+1},d_{n+1}\right)= \argmin_{\alpha \in \R,\, d \in \mathcal{D}} \left( \left\| \boldsymbol{R}_{n}-\alpha\mathcal{P}_{\mathcal{W}_{n}}\mathcal{F}d \right\|_{\R^l}^{2} + \lambda \left\| F_{n}-\alpha B_{n}(d)+\alpha d \right\|_{\h}^{2} \right). \label{rofmp_task}
\end{align}
Here, 
\begin{align}
\sum_{i=1}^n \beta_i^{(n)}(d)  \mathcal{F} d_i = \mathcal{P}_{\mathcal{V}_n} (\mathcal{F} d)  \label{rofmp1}
\end{align}
and, thereby, we set
\begin{align}
B_n (d)  \coloneqq \sum_{i=1}^n  \beta_i^{(n)} (d) d_i  \in\h. \label{rofmp2}
\end{align}

The updated coefficients for the expansion at step $n+1$ are given by
\begin{align}
\alpha_{i}^{(n+1)} \coloneqq & \alpha_{i}^{(n)}-\alpha_{n+1}\beta_{i}^{(n)} (d_{n+1}), \quad i=1,\dots,n, \label{alpha} \\
\alpha_{n+1}^{(n+1)} \coloneqq & \alpha_{n+1}. \label{alpha2}
\end{align}
\newpage
The ROFMP algorithm can be summarized as follows.

\begin{algorithm} \label{rofmp_algo}
Let a dictionary $\mathcal{D} \subset \mathcal{H}$, a data vector $\boldsymbol{y} \in \R^l$ and an operator $\mathcal{F}:\h \rightarrow \R^l$ (linear and continuous) be given.
\begin{description}
\item[(1) Initialization] 
Set $n \coloneqq 0$, $F_{0} \coloneqq 0$ and $\boldsymbol{R}_{0} \coloneqq \boldsymbol{y}$, choose a stopping criterion (we stop, if $\left\|\boldsymbol{R}_{n+1}\right\|< \varrho$ for a given $\varrho > 0$ or $\alpha_{n+1} < \delta$ for a given $\delta>0$ or $n+1 > N$ for a given $N \in \N$, see also \cref{sec41}), and choose a regularization parameter $\lambda \in \R^{+}$.
\item[(2) Iteration]
Choose a function
\begin{align}
d_{n+1} \coloneqq \argmax\limits_{d \in \mathcal{D}} \frac{\left(\left\langle \boldsymbol{R}_{n},\mathcal{P}_{\mathcal{W}_{n}}\mathcal{F}d\right\rangle_{\R^l} + \lambda \left( \left\langle F_{n},B_{n}(d)\right\rangle_{\h} - \left\langle F_{n},d\right\rangle_{\h} \right)\right)^2}{ \left\|\mathcal{P}_{\mathcal{W}_{n}} \mathcal{F}d\right\|_{\R^l}^2 + \lambda \left( \left\|d-B_{n}(d)\right\|_{\h}^2  \right)},
\end{align}
and calculate the corresponding coefficient
\begin{align}
\alpha_{n+1} \coloneqq \frac{\left\langle \boldsymbol{R}_{n},\mathcal{P}_{\mathcal{W}_{n}}\mathcal{F}d_{n+1}\right\rangle_{\R^l} + \lambda \left( \left\langle F_{n},B_{n}(d_{n+1})\right\rangle_{\h} - \left\langle F_{n},d_{n+1}\right\rangle_{\h} \right)}{ \left\|\mathcal{P}_{\mathcal{W}_{n}} \mathcal{F}d_{n+1}\right\|_{\R^l}^2 + \lambda \left\|d_{n+1}-B_{n}(d_{n+1})\right\|_{\h}^2},
\end{align}
where $B_{n}(d)$ is defined according to \eqref{rofmp1} and \eqref{rofmp2}. With the updated coefficients $\alpha_{i}^{(n+1)}\coloneqq\alpha_{i}^{(n)}-\alpha_{n+1}\beta_{i}^{(n)}(d_{n+1})$ for $i=1,\dots,n$, we set $\alpha_{n+1}^{(n+1)}\coloneqq \alpha_{n+1}$ and build $F_{n+1}\coloneqq \sum_{i=1}^{n+1} \alpha_{i}^{(n+1)}d_{i}$. Finally, update the residual $\boldsymbol{R}_{n+1} \coloneqq \boldsymbol{R}_{n}-\alpha_{n+1}\mathcal{P}_{\mathcal{W}_{n}}\mathcal{F}d_{n+1}$.
\item[(3) Stopping criterion] 
$F_{n+1}$ is the output, if the stopping criterion is fulfilled. Otherwise, increase $n$ and go to step 2.
\end{description}
\end{algorithm}
For practical details of the implementation, see \cite{telschow}.
\begin{remark} \label{stagn}
If we  choose $d_{i} \in \mathcal{D}$ and $\alpha_{i}$ as in \cref{rofmp_algo} and update $\alpha$ as in \eqref{alpha} and \eqref{alpha2}, we obtain for the regularized case $(\lambda >0)$ that $\boldsymbol{R}_{n}$ is, in general, not orthogonal to $\mathcal{V}_{n}$ for all $n \in \N_{0}$, that means
\begin{align}
\boldsymbol{R}_{n} \notin \mathcal{W}_{n}.
\end{align}
\end{remark}
In \cite{telschow}, it was shown that, with the assumptions from \cref{stagn}, there exists a number $N:=N(\lambda)$ such that
\begin{align}
\boldsymbol{R}_{n}=\boldsymbol{R}_{N} \mathrm{~for~all~} n \geq N.
\end{align}
That means we get a stagnation of the residual.
This is a problem for the ROFMP, because we cannot reconstruct a certain part of the signal which lies in $\mathcal{V}_{n}$. Therefore, we have to modify the algorithm to an iterated Tikhonov-Phillips regularization. That means we run the algorithm for a given number of iterations (in our case $K>0$), then break up the process and start the algorithm again with the previous residual $\boldsymbol{R}_{K}$. This is called the restart or repetition. For this process, we first need an additional notation: we add a further subscript $j$ to the expansion $F_{n}$. 
Note that we have two levels of iterations here. The upper level is associated to the restart procedure and is enumerated by the second subscript $j$. The lower iteration level is the previously described ROFMP iteration with the first subscript $n$.
We denote the current expansion by  
\begin{align}
F_{n,j} \coloneqq F_{K,j-1}+\sum\limits_{i=1}^{n} \alpha_{i,j}^{(n)} d_{i,j},
\end{align}
where $F_{0,1} \coloneqq 0$ and $F_{0,j} \coloneqq F_{K,j-1}$. In analogy to the previous definitions, the residual can be defined in the following way:
\begin{align}
\boldsymbol{R}_{n,j}  \coloneqq  \boldsymbol{y}-\mathcal{F}F_{n,j}, \, 1 \leq n \leq K,\,j \geq 1 \mathrm{~and~} \boldsymbol{R}_{0,j} \coloneqq  \boldsymbol{y}-\mathcal{F}F_{K,j-1}=\boldsymbol{R}_{K,j-1}.
\end{align}
That means, after $K$ iterations, we keep the previously chosen coefficients fixed and restart the ROFMP with the residual of the step before. All in all, we have to solve
\begin{align}
\left(\alpha_{n+1,j},d_{n+1,j}\right) = \argmin_{\alpha \in \R,\, d \in \mathcal{D}} \left( \left\| \boldsymbol{R}_{n,j} -\alpha \mathcal{P}_{\mathcal{W}_{n,j}}\mathcal{F}d \right\|_{\R^l}^2 + \lambda\left\|F_{n+1,j}\right\|_{\h}^2 \right)
\end{align}
and update the coefficients in the following way
\begin{align}
\alpha_{i,j}^{(n+1)} := \alpha_{i,j}^{(n)} - \alpha_{n+1,j} \beta_{i,j}^{(n)} (d_{n+1,j}),\quad i=1,\dots,n.
\end{align}
We summarize for the expansion $F_{K,m}$, which is the approximation produced by the ROFMP after $m$ restarts:
\begin{align}
T_{m} := F_{K,m} = \sum\limits_{j=1}^{m} \sum\limits_{i=1}^{K} \alpha_{i,j}^{(K)} d_{i,j}.
\end{align}
In analogy to the RFMP, we obtain a similar convergence result for the ROFMP. That is, under certain technical conditions, the sequence $(T_m)_m$ converges in the Sobolev space $\h$. For further details, we refer to \cite{michel_telschow_pp, telschow}.

%Methdos Intro
\section{Parameter choice methods} \label{methods}

The choice of the regularization parameter $\lambda$ is crucial for the RFMP and the ROFMP, as for every other regularization method.
In this section, we briefly summarize the parameter choice methods which we test for the RFMP and the ROFMP.
This section is basically conform to \cite{bauer_gutting_lukas, bauer_lukas}.
\subsection{Introduction}

The Earth Gravitational Model 2008 (EGM2008, see \cite{egm}) is a spherical harmonics model of the gravitational potential of the Earth up to degree $2190$ and order $2159$. We use this model up to degree $100$ for the solution $F$ in our numerical tests.
For checking the parameter choice methods, we generate different test cases that means test scenarios which vary in the satellite height, the noise-to-signal ratio and the data grid. 
Based on the chosen function $F$, our dictionary contains all spherical harmonics up to degree $100$ that means our approximation $F$ from the algorithm has the following representation
\begin{align}
F = \sum\limits_{n=0}^{100} \sum\limits_{j=-n}^n \alpha_{n,j} Y_{n,j},\mathrm{~where~not~all~} \alpha_{n,j} \mathrm{~vanish}. \label{approx}
\end{align} 
This is a strong limitation, but higher degrees would essentially enlarge the computational expenses. 

Moreover, for the stabilization of the solution, we use the norm of the Sobolev space $\mathcal{H}:=\mathcal{H}((a_n)_n;\Omega)$ which is constructed with 
\begin{align}
a_{n} \coloneqq \left(n+\frac{1}{2}\right)^2, n \in \N_0, \label{sequence}
\end{align}
see \cite{freeden_gervens_schreiner}.
This Sobolev space contains all functions $F$ on $\Omega$ which fulfil
\begin{align}
\sum_{n=0}^\infty \sum_{j=-n}^n a_n^2 \left\langle F,Y_{n,j} \right\rangle_{\LO}^2 < \infty. 
\end{align}
The inner product of functions $F, G \in \mathcal{H}$ is given by
\begin{align}
\left\langle  F,G \right\rangle_{\h} \coloneqq \sum_{n=0}^{\infty} \sum_{j=-n}^n a_n^2 \left\langle F,Y_{n,j} \right\rangle_{\LO} \left\langle G,Y_{n,j} \right\rangle_{\LO}. 
\end{align}

The particular sequence $(a_n)_n = ((n+\frac{1}{2})^2)$ was chosen, because preliminary numerical experiments showed that the associated regularization term yielded results with an appropriate smoothness.

In our test scenarios, we use a finite set $\{\lambda_k\}_{k=1,\dots,100}$ of $100$ regularization parameters (for details, see \cref{lambda}). 
The approximate solution of the inverse problem as an output of the RFMP/ROFMP corresponding to the regularization parameter $\lambda_k$ and the data vector $\boldsymbol{y}$ is denoted by $x_k$. This notation is introduced to avoid confusions with the functions $F_n$ which occur at the $n$-th step of the iteration within the RFMP. 

In practice, we deal with noisy data $\boldsymbol{y}^\eps$ where the noise level $\eps$ is defined by
\begin{align}
\eps \coloneqq \mathrm{N2S} \cdot \left\|\boldsymbol{y}\right\|_{\R^l} / \sqrt{l}, \label{eps}
\end{align}
where $l$ is the length of the data vector $\boldsymbol{y}$ and N2S is called the noise-to-signal ratio.
The corresponding result of the RFMP/ROFMP for the regularization parameter $\lambda_k$ and the noisy data vector $\boldsymbol{y}^\eps$ is called $x_k^\eps$. Due to the convergence results for the RFMP/ROFMP (see \eqref{normalequ}), we introduce the linear regularization operators $\mathcal{R}_k:\R^l \rightarrow \h$, 
\begin{align}
\mathcal{R}_k \coloneqq (\mathcal{F}^* \mathcal{F}+\lambda_k \mathcal{I})^{-1} \mathcal{F}^*
\end{align} 
and assume $x_k$ to be $\mathcal{R}_k \boldsymbol{y}$ and $x_k^\eps$ to be $\mathcal{R}_k \boldsymbol{y^\eps}$, though this could certainly only be guaranteed for an infinite number of iterations.

Due to the importance of the regularization parameter, we summarize in the next section some methods for the choice of this parameter $\lambda$.
For the comparison of the methods, we have to define the optimal regularization parameter $\lambda_{k_{\mathrm{opt}}}$. We do this by minimizing the difference between the exact solution $x$ and the regularized solution $x_k^{\eps}$ corresponding to the parameter $\lambda_k$ and noisy data.
\begin{align}
k_{\mathrm{opt}} \coloneqq \argmin_{k \in \{1,\dots,100\}} \| x-x_k^\eps \| _{\LO}.
\end{align}
Then we evaluate the results by computing the so-called \emph{inefficiency} by
\begin{align}
\frac{\left\|x-x_{k_{*}}^\eps \right\|_{\LO}}{\left\|x-x_{k_{\mathrm{opt}}}^\eps\right\|_{\LO}}, \label{inefficiency}
\end{align}
where $\lambda_{k_*}$ is the regularization parameter selected by the considered parameter choice method.
For the computation of the inefficiency, we use the $\LO$-norm, since our numerical results led to a better distinction of the different inefficiencies than by using the $\mathcal{H}$-norm. However, the tendency regarding 'good' and 'bad' parameters were the same in both cases.
The closer the obtained inefficiency is to 1, the better the parameter choice method performs.

%Parameter choice methods

The norms which occur in the several parameter choice methods can be computed with the help of the singular value decomposition. However, we use the singular values of $\Psi$ (see \eqref{psi} and \eqref{task_points}) for this purpose, because the singular value decomposition of $\mathcal{F}$ is unavailable. This certainly causes an inaccuracy in our calculations, but appears to be unavoidable for the sake of practicability.

\subsection{Parameter Choice Methods}
\cref{tab:methods} shows the different parameter choice methods we tested. The tuning parameters are chosen in accordance to \cite{bauer_gutting_lukas, bauer_lukas}. For the choice of the maximal index $\hat{K}$, see \cref{max_lambda}.
\begin{sidewaystable}
\begin{center}
\tabulinesep=1mm
\begin{tabu}[H]{c|c|c}
Name & Selection criterion & Specifications \\ 
\hline
Discrepancy Principle (DP) & Choose the first $k$ such that & Tuning parameter $\tau >1$. \\
(References: \cite{morozov1966, morozov1984,phillips}) & $\| \mathcal{F} x_k^\eps-\boldsymbol{y}^\eps \|_{\R^l} \leq \tau \eps \sqrt{l}$. & (We choose $\tau =1.5$.)  \\
\hline
Transformed Discrepancy Principle (TDP) & Choose the first $k$ such that & Tuning parameter $b>\gamma = ((1/4)^{(1/4)}(3/4)^{(3/4)})^2$, \\
(References: \cite{raus1990, raus1992}) & $\| \mathcal{R}_{k} (\mathcal{F} x_k^\eps -\boldsymbol{y}^\eps)  \|_{\mathcal{H}} \leq b \hat{\eps} \frac{\sqrt{l}}{\sqrt{\lambda_k}}$. & estimate $\hat{\eps}$ of $\eps$.  (We choose $b=1.5\gamma$ and $\hat{\eps}=\eps$.)  \\
\hline
Quasi-optimality Criterion (QOC) & $k_{*}=\argmin\limits_{k \leq \hat{K}} \| x_k^\eps-x_{k+1}^\eps \|_{\mathcal{H}}$ &  \\
(References: \cite{tikhonov_arsenin, tikhonov_glasko}) & & \\
\hline
L-curve Method (LC) & $k_{*}=\argmin\limits_{k\leq \hat{K}} \{ \| \mathcal{F} x_k^\eps -\boldsymbol{y}^\eps \|_{\R^l} \cdot \| x_k^\eps \|_{\mathcal{H}} \}$ & \\
(References: \cite{hansen1992, hansen1998, hansen_oleary}) & & \\
\hline
Extrapolated Error Method (EEM) & $k_{*} = \argmin\limits_{k \leq \hat{K}} \Big\{ \frac{\| \mathcal{F}x_k^\eps-\boldsymbol{y}^{\eps}  \|_{\R^l}^2}{\| \mathcal{F}^* (\mathcal{F}x_k^\eps-\boldsymbol{y}^{\eps}) \|_{\mathcal{H}}} \Big\}$ & \\
(References: \cite{brezinski2008, brezinski2009}) & & \\
\hline 
Residual Method (RM) & $k_{*}=\argmin\limits_{k \leq \hat{K}} \Big\{ \frac{\| \mathcal{F} x_k^\eps -\boldsymbol{y}^\eps \|_{\R^l}}{(\mathrm{tr} \mathcal{B}^* \mathcal{B})^{1/4}} \Big\} $, & \\
(References: \cite{bauer_mathe}) & where $\mathcal{B}=\mathcal{F}(\mathcal{I}-\mathcal{R}_{k}\mathcal{F})$.  & \\
\hline 
Generalized Maximum Likelihood (GML) & $k_* =\argmin\limits_{k \leq \hat{K}} \Big\{ \frac{\| \mathcal{F} x_k^\eps -\boldsymbol{y}^\eps \|_{\R^l}^2}{(\mathrm{det}^+(\mathcal{I}-\mathcal{F}\mathcal{R}_{k}))^{1/l_1}} \Big\}$ & $l_1=\mathrm{rank} (\mathcal{I}-\mathcal{F}\mathcal{R}_{k})$. (In our case $l_1=l$.)  \\
(References: \cite{wahba1985}) &  & $\mathrm{det}^+$ is the product of the nonzero eigenvalues. \\
\hline 
Generalized Cross Validation (GCV) & $k_{*}=\argmin\limits_{k \leq \hat{K}} \Big\{ \frac{\left\| \mathcal{F} x_k^\eps-\boldsymbol{y}^{\eps} \right\|_{\R^l}^2}{(l^{-1}\mathrm{tr}(\mathcal{I}-\mathcal{F}\mathcal{R}_{k}))^2} \Big\}$ &  \\
(References: \cite{wahba1977}) & & \\
\hline
Robust GCV (RGCV) & $k_{*}=\argmin\limits_{k \leq \hat{K}} \Big\{ \frac{\left\| \mathcal{F} x_k^\eps-\boldsymbol{y}^{\eps} \right\|_{\R^l}^2}{(l^{-1}\mathrm{tr}(\mathcal{I}-\mathcal{F}\mathcal{R}_{k}))^2} $ & Robustness parameter $\gamma \in (0,1)$. \\
(References: \cite{lukas2006, robinson_moyeed}) & $ \times \left( \gamma + (1-\gamma)l^{-1} \mathrm{tr}((\mathcal{F}\mathcal{R}_{k})^2) \right) \Big\}$ & (We choose $\gamma=0.1$.) \\
\hline 
Strong RGCV (SRGCV) & $k_{*}=\argmin\limits_{k \leq \hat{K}} \Big\{ \frac{\left\| \mathcal{F} x_k^\eps-\boldsymbol{y}^{\eps} \right\|_{\R^l}^2}{(l^{-1}\mathrm{tr}(\mathcal{I}-\mathcal{F}\mathcal{R}_{k}))^2} $ & Robustness parameter $\gamma \in (0,1)$. \\
(References: \cite{lukas2008}) & $ \times \left( \gamma + (1-\gamma)l^{-1} \mathrm{tr}((\mathcal{F}\mathcal{R}_{k})^2) \right) \Big\}$  & (We choose $\gamma=0.95$.) \\
\hline 
Modified Generalized Cross Validation (MGCV) &  $k_{*}=\argmin\limits_{k \leq \hat{K}} \Big\{ \frac{\left\| \mathcal{F} x_k^\eps-\boldsymbol{y}^{\eps} \right\|_{\R^l}^2}{(l^{-1}\mathrm{tr}(\mathcal{I}-c\mathcal{F}\mathcal{R}_{k}))^2} \Big\}$ & Stabilization parameter $c>1$. \\
(References: \cite{cummins, vio}) & & (We choose $c=3$.)
\end{tabu}
\caption{The parameter choice methods and their specifications.}
\label{tab:methods}
\end{center}
\end{sidewaystable}

%Specifications
\section{Evaluation} \label{spec_methods}
\subsection{Specifications for the algorithm} \label{sec41}
In \cref{rfmp_algorithm,rofmp_algorithm}, we mentioned that we need to define stopping criteria for our algorithm.
We state the following stopping criteria for the RFMP and ROFMP (see also \cref{rfmp_algo,rofmp_algo}).
\begin{itemize}
\item $\left\|\boldsymbol{R}_{n+1}\right\|_{\R^l} <\varrho$ for a given $\varrho> 0$ (in our case, this is the N2S),
\item $n+1 > N$ for a given $N \in \N$ (in our case, $N=10000$ because of our computing capacity),
\item $\alpha_{n+1} < \delta$ for a given $\delta>0$ (in our case $\delta=10^{-6}$). 
\end{itemize}
In the case of the ROFMP, we choose $K=200$ for the restart.

\subsection{The data grids}
\begin{figure}[!ht]
\begin{minipage}[t]{0.49\textwidth}
\includegraphics[trim= 30mm 100mm 15mm 114mm ,clip, width=\textwidth]{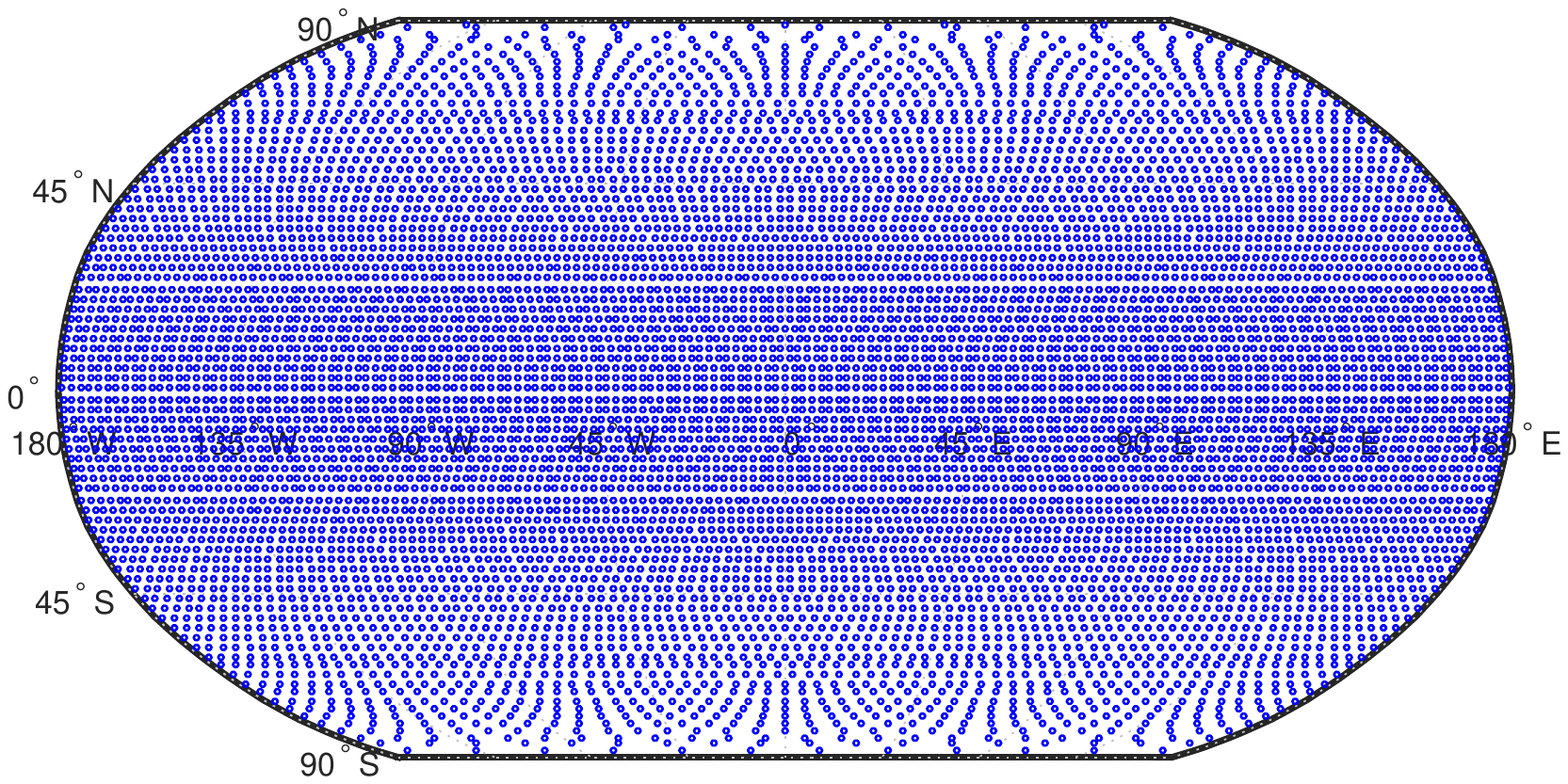}
\end{minipage}
\begin{minipage}[t]{0.49\textwidth}
\includegraphics[trim= 30mm 100mm 15mm 114mm ,clip, width=\textwidth]{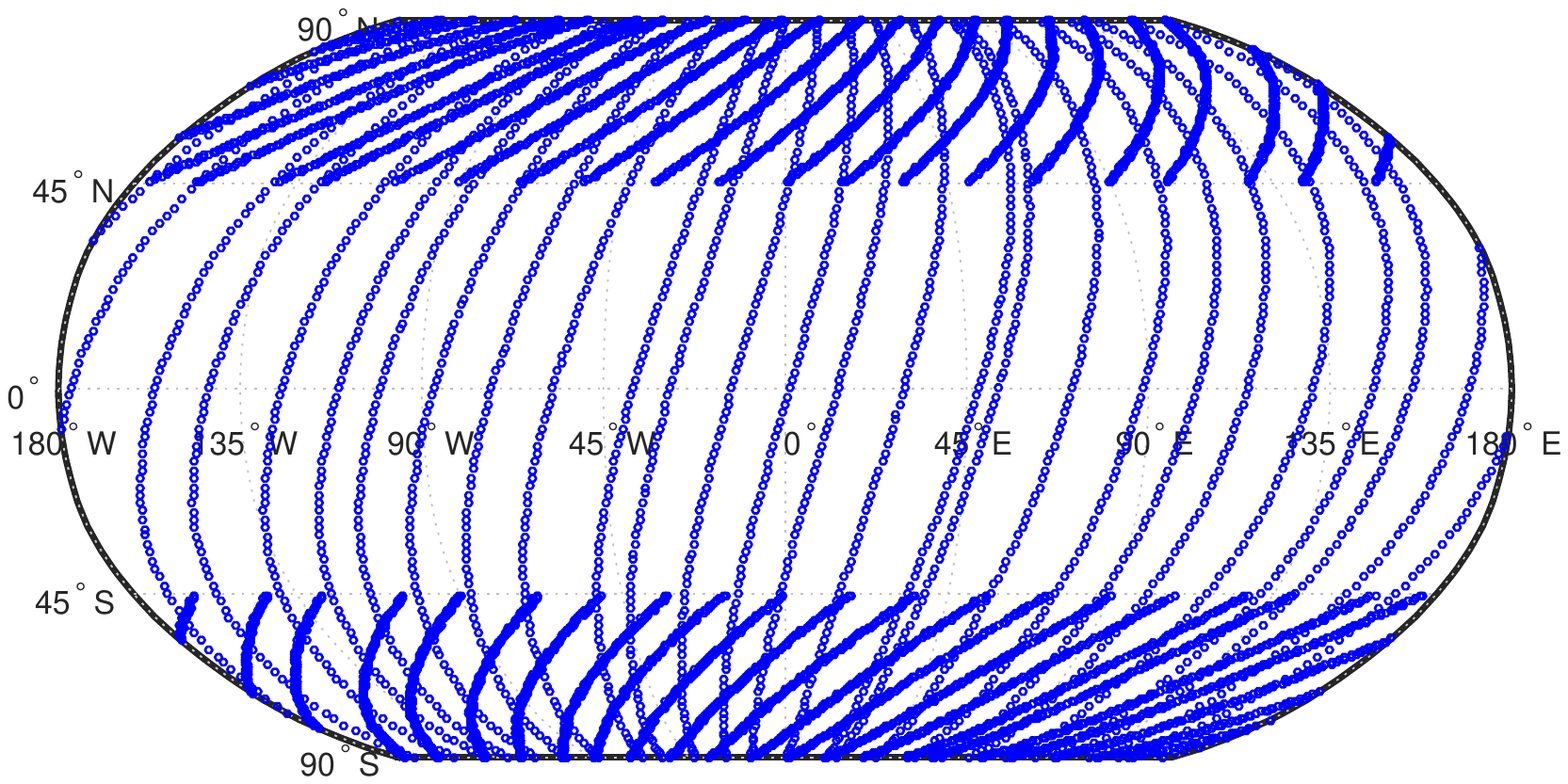}
\end{minipage}
\caption{Reuter grid with 8514 points (left) and scattered grid with 8500 points (right).}
\label{fig:grids}
\end{figure}

\cref{fig:grids} shows two data grids which we use for our experiments. First of all, the Reuter grid (see \cite{reuter}) is an example of a regular data grid on the sphere. Second, we have a set of irregularly distributed data points on a grid which we refer to as the scattered grid in the following and which was first used in \cite{telschow}. The latter grid tries to imitate the distribution of measurements along the tracks of a satellite. It possesses additional shorter tracks and, thus, a higher accumulation of data points at the poles and only fewer tracks in a belt around the equator.

\subsection{Noise generation}
For our various scenarios, we get our noisy data if we add white noise to our data values or we add coloured noise that is obtained by an autoregression process. Additionally, we test some local noise.

\subsubsection{White noise scenario}
For white noise, we add Gaussian noise corresponding to a certain noise-to-singal ratio N2S to the particular value of each datum, that means we get our noisy data by
\begin{align}
y_{i}^{\eps}=(1+\mathrm{N2S} \cdot  \epsilon_{i})y_{i}=\mathrm{N2S} \cdot y_{i} \epsilon_{i}+y_{i},\quad i=1,\dots,l, \label{noise_gen}
\end{align}
where $y_{i}$ are the components of $\boldsymbol{y}$ and $\epsilon_{i} \sim \mathcal{N}(0,1)$, that means every $\epsilon_{i}$ is a standard normally distributed random variable. 

\subsubsection{Coloured noise scenario}
Since our scattered grid tries to imitate tracks of satellites, we can assume that we have a chronology of the data points for each track.
To obtain some sort of coloured noise, we use an autoregression process of order $1$ (briefly: AR(1)-process, see \cite{brockwell_davis}) with whom we simulate correlated noise.

A stochastic process $\{\epsilon_{i},\, i \in \mathbb{Z}\}$ is called an autoregressive process of order $1$, if $\epsilon_{i}=\alpha \epsilon_{i-1}+\eps_{i},\  \left|\alpha\right| < 1$, where $\eps_{i} \sim \mathcal{N}(0,1)$. In the case of our simulation, we start with $\epsilon_{1} \sim \mathcal{N}(0,1)$ and run the recursion for a fixed $\alpha \in (-1,1)$, which we determined at random.
%Moreover, we omit $\epsilon_{1},\dots,\epsilon_{200}$ to obtain a stationary process.

\begin{figure}[!ht]
\begin{minipage}[t]{0.49\textwidth}
\includegraphics[trim= 18mm 70mm 14mm 90mm ,clip, width=\textwidth]{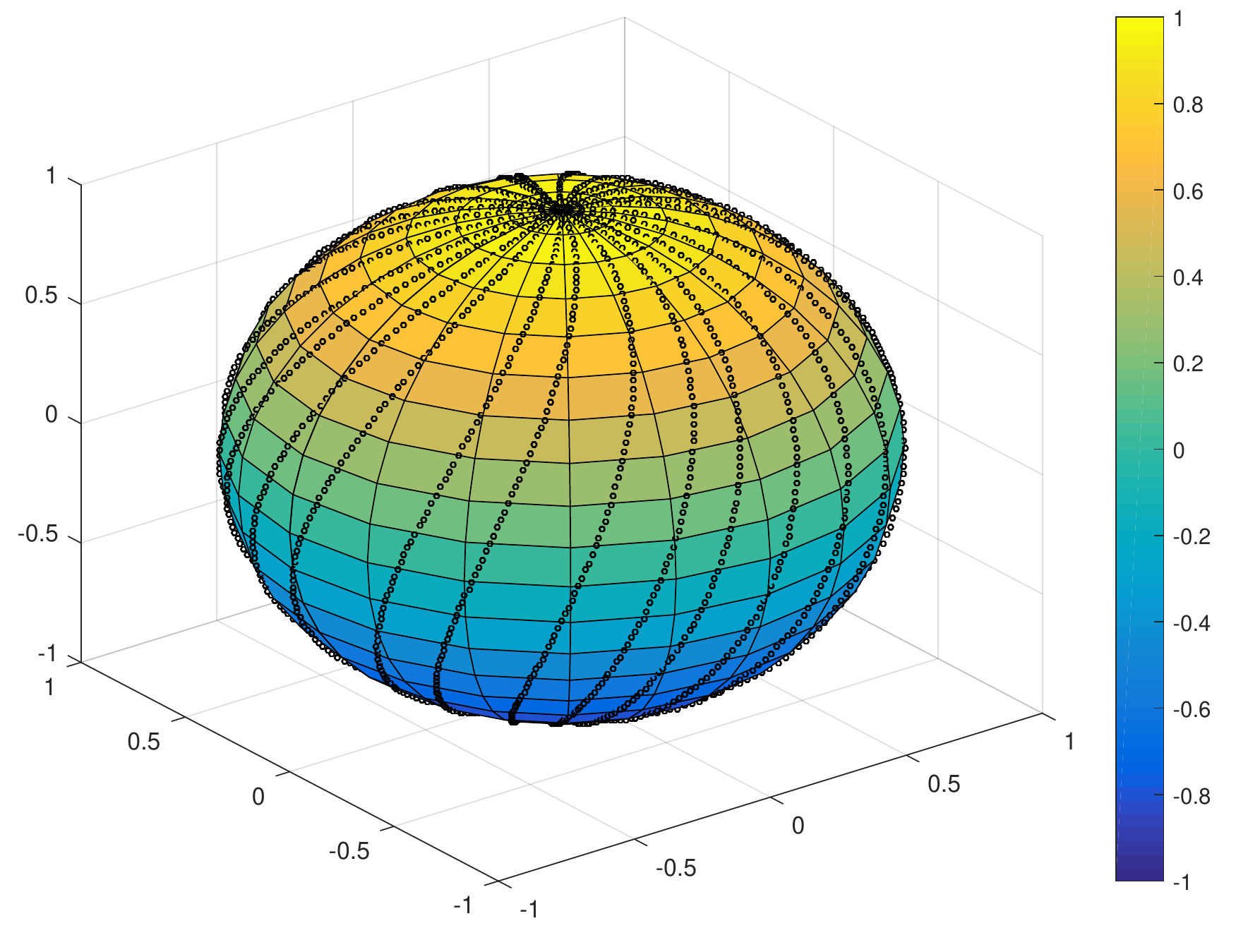}
%\fbox{\includegraphics[trim= 18mm 70mm 14mm 90mm ,width=\textwidth]{Figures/Punkte_sphere/scattersphere1.pdf}}
%\caption{One track set of the scattered grid.}
%\label{fig:scatteredsphere1}
\end{minipage}
%\hspace*{0.2cm}
\begin{minipage}[t]{0.49\textwidth}
\includegraphics[trim= 18mm 70mm 14mm 90mm ,clip, width=\textwidth]{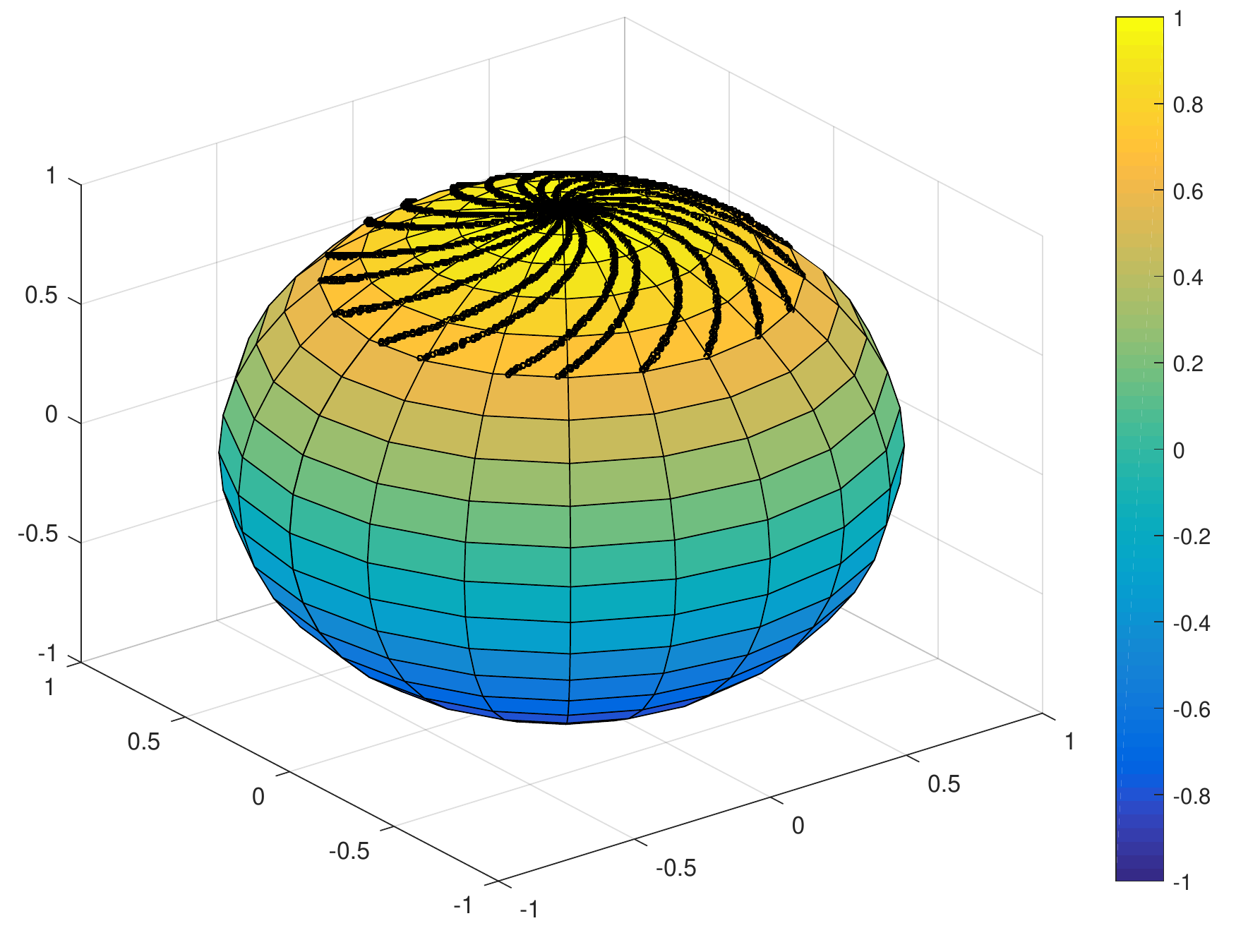}
%\fbox{\includegraphics[trim= 18mm 70mm 14mm 90mm ,width=\textwidth]{Figures/Punkte_sphere/scattersphere2.pdf}}
%\caption{Track set of the scattered grid at the north pole (for the south pole we have the same point distribution.}
%\label{fig:scatteredsphere2}
\end{minipage}
\caption{The track sets of the scattered grid (left and right). For the South pole, we have an analogous point distribution.}
\label{fig:tracks}
\end{figure}

For each track of the scattered grid, we apply this autoregression process (for the tracks, see \cref{fig:tracks}) and obtain $y_i^{\eps}$ as in \eqref{noise_gen} using the $\epsilon_i$ from above.

\subsubsection{Local noise scenario}
\begin{figure}[!ht]
\begin{minipage}[t]{0.45\textwidth}
\includegraphics[trim= 30mm 100mm 15mm 114mm ,clip, width=\textwidth]{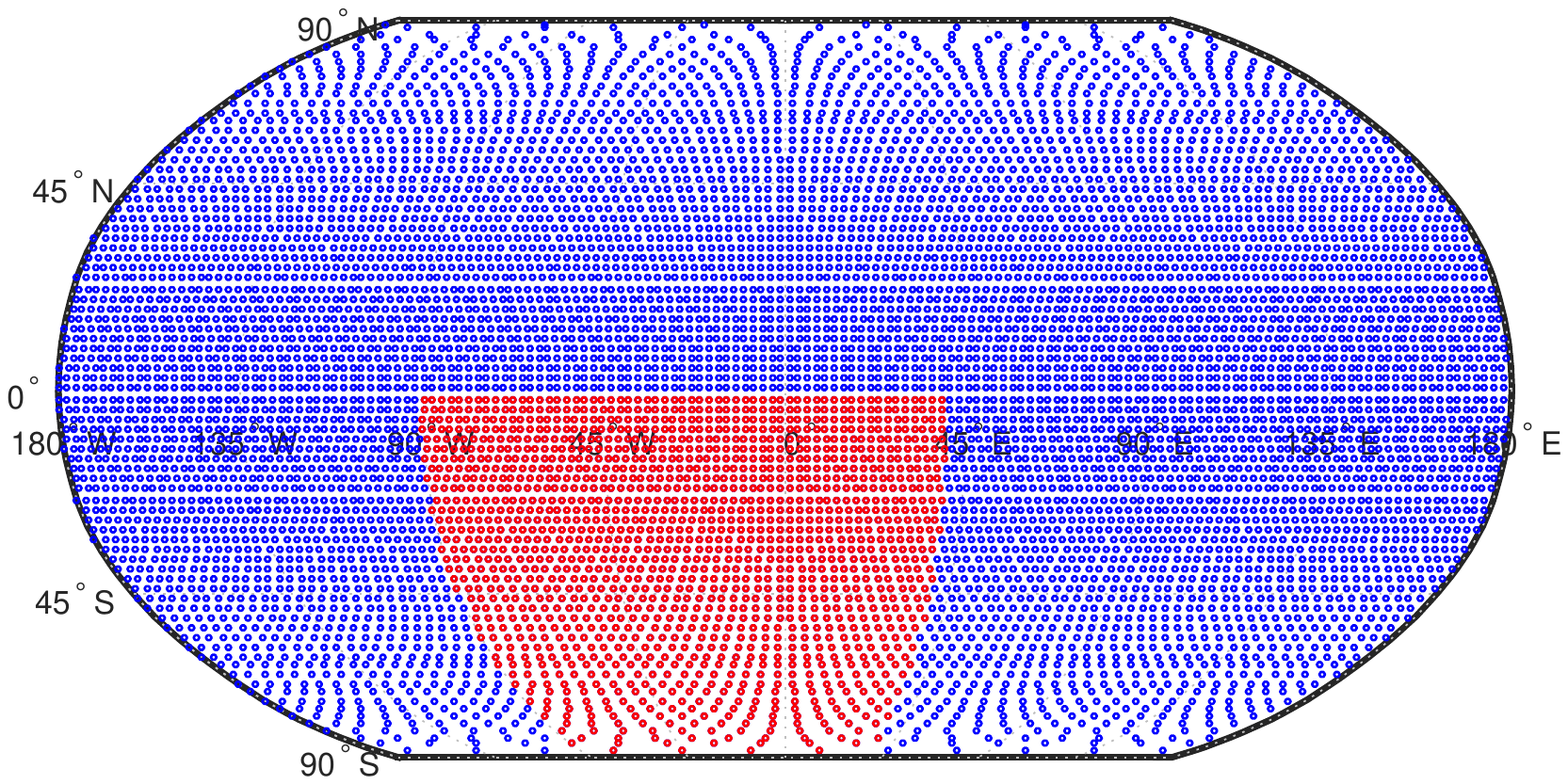}
\end{minipage}
\hspace*{0.2cm}
\begin{minipage}[t]{0.45\textwidth}
\includegraphics[trim= 30mm 100mm 15mm 114mm ,clip, width=\textwidth]{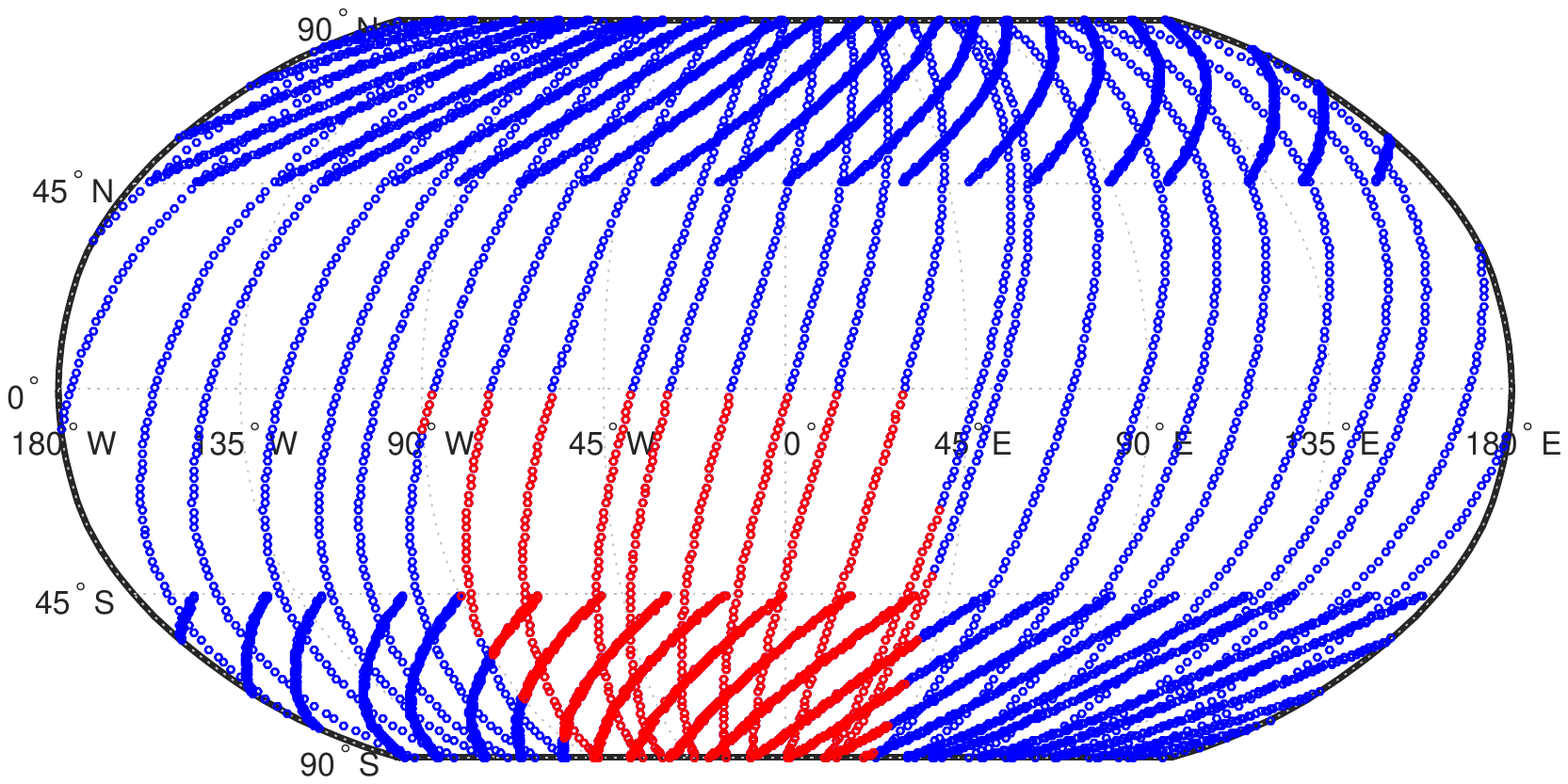}
\end{minipage}
\caption{Reuter grid (left) and scattered grid (right). The values of the data points in the red area contain an N2S of $5\%$ and the values of the data points in the blue area an N2S of $1\%$ for the local noise scenario.}
\label{fig:localnoise}
\end{figure}

For the local noise, we choose a certain area and add white noise with an $\mathrm{N2S}=5\%$ relative to the particular value to each data point. To the values of the remaining data points we add white noise with an $\mathrm{N2S}=1\%$.
We choose this area as illustrated in \cref{fig:localnoise}. The choice of this area is a very rough approximation of the domain of the South Atlantic Anomaly, where a dip in the Earth's magnetic field exists (see e.g. \cite{saa_heirtzler}). Since only a few points of our grid would have been in the actual domain, we extended the area towards the South pole.

\cref{tab:testcases} shows our different test cases for the RFMP and ROFMP.
\begin{table}[H]
\begin{center}
\begin{tabular}{c|c|c|c|c}
height & N2S & noise & grid & shortcut\\
\hline 
500km & 5\% & white & scattered & (500,5,wn,S) \\
500km & 5\% & coloured & scattered & (500,5,cn,S) \\
500km & 5\% & white & Reuter & (500,5,wn,R) \\
\hline
500km & 1\% & white & scattered & (500,1,wn,S) \\
500km & 1\% & coloured & scattered & (500,1,cn,S) \\
500km & 1\% & white & Reuter & (500,1,wn,R) \\
\hline
300km & 5\% & white & scattered & (300,5,wn,S) \\
300km & 5\% & coloured & scattered & (300,5,cn,S)\\
300km & 5\% & white & Reuter & (300,5,wn,R) \\
\hline
500km & 5\%/1\% & local & scattered & (500,5,ln,S) \\
500km & 5\%/1\% & local & Reuter & (500,5,ln,R) \\
\end{tabular}
\caption{Overview of the implemented test cases.}
\label{tab:testcases}
\end{center}
\end{table}

\subsection{Regularization parameters} \label{lambda}
We constructed the admissible values $\lambda_k$ for the parameter choice as a monotonically decreasing sequence with 100 values from $\lambda_1=1$ to $\lambda_{100}=10^{-14}$ and a logarithmically equal spacing in the following way
\begin{align}
\lambda_{k}=\lambda_{0}q_{\lambda}^k,\quad k\in \{1,\dots,100\}. \label{lambda2}
\end{align}
Here, $\lambda_{0}=1.3849$ and $q_{\lambda}=0.7221$. The test scenarios are chosen such that the parameter range lies between $1$ and $10^{-14}$ and includes the optimal parameter away from the boundaries $\lambda_1$ and $\lambda_{100}$. For the choice of the parameters $\lambda_{k}$, we refer to \cite{bauer_gutting_lukas, bauer_lukas}. 

\subsection{Maximal index} \label{max_lambda}

Most parameter choice methods either increase the index $k$ until a certain condition is satisfied or minimize a certain function for all regularization parameters $\lambda$, i.e. after our discretization (see \eqref{lambda2}) they minimize for all $k$  (see \cref{tab:methods}).
For some methods like the quasi-optimality criterion, the values of $k$ have to be constrained by a suitable maximal index $\hat{K}$ which must be chosen such that $k_{\mathrm{opt}}<\hat{K}$.
To increase computational efficiency, such a maximal index can be used for other methods as well without changing their performance. As in \cite{bauer_gutting_lukas, bauer_lukas}, we define this maximal index by 
\begin{align}
\hat{K}=\max\left\{ k \, \big| \, \rho(k) < 0.5\rho(\infty) \right\},
\end{align}
where $\mathbb{E}\|x_k-x_{k}^{\eps}\|^2=\eps^2 \rho^2(k)$ is the variance of the regularized solution corresponding to noisy data and $\eps^2 \rho^2(\infty)$ is its largest value.  It is well-known that, in the case of white noise, $\rho (k)$ for the Tikhonov-Phillips regularization is generally given by
\begin{align}
\rho^2 (k)=\sum_{n} \left( \dfrac{\sigma_n}{\sigma_n^2+\lambda_k} \right)^2. 
\end{align}
Since our singular values occur with a multiplicity of $2n+1$ and we restrict our tests to $n=0,\dots,100$, the sum above in our tests is given by
\begin{align}
\rho^2 (k)=\sum_{n=0}^{100}  (2n+1) \left( \dfrac{\sigma_n}{\sigma_n^2+\lambda_k} \right)^2.
\end{align}
For any coloured noise, we use the estimate (cf. \cite{bauer_lukas})
\begin{align}
\eps^2 \rho^2 (k) \approx 2^{-1}  \left\| x_{k,1}^{\eps}-x_{k,2}^{\eps} \right\|_{\h}^2,
\end{align}
with two independent data sets $\boldsymbol{y}_1^\eps$, $\boldsymbol{y}_2^\eps$ for the same regularization parameter $\lambda_k$. Note that $x_{k,1}^\eps$, $x_{k,2}^\eps$ are the regularized solutions corresponding to the parameter $\lambda_k$ and the noisy data sets $\boldsymbol{y}_1^\eps$, $\boldsymbol{y}_2^\eps$.

%Comparison
\section{Comparison of the methods} \label{comparison}
For the error comparison, we compute the inefficiency (see \eqref{inefficiency}) in each scenario (see \cref{tab:testcases} for an overview) for each parameter choice method and compare the inefficiencies. We generate 32 data sets for each of the eleven scenarios, i.e.  we run each algorithm for 352 times for a single regularization parameter. \crefrange{fig:dp}{fig:mgcv} show the inefficiencies, collected based on the parameter choice methods. The red middle band in the box is the median and the red $+$ symbol shows outliers. The boxplots of our results are plotted at a logarithmic scale.

\subsection{Discrepancy principle (DP)}
\begin{figure}[!ht]
\begin{minipage}[t]{0.49\textwidth}
\includegraphics[trim= 20mm 65mm 19mm 92mm ,clip, width=\textwidth]{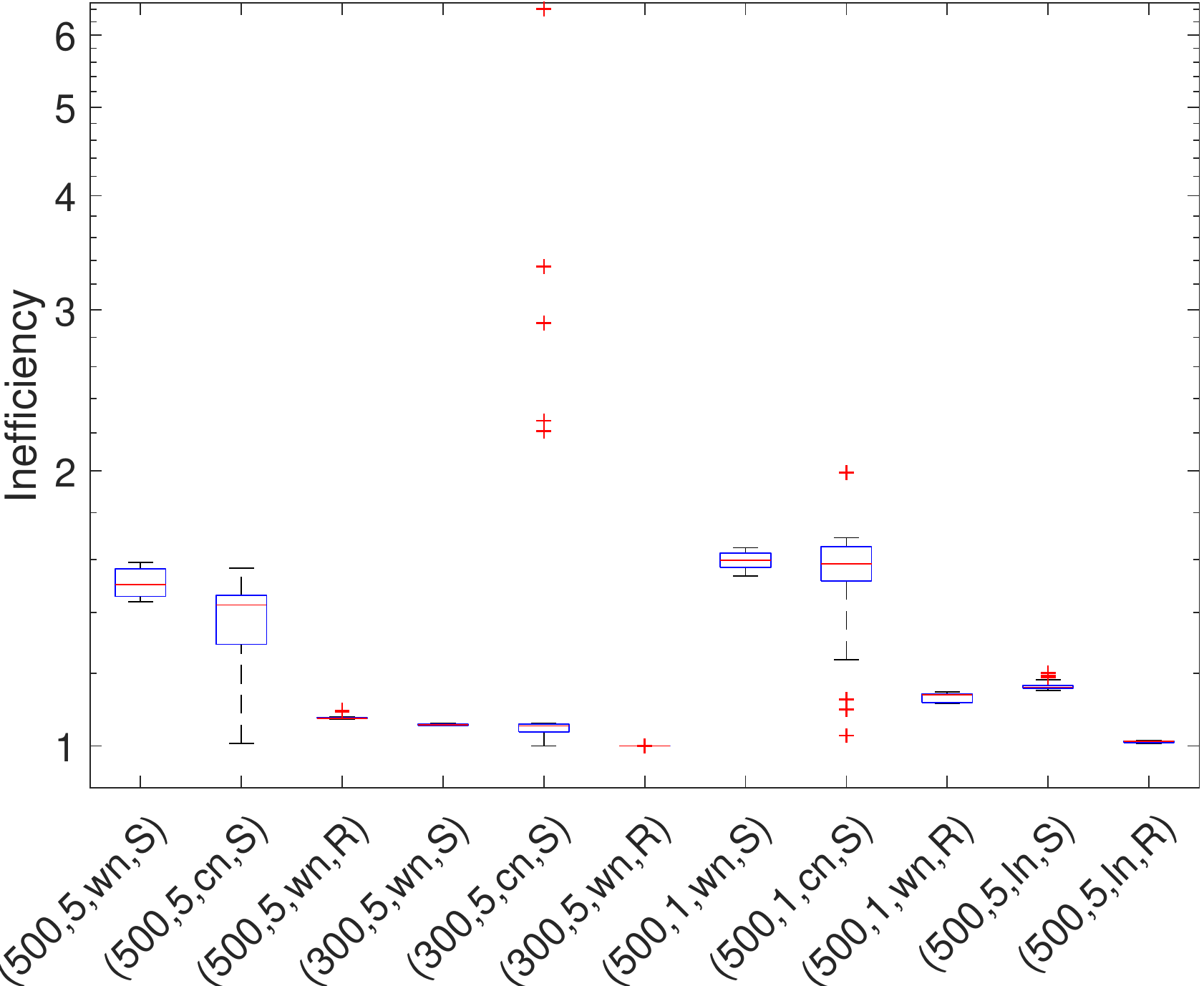}
\end{minipage}
\hspace*{0.1cm}
\begin{minipage}[t]{0.49\textwidth}
\includegraphics[trim= 20mm 65mm 19mm 92mm ,clip, width=\textwidth]{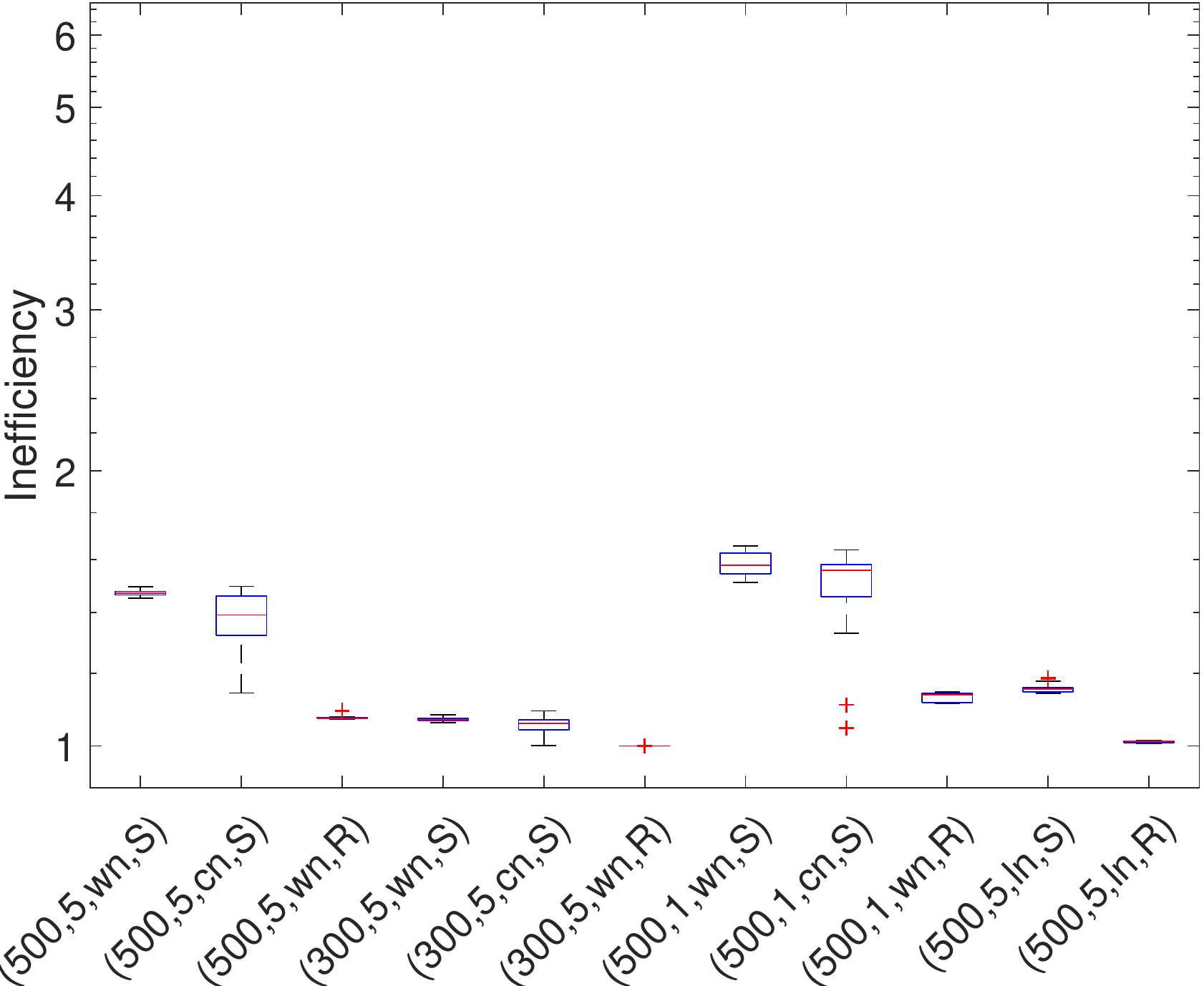}
\end{minipage}
\caption{DP for the RFMP (left) and the ROFMP (right).}
\label{fig:dp}
\end{figure}
We can see from \cref{fig:dp} that the DP leads to results which are in the range  from good to acceptable in all test cases. It yields better results with a more uniformly distributed grid.

\subsection{Transformed discrepancy principle (TDP)}
\begin{figure}[!ht]
\begin{minipage}[t]{0.49\textwidth}
\includegraphics[trim= 20mm 65mm 19mm 92mm ,clip, width=\textwidth]{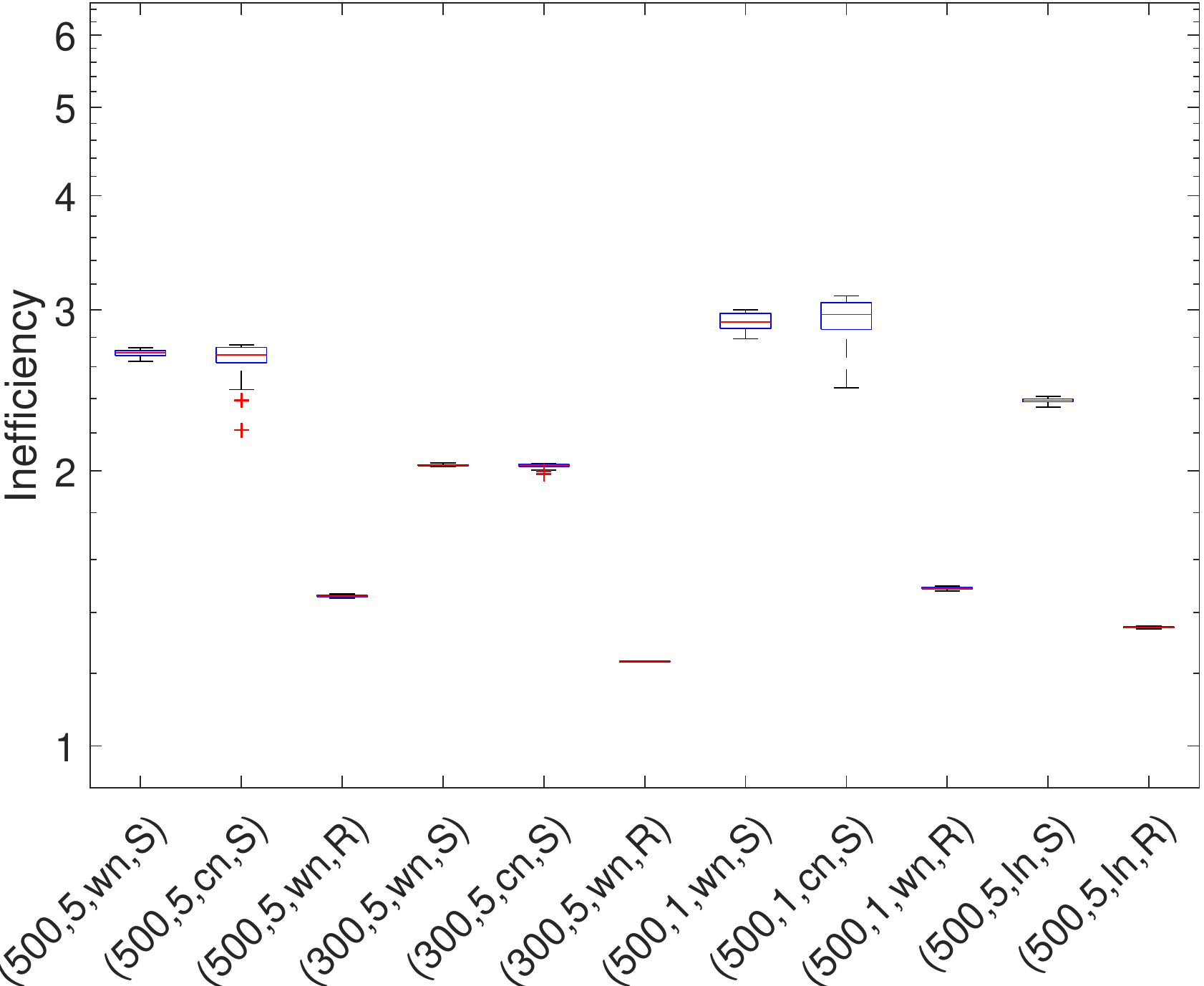}
\end{minipage}
\hspace*{0.1cm}
\begin{minipage}[t]{0.49\textwidth}
\includegraphics[trim= 20mm 65mm 19mm 92mm ,clip, width=\textwidth]{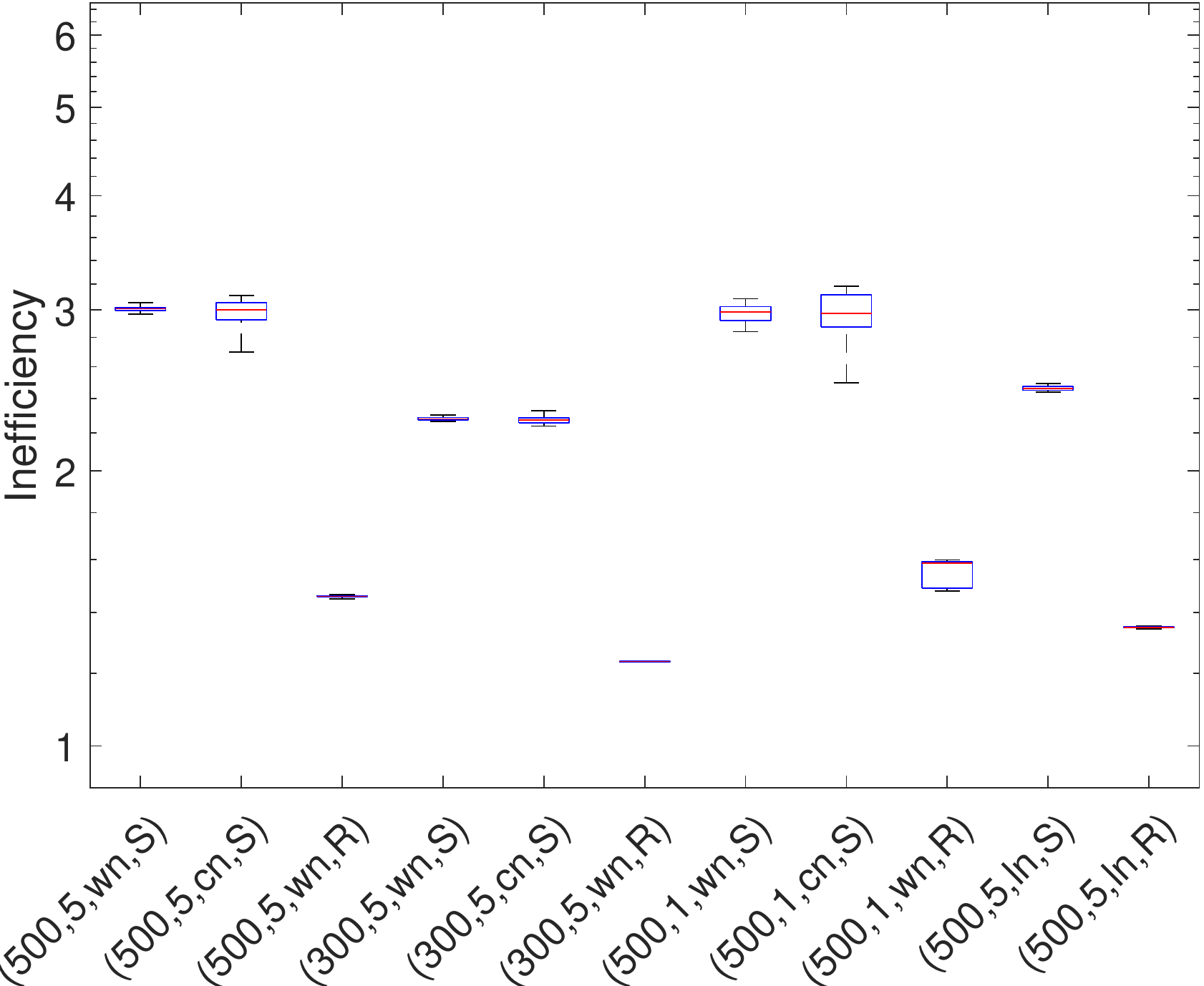}
\end{minipage}
\caption{TDP for the RFMP (left) and the ROFMP (right).}
\label{fig:tdp}
\end{figure}
The results for the TDP (see \cref{fig:tdp}) for all test cases are rather poor. We can remark that the results get better with a more uniformly distributed data grid. Furthermore, the coloured noise leads to slightly bigger boxes than the white noise.

\subsection{Quasi-optimality criterion (QOC)}
\begin{figure}[!ht]
\begin{minipage}[t]{0.49\textwidth}
\includegraphics[trim= 20mm 65mm 19mm 92mm ,clip, width=\textwidth]{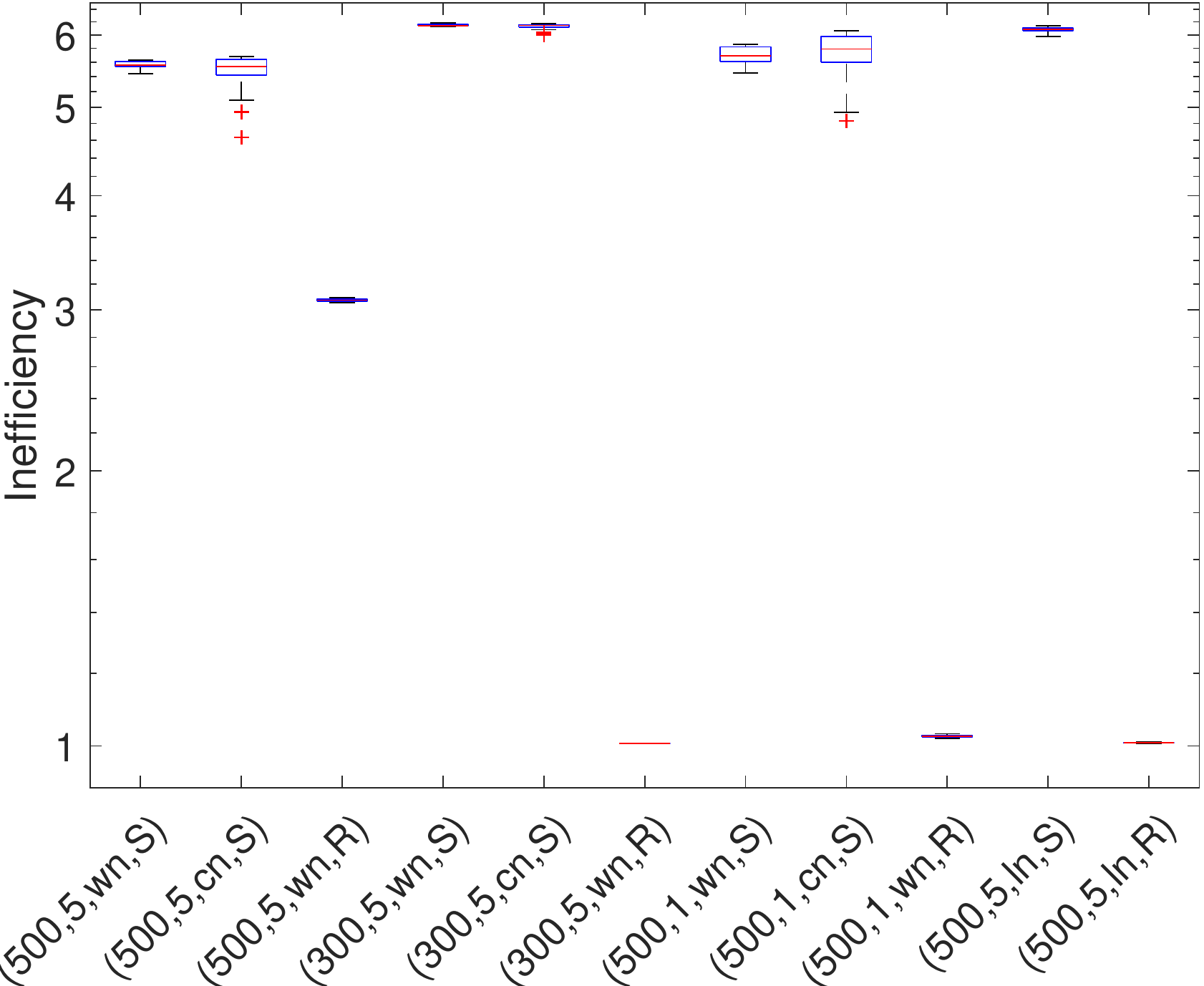}
\end{minipage}
\hspace*{0.1cm}
\begin{minipage}[t]{0.49\textwidth}
\includegraphics[trim= 20mm 65mm 19mm 92mm ,clip, width=\textwidth]{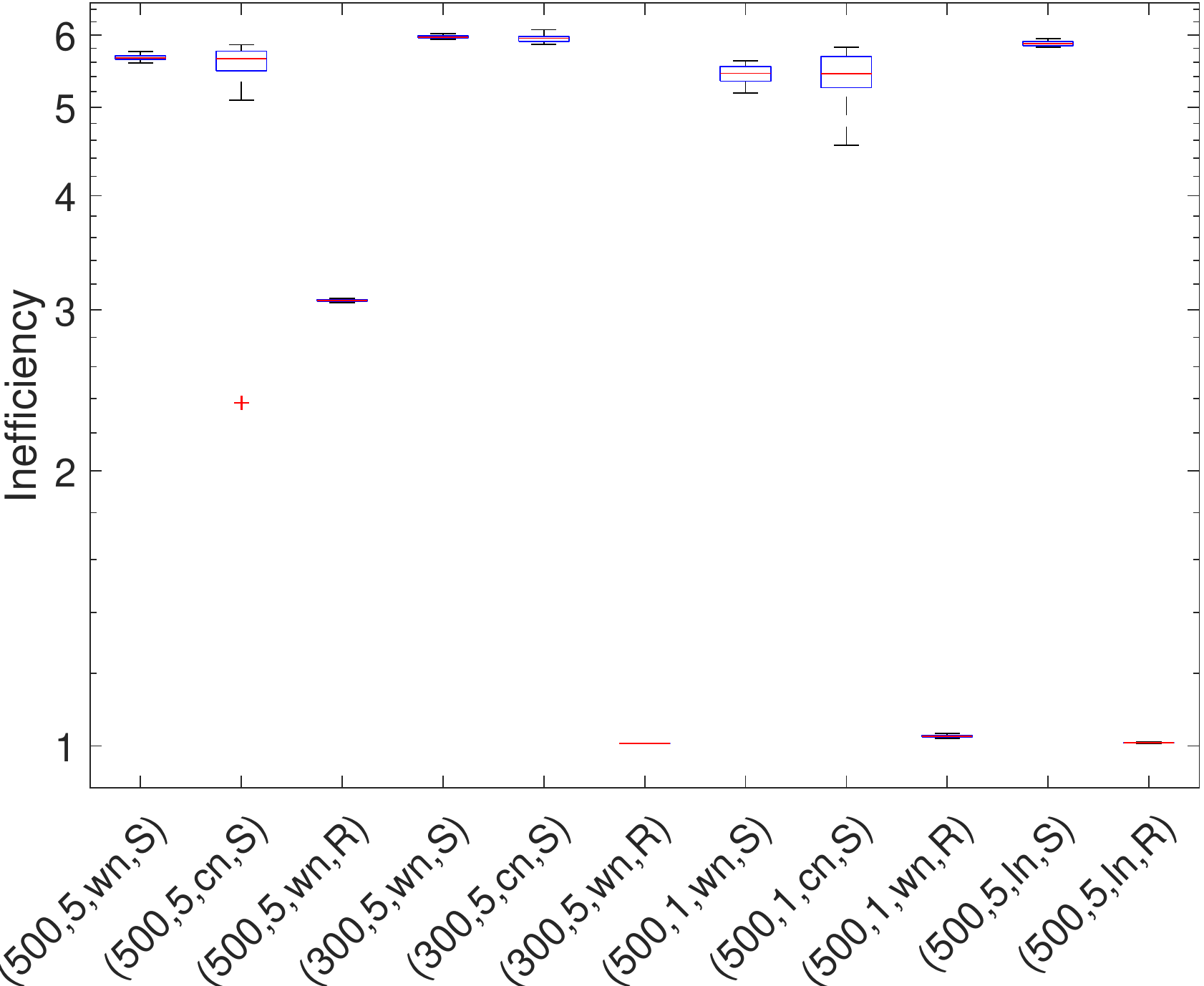}
\end{minipage}
\caption{QOC for the RFMP (left) and the ROFMP (right).}
\label{fig:qoc}
\end{figure}
In \cref{fig:qoc}, the inefficiencies of the QOC show that the performance of this method is rather poor. In the case of the Reuter grid, the results reach from good to mediocre in contrast to the scattered grid.

\subsection{L-curve method (LC)}
\begin{figure}[!ht]
\begin{minipage}[t]{0.49\textwidth}
\includegraphics[trim= 20mm 65mm 19mm 92mm ,clip, width=\textwidth]{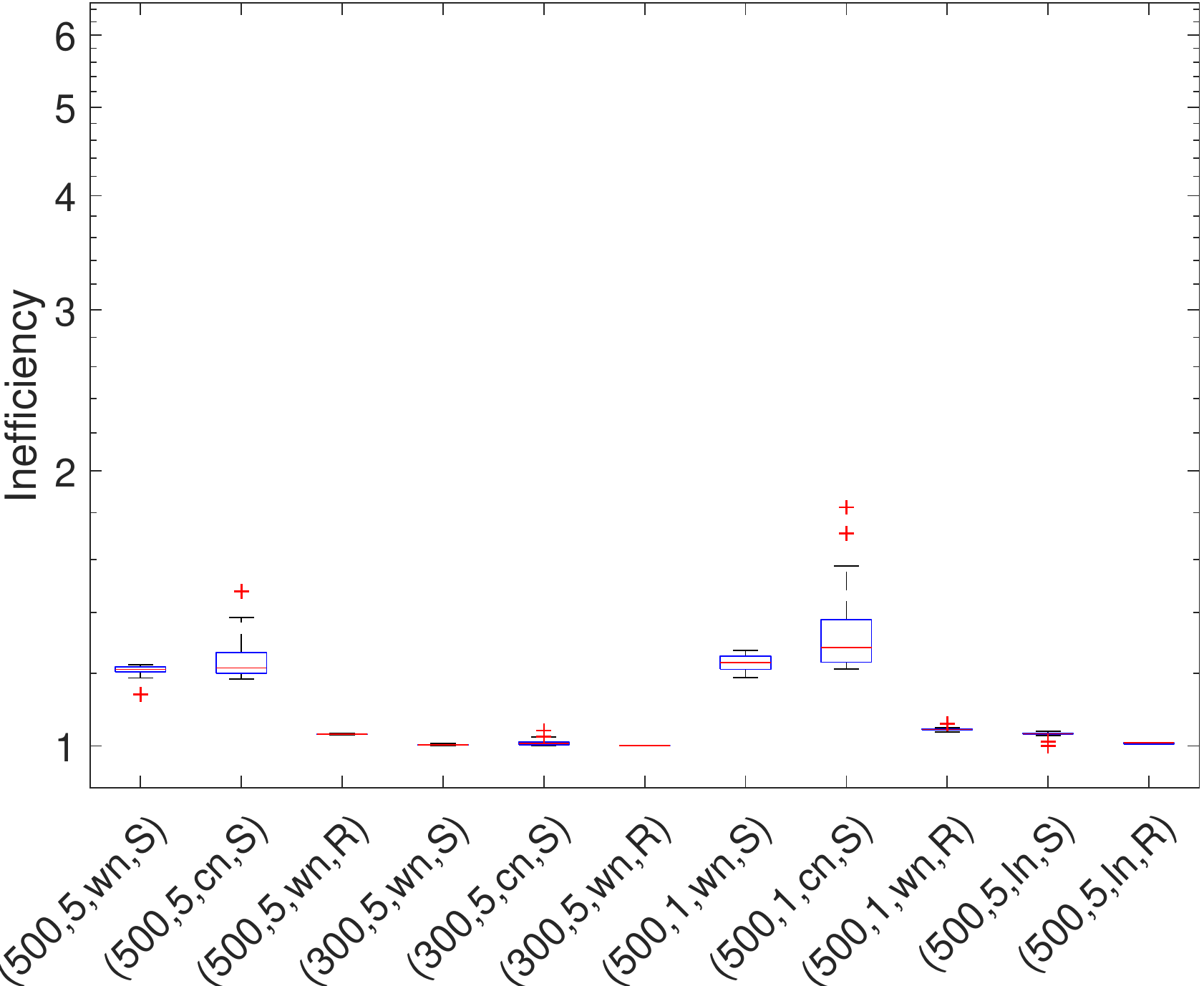}
\end{minipage}
\hspace*{0.1cm}
\begin{minipage}[t]{0.49\textwidth}
\includegraphics[trim= 20mm 65mm 19mm 92mm ,clip, width=\textwidth]{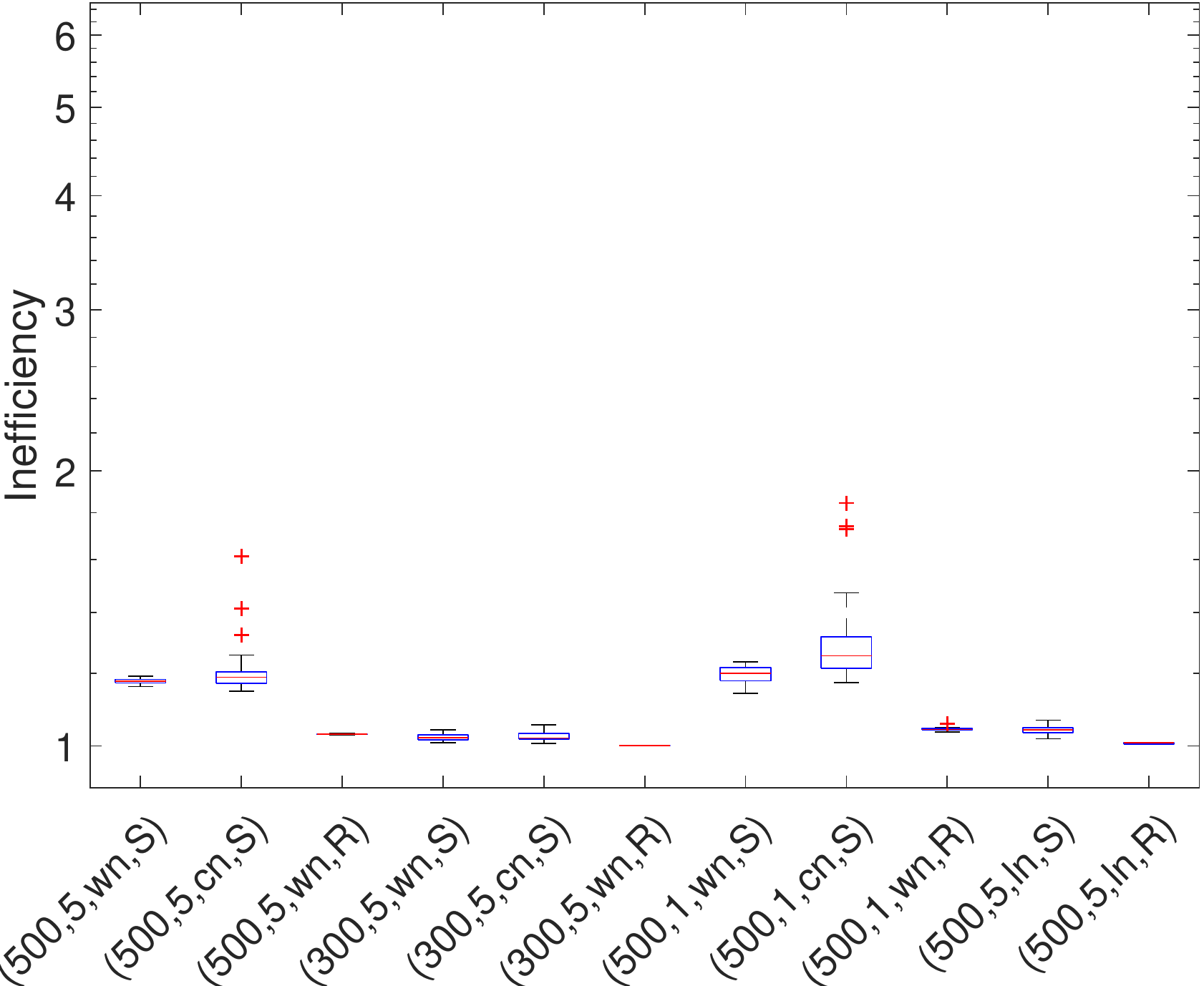}
\end{minipage}
\caption{LC for the RFMP (left) and the ROFMP (right).}
\label{fig:lc}
\end{figure}
The LC (see \cref{fig:lc}) yields good results in all test cases. We can remark that there are a few outliers and bigger boxes for the test cases with coloured noise and the scattered grid.

\subsection{Extrapolated Error method (EEM)}
\begin{figure}[!ht]
\begin{minipage}[t]{0.49\textwidth}
\includegraphics[trim= 20mm 65mm 19mm 92mm ,clip, width=\textwidth]{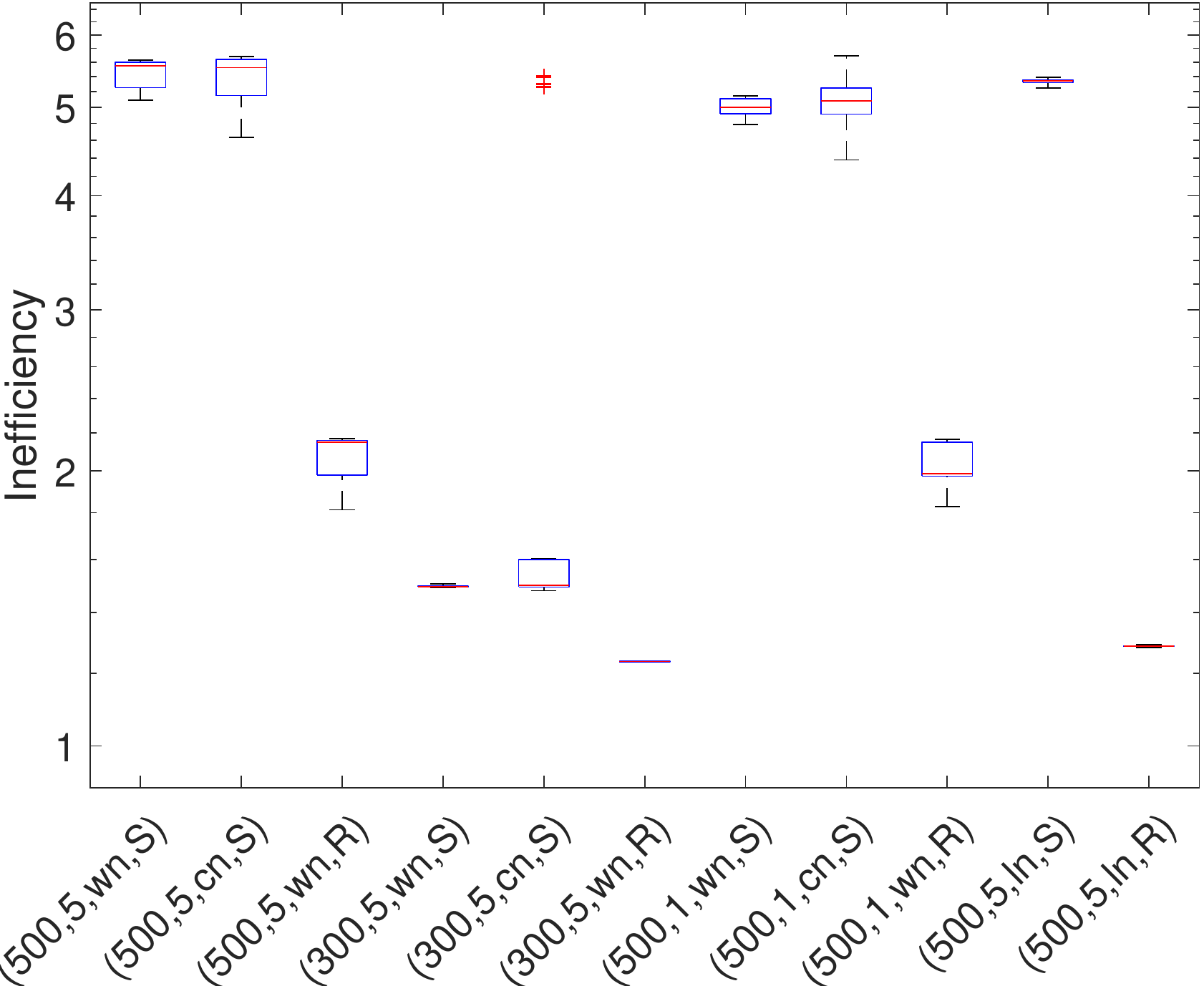}
\end{minipage}
\hspace*{0.1cm}
\begin{minipage}[t]{0.49\textwidth}
\includegraphics[trim= 20mm 65mm 19mm 92mm ,clip, width=\textwidth]{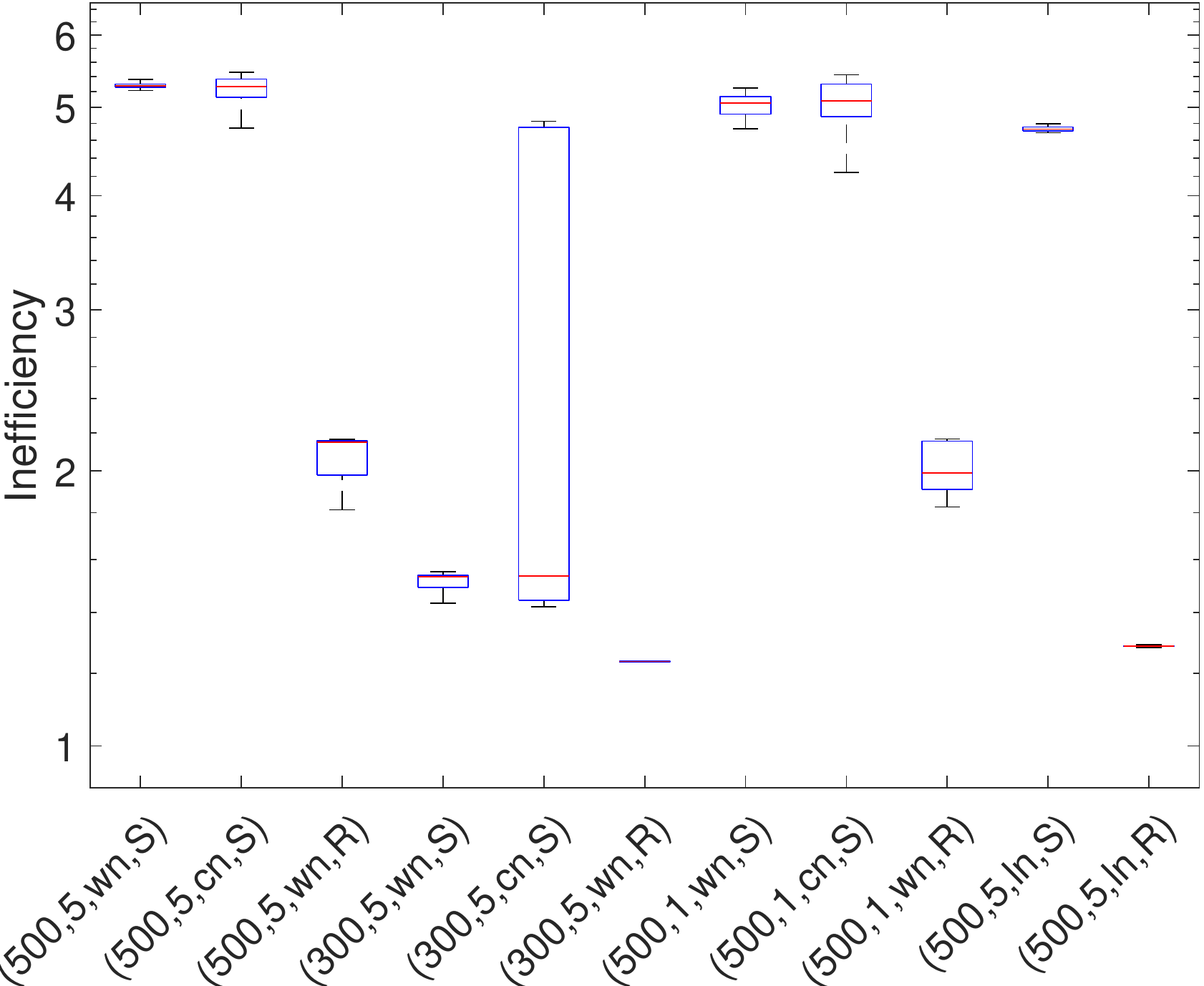}
\end{minipage}
\caption{EEM for the RFMP (left) and the ROFMP (right).}
\label{fig:eem}
\end{figure}
The EEM yields acceptable to rather poor results (see \cref{fig:eem}). We cannot observe any dependency on the grid or the kind of noise related to the acceptable results. Moreover, in the test case with a height of $300$km and an N2S of $5$\% with coloured noise we have some outliers for the RFMP and a large box for the ROFMP.

\subsection{Residual method (RM)}
\begin{figure}[!ht]
\begin{minipage}[t]{0.49\textwidth}
\includegraphics[trim= 20mm 65mm 19mm 92mm ,clip, width=\textwidth]{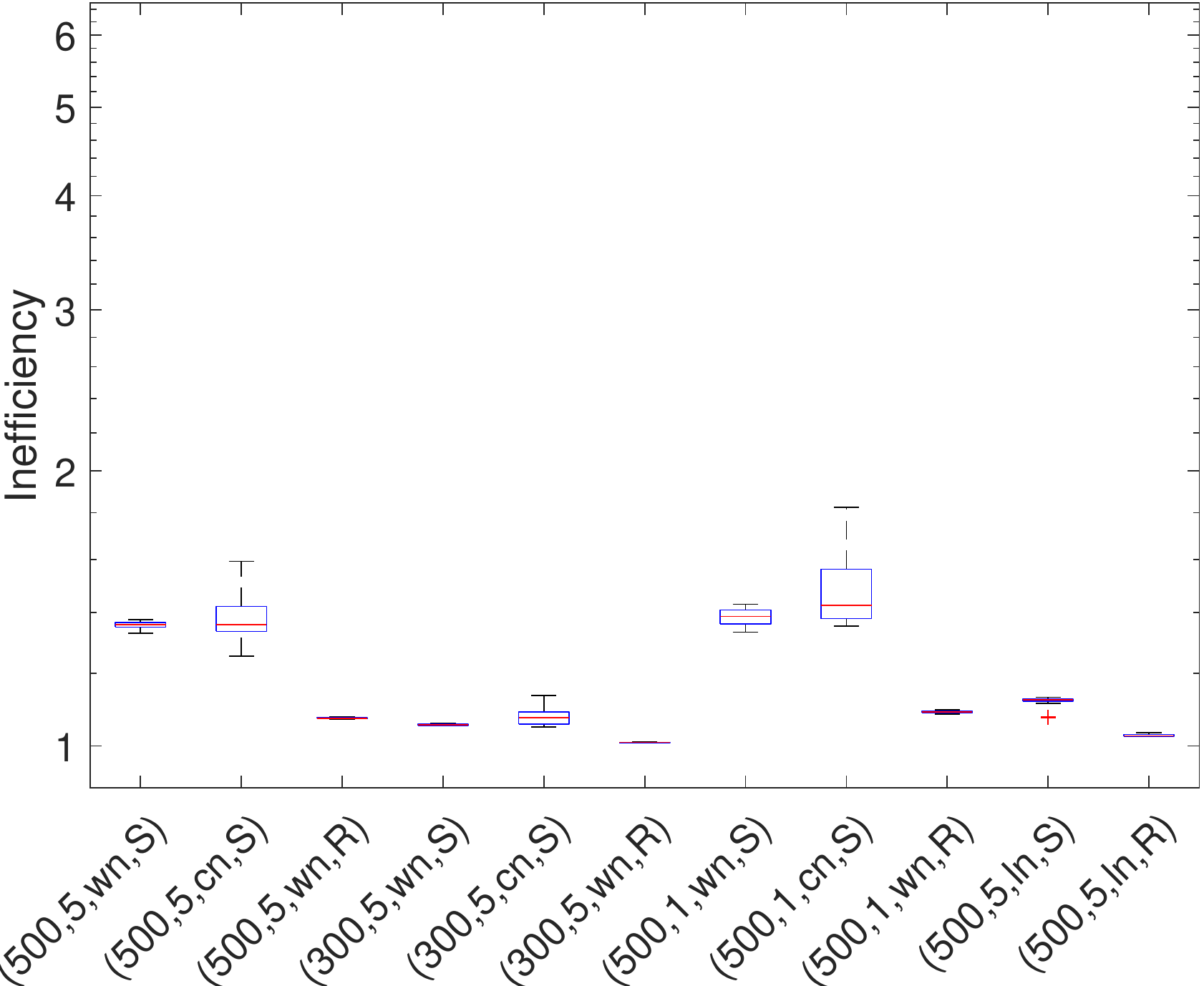}
\end{minipage}
\hspace*{0.1cm}
\begin{minipage}[t]{0.49\textwidth}
\includegraphics[trim= 20mm 65mm 19mm 92mm ,clip, width=\textwidth]{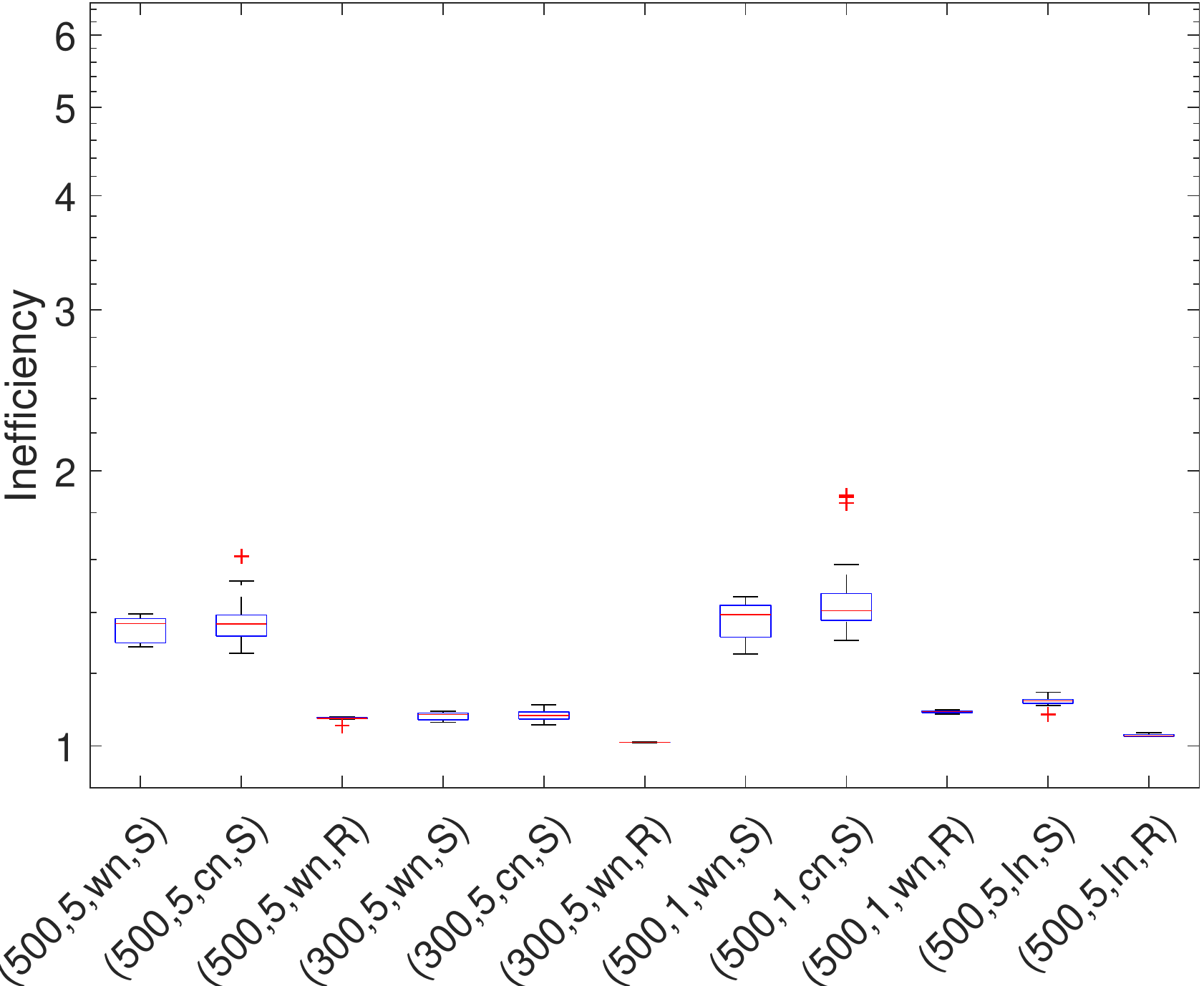}
\end{minipage}
\caption{RM for the RFMP (left) and the ROFMP (right).}
\label{fig:rm}
\end{figure}
The results for the RM (see \cref{fig:rm}) are good to acceptable in all test cases. We only have a few minor outliers.

\subsection{Generalized maximum likelihood (GML)}
\begin{figure}[!ht]
\begin{minipage}[t]{0.49\textwidth}
\includegraphics[trim= 20mm 65mm 19mm 92mm ,clip, width=\textwidth]{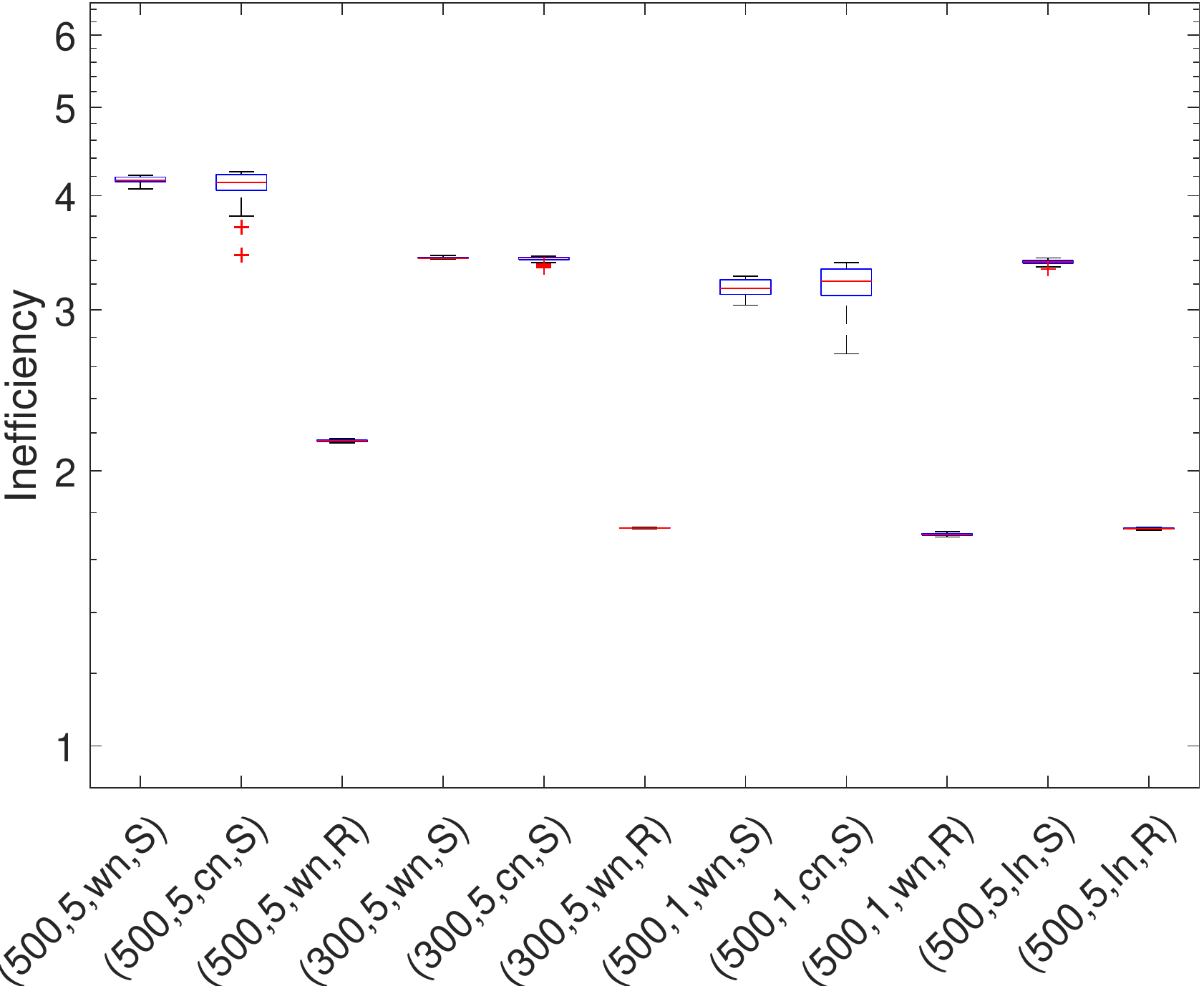}
\end{minipage}
\hspace*{0.1cm}
\begin{minipage}[t]{0.49\textwidth}
\includegraphics[trim= 20mm 65mm 19mm 92mm ,clip, width=\textwidth]{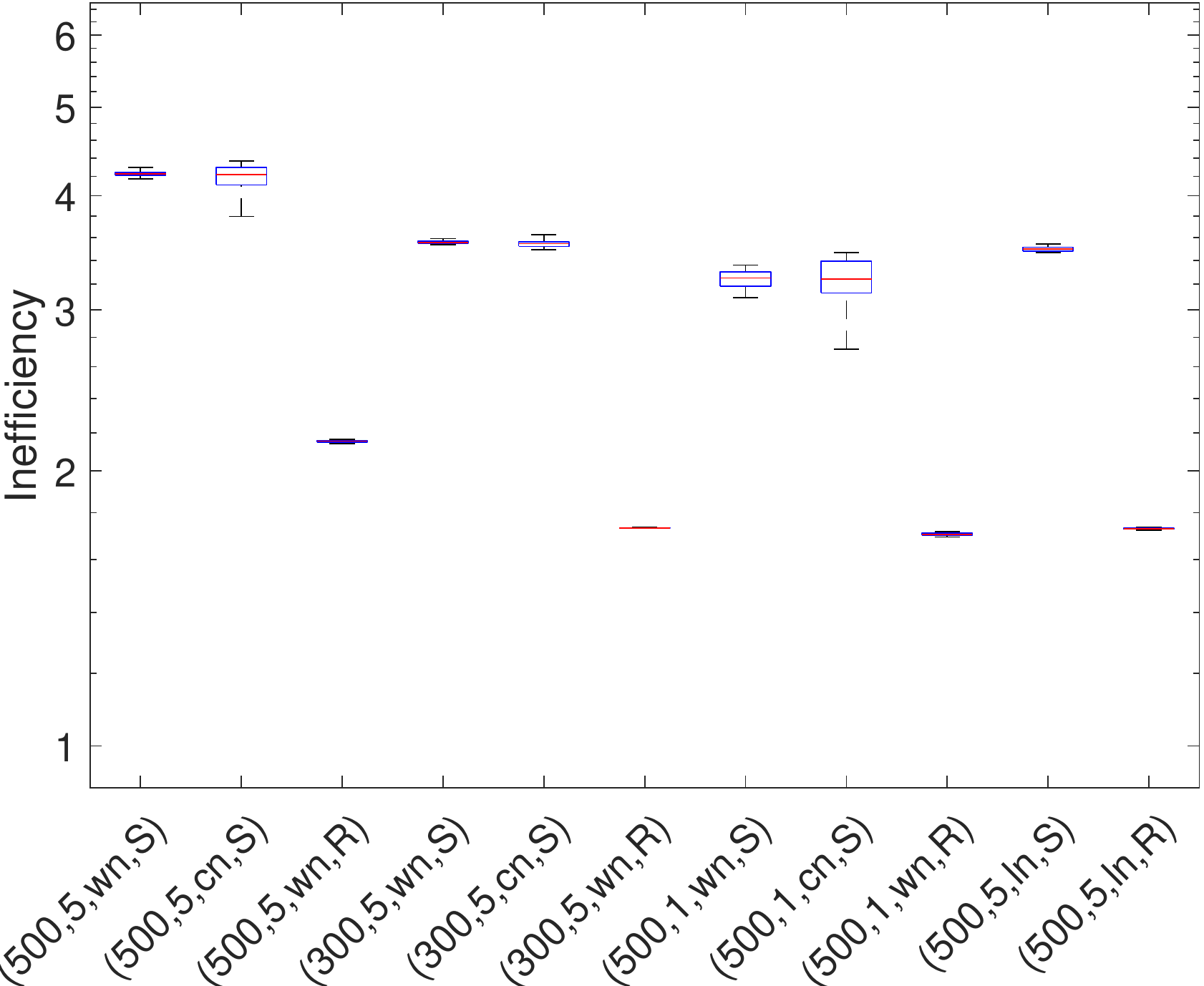}
\end{minipage}
\caption{GML for the RFMP (left) and the ROFMP (right).}
\label{fig:gml}
\end{figure}
In \cref{fig:gml}, we can see that the GML leads to acceptable results only in the case of the Reuter grid. In all cases of the scattered grid, its performance is rather bad. 

\subsection{Generalized cross validation (GCV)}
\begin{figure}[!ht]
\begin{minipage}[t]{0.49\textwidth}
\includegraphics[trim= 20mm 65mm 19mm 92mm ,clip, width=\textwidth]{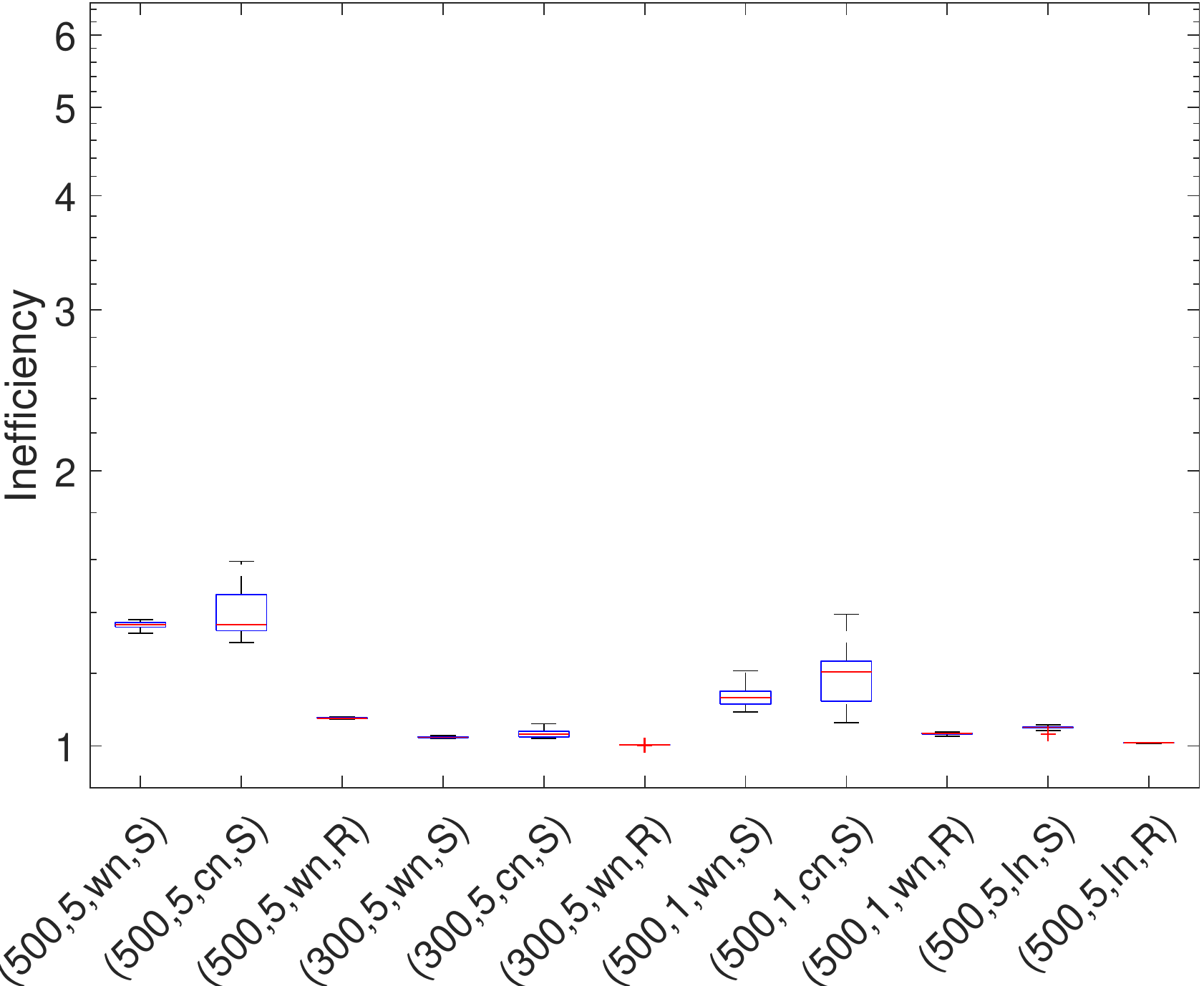}
\end{minipage}
\hspace*{0.1cm}
\begin{minipage}[t]{0.49\textwidth}
\includegraphics[trim= 20mm 65mm 19mm 92mm ,clip, width=\textwidth]{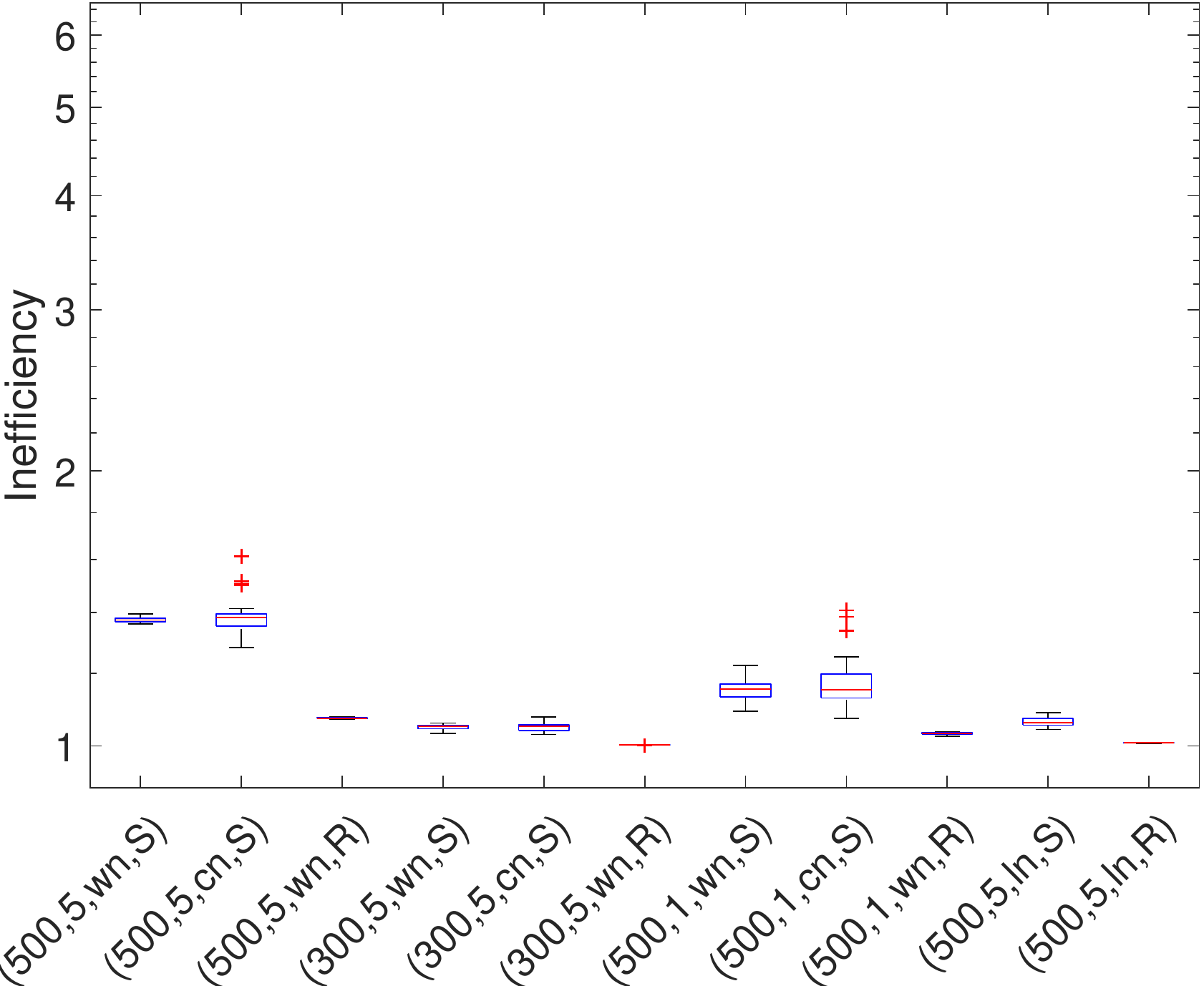}
\end{minipage}
\caption{GCV for the RFMP (left) and the ROFMP (right).}
\label{fig:gcv}
\end{figure}
From \cref{fig:gcv}, we can observe that the GCV yields good results in all test cases. We only have, in the case of the ROFMP, some minor outliers. It yields the best results with a more regularly distributed data grid.

\subsection{Robust generalized cross validation (RGCV)}
\begin{figure}[!ht]
\begin{minipage}[t]{0.49\textwidth}
\includegraphics[trim= 20mm 65mm 19mm 92mm ,clip, width=\textwidth]{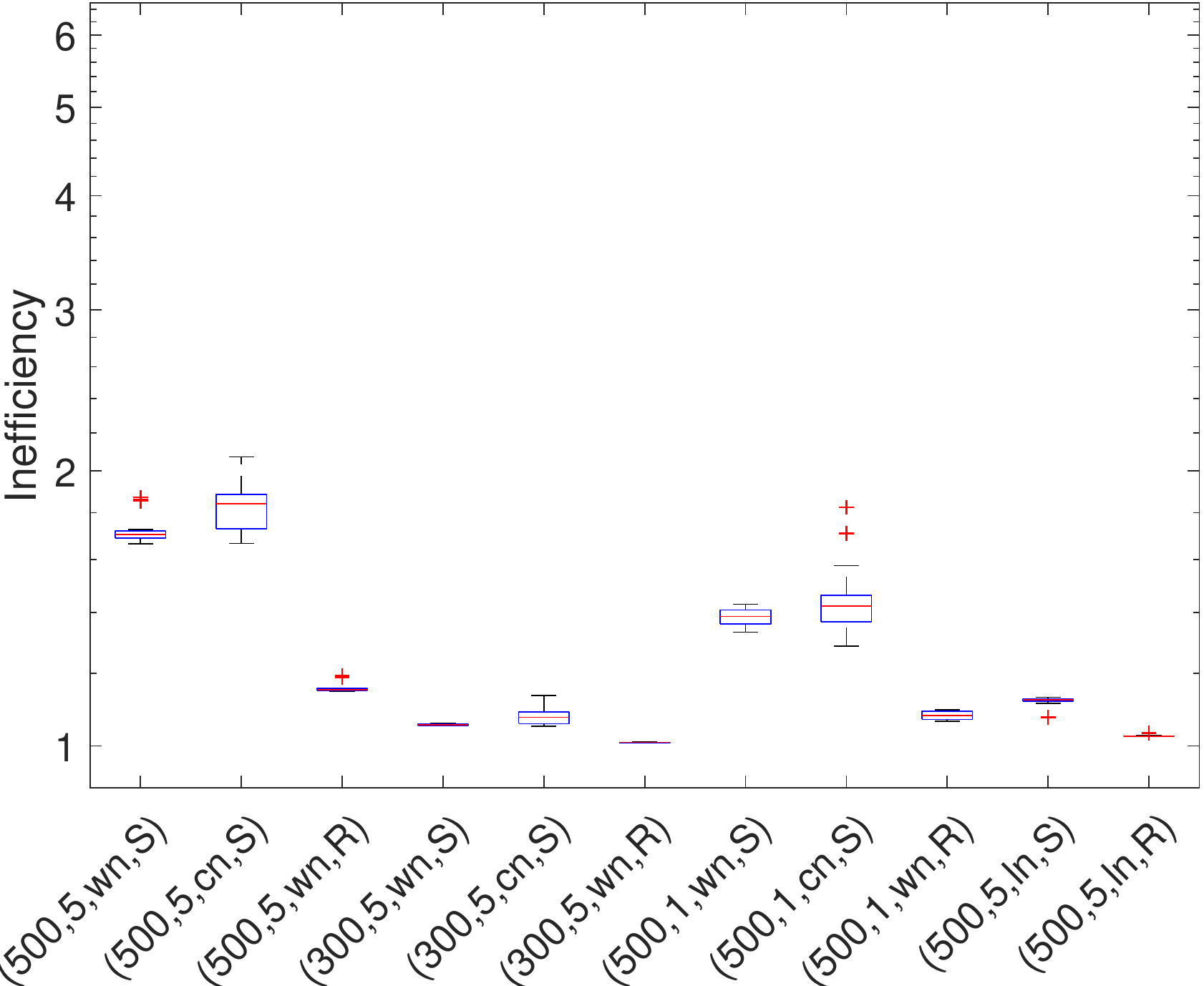}
\end{minipage}
\hspace*{0.1cm}
\begin{minipage}[t]{0.49\textwidth}
\includegraphics[trim= 20mm 65mm 19mm 92mm ,clip, width=\textwidth]{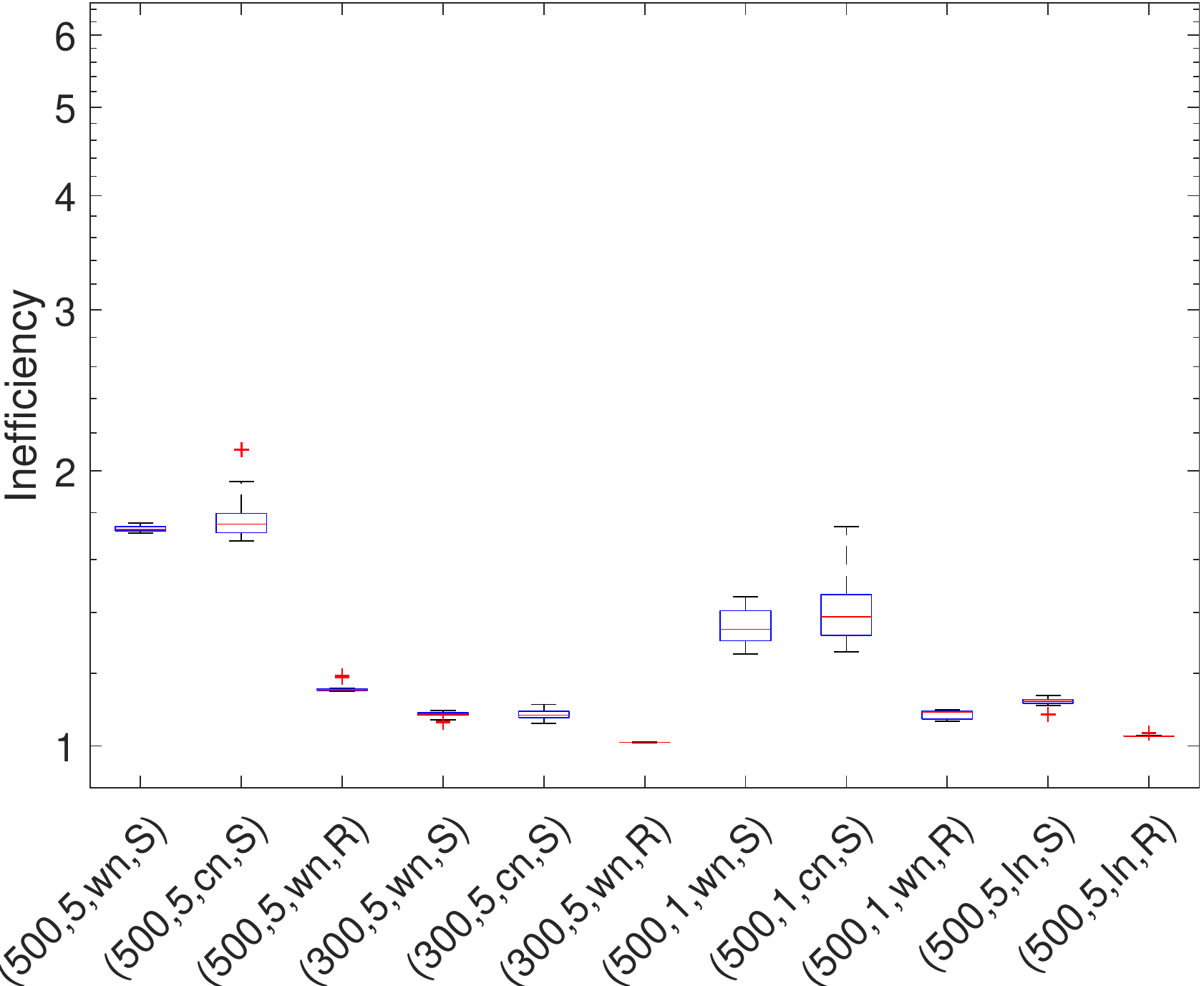}
\end{minipage}
\caption{RGCV for the RFMP (left) and the ROFMP (right).}
\label{fig:rgcv}
\end{figure}

The RGCV yields good to acceptable results (see \cref{fig:rgcv}) which get slightly worse and show a larger variance for a higher N2S or coloured noise scenarios.

\subsection{Strong robust generalized cross validation (SRGCV)}
\begin{figure}[!ht]
\begin{minipage}[t]{0.49\textwidth}
\includegraphics[trim= 20mm 65mm 19mm 92mm ,clip, width=\textwidth]{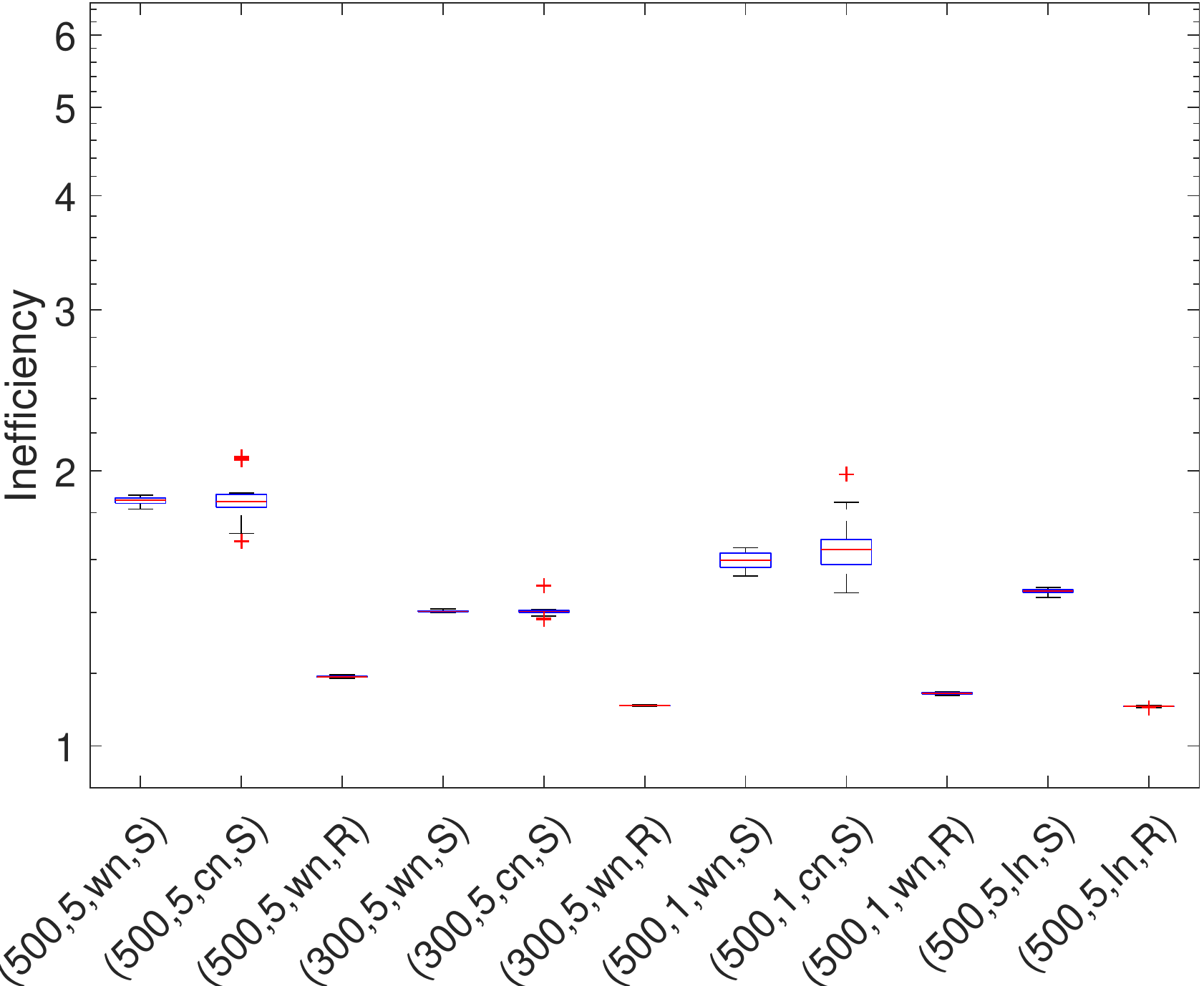}
\end{minipage}
\hspace*{0.1cm}
\begin{minipage}[t]{0.49\textwidth}
\includegraphics[trim= 20mm 65mm 19mm 92mm ,clip, width=\textwidth]{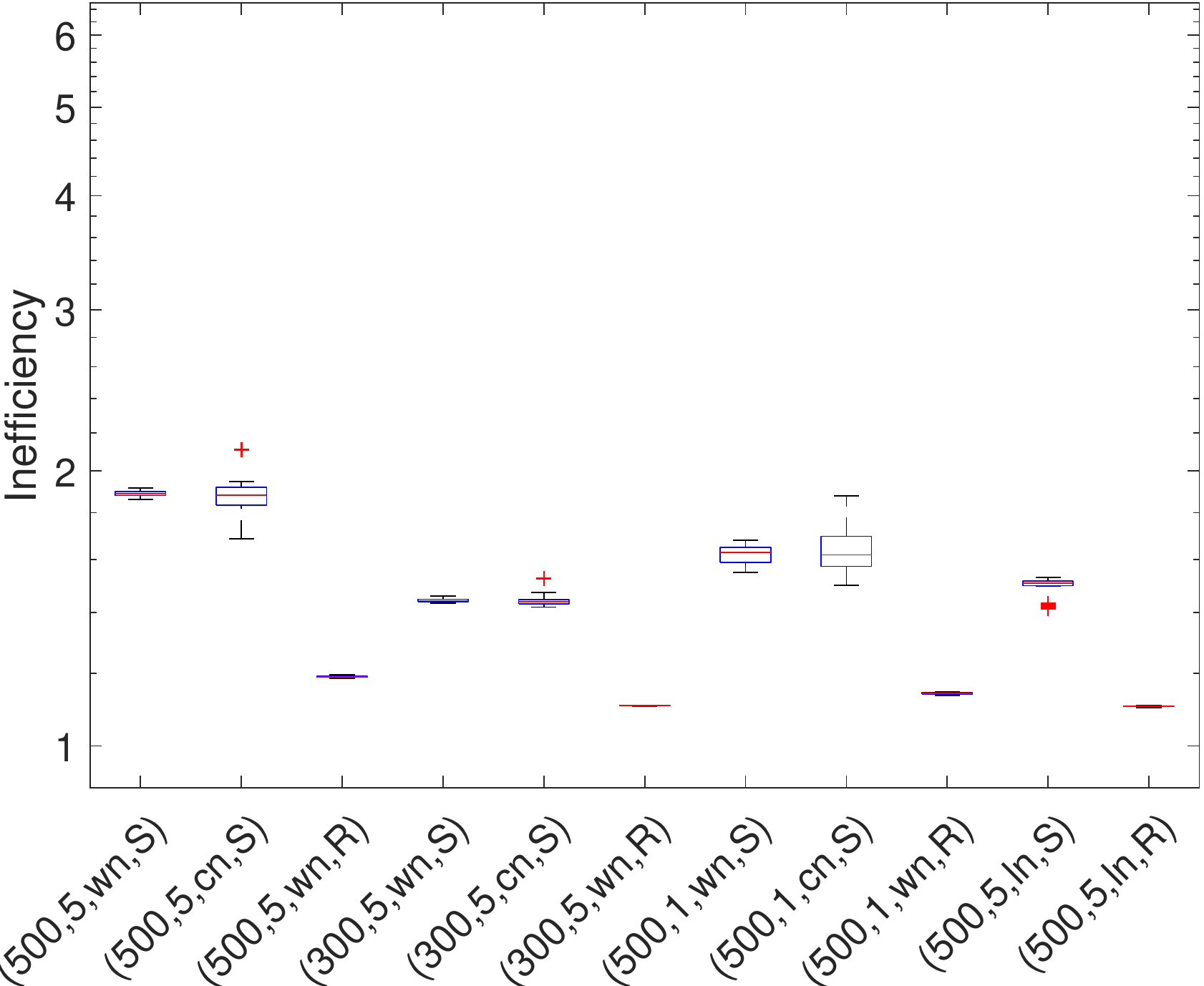}
\end{minipage}
\caption{SRGCV for the RFMP (left) and the ROFMP (right).}
\label{fig:srgcv}
\end{figure}
The SRGCV (see \cref{fig:srgcv}) has good to acceptable results in all the test cases which are a little bit worse than for the RGCV. The Reuter grid leads to good results whereas the scattered grid seems to be more difficult to handle by the method. 

\subsection{Modified generalized cross validation (MGCV)}
\begin{figure}[!ht]
\begin{minipage}[t]{0.49\textwidth}
\includegraphics[trim= 20mm 65mm 19mm 92mm ,clip, width=\textwidth]{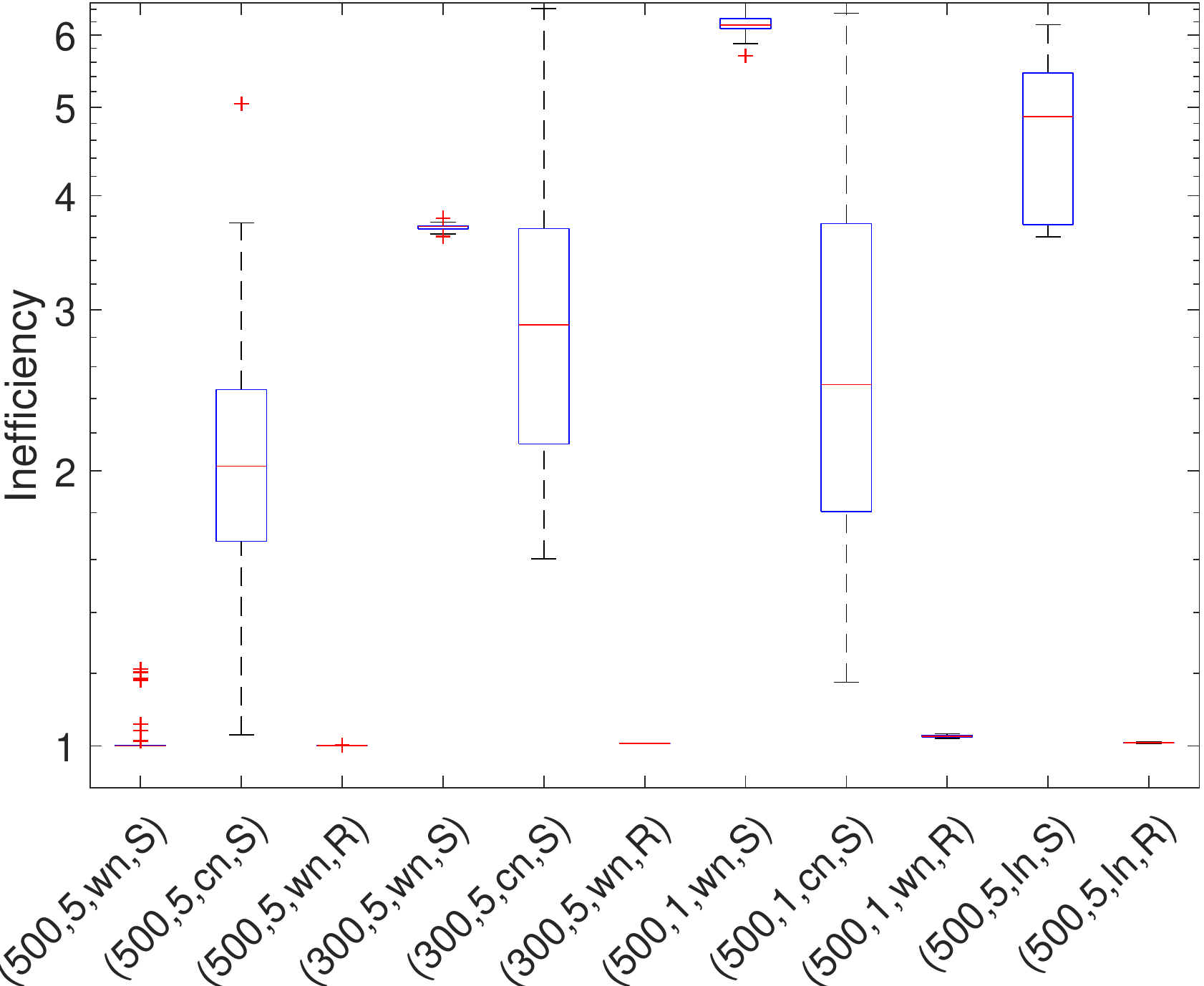}
\end{minipage}
\hspace*{0.1cm}
\begin{minipage}[t]{0.49\textwidth}
\includegraphics[trim= 20mm 65mm 19mm 92mm ,clip, width=\textwidth]{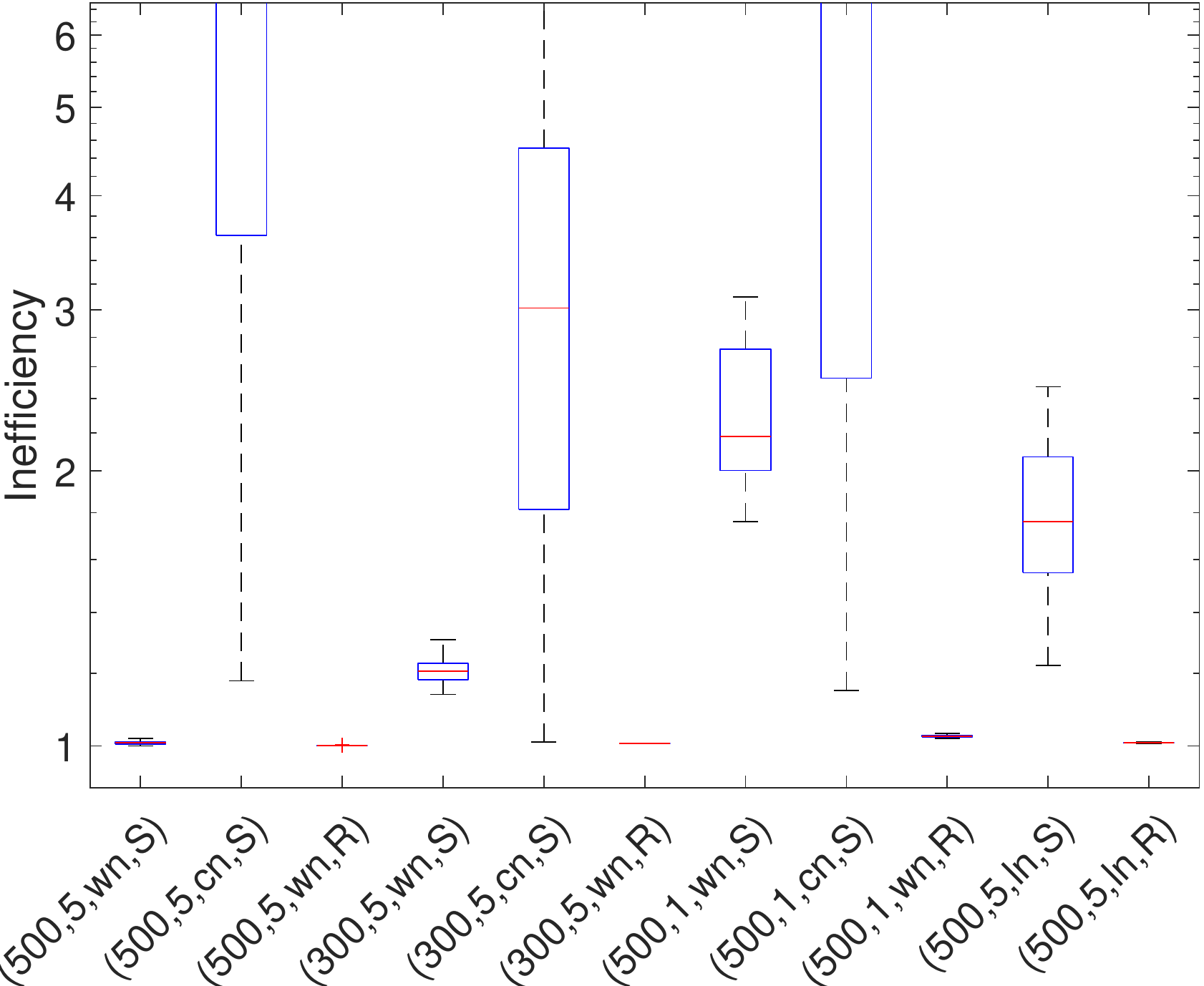}
\end{minipage}
\caption{MGCV for the RFMP (left) and the ROFMP (right).}
\label{fig:mgcv}
\end{figure}
The inefficiencies for the MGCV (see \cref{fig:mgcv}) for the test cases with white noise and the Reuter grid are good. In particular, in several of the cases with coloured noise the boxes are so big that they partially do not fit in the figure. Obviously, we get here a very large distribution of the inefficiencies. These cases seem to be very hard to handle for this method.

%Plots
\subsection{Plots of the results}
In this section, we show briefly the approximations of the gravitational potential which we obtain by the RFMP and the ROFMP for one typical noisy data set considering a good or a rather poor parameter choice.

\begin{figure}[!ht]
\centering
\begin{minipage}[t]{0.48\textwidth}
\includegraphics[trim=98mm 150mm 142mm 120mm, clip, width=\textwidth]{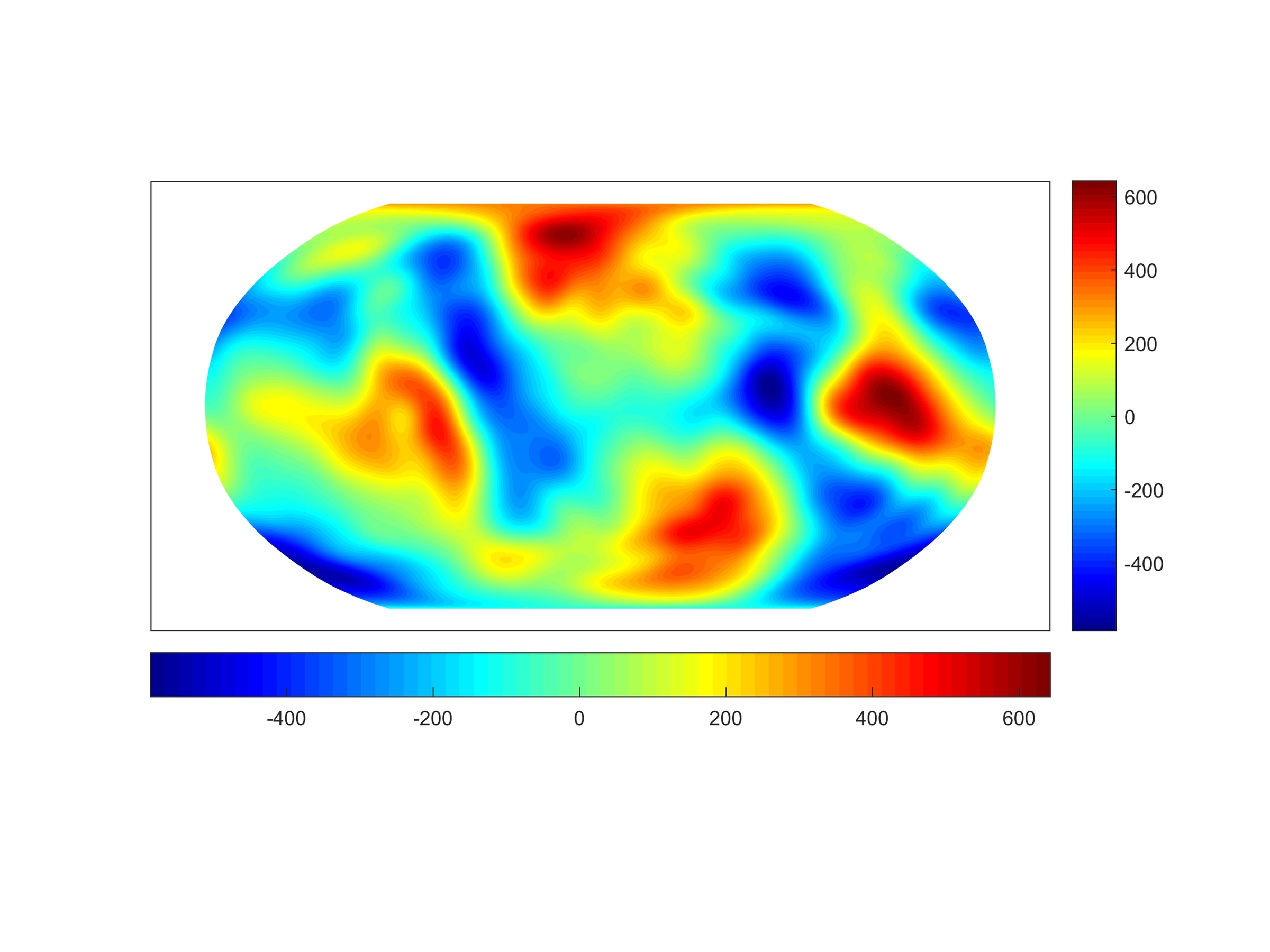}
\end{minipage}
\hspace*{0.2cm}
\begin{minipage}[t]{0.48\textwidth}
\includegraphics[trim=98mm 150mm 142mm 120mm, clip, width=\textwidth]{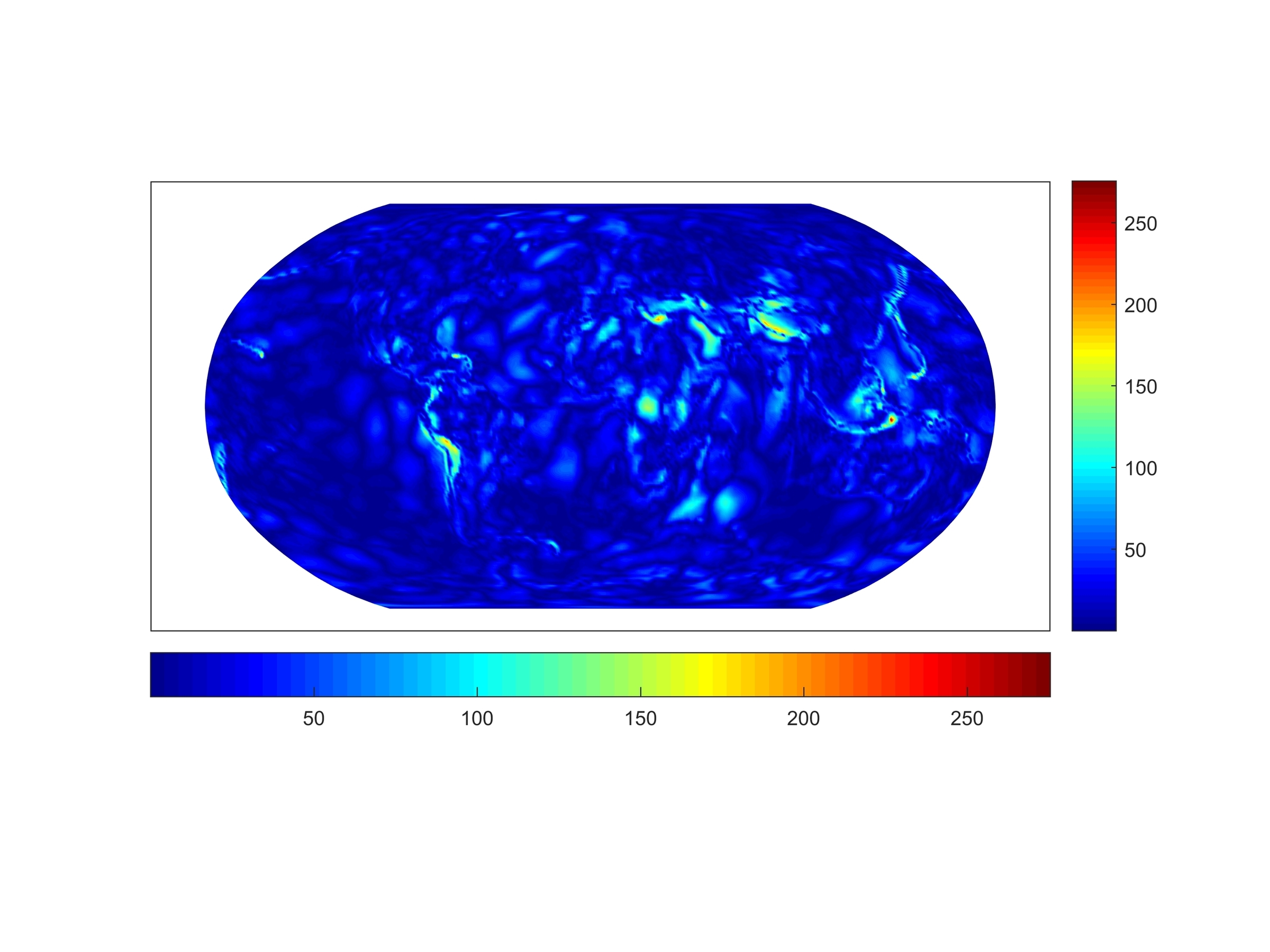}
\end{minipage}
\caption{The approximation from the RFMP for the best parameter (left) and the difference to the EGM2008 up to degree $100$ (right). Values in $\mathrm{m^2/s^2}$.}
\label{fig:rfmp_29}
\end{figure}

\begin{figure}[!ht]
\centering
\begin{minipage}[t]{0.48\textwidth}
\includegraphics[trim=98mm 150mm 142mm 120mm, clip, width=\textwidth]{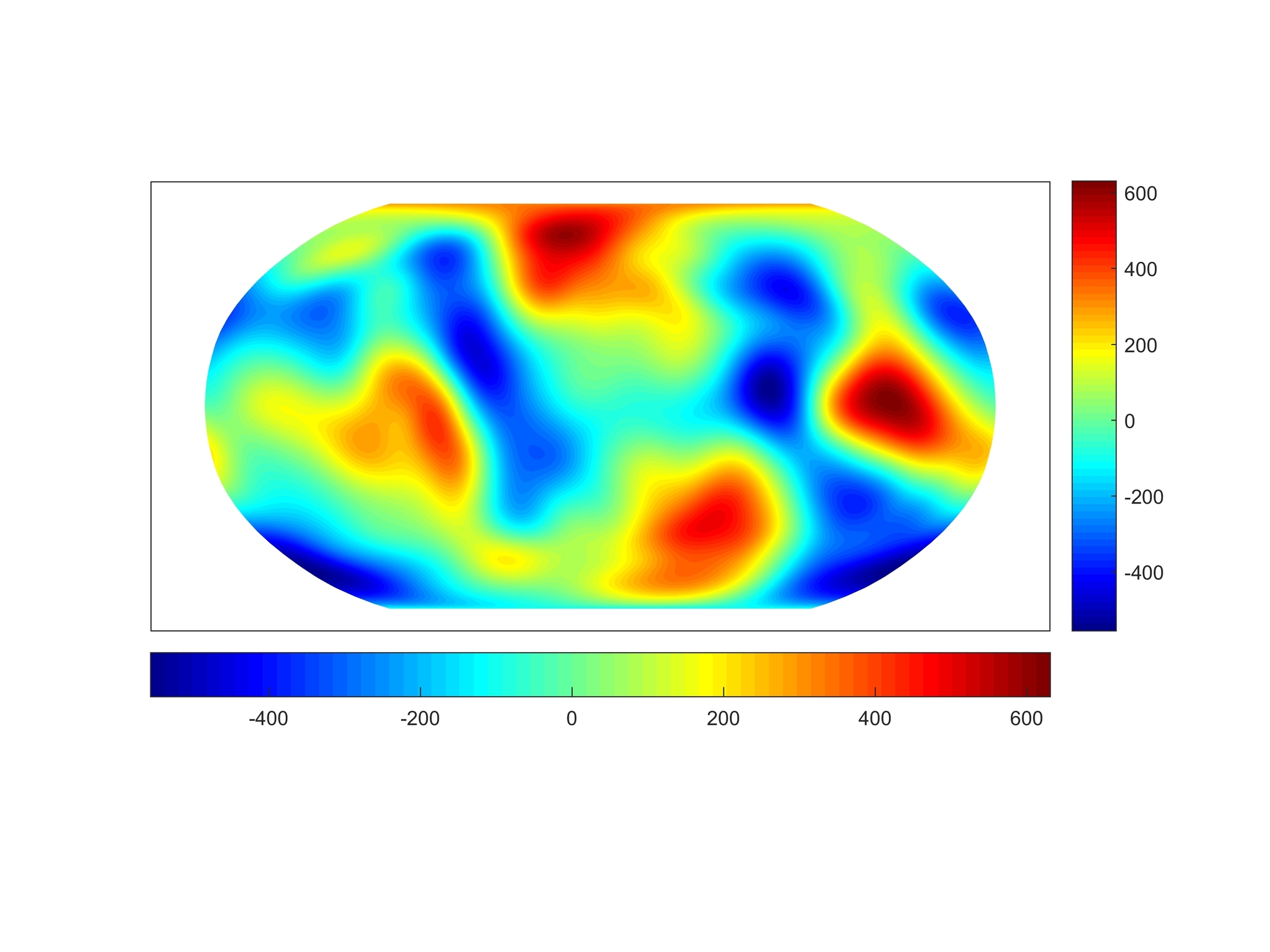}
\end{minipage}
\hspace*{0.2cm}
\begin{minipage}[t]{0.48\textwidth}
\includegraphics[trim=98mm 150mm 142mm 120mm, clip, width=\textwidth]{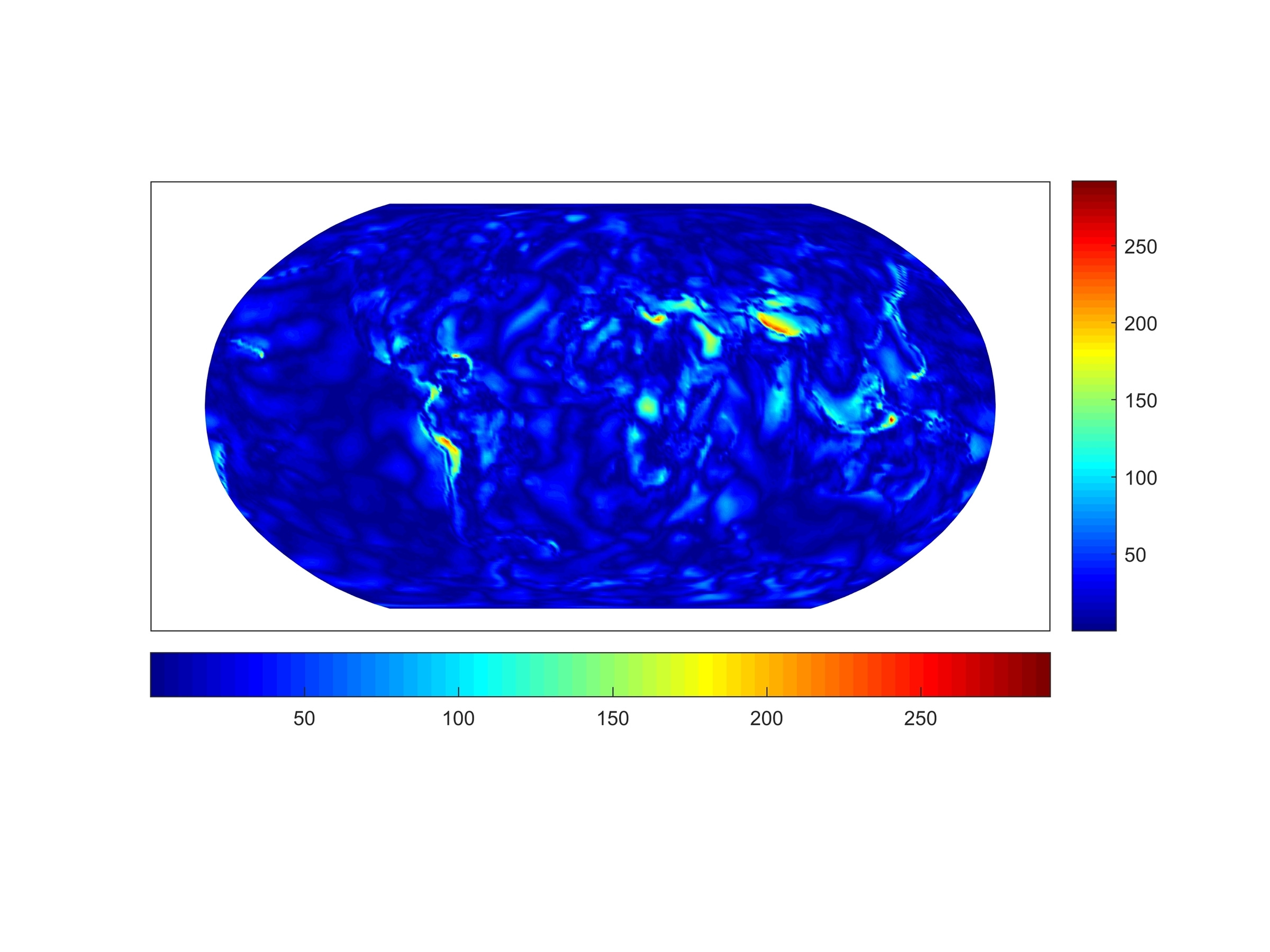}
\end{minipage}
\caption{The approximation from the RFMP for the parameter chosen by the GCV (left) and the difference to the EGM2008 up to degree $100$ (right). Values in $\mathrm{m^2/s^2}$. The inefficiency amounts to $1.16$.}
\label{fig:rfmp_22}
\end{figure}

\begin{figure}[!ht]
\centering
\begin{minipage}[t]{0.48\textwidth}
\includegraphics[trim=98mm 150mm 142mm 120mm, clip, width=\textwidth]{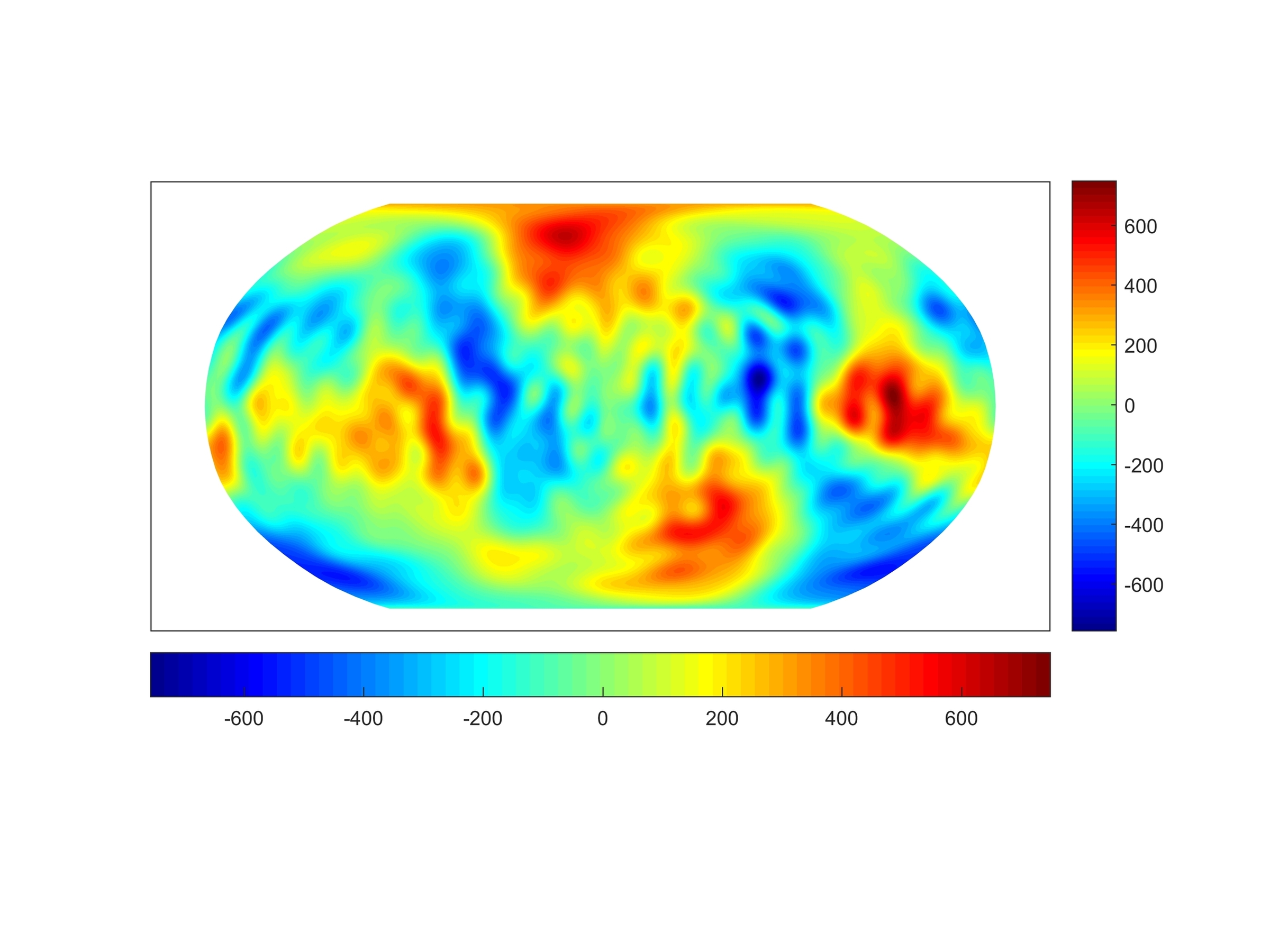}
\end{minipage}
\hspace*{0.2cm}
\begin{minipage}[t]{0.48\textwidth}
\includegraphics[trim=98mm 150mm 142mm 120mm, clip, width=\textwidth]{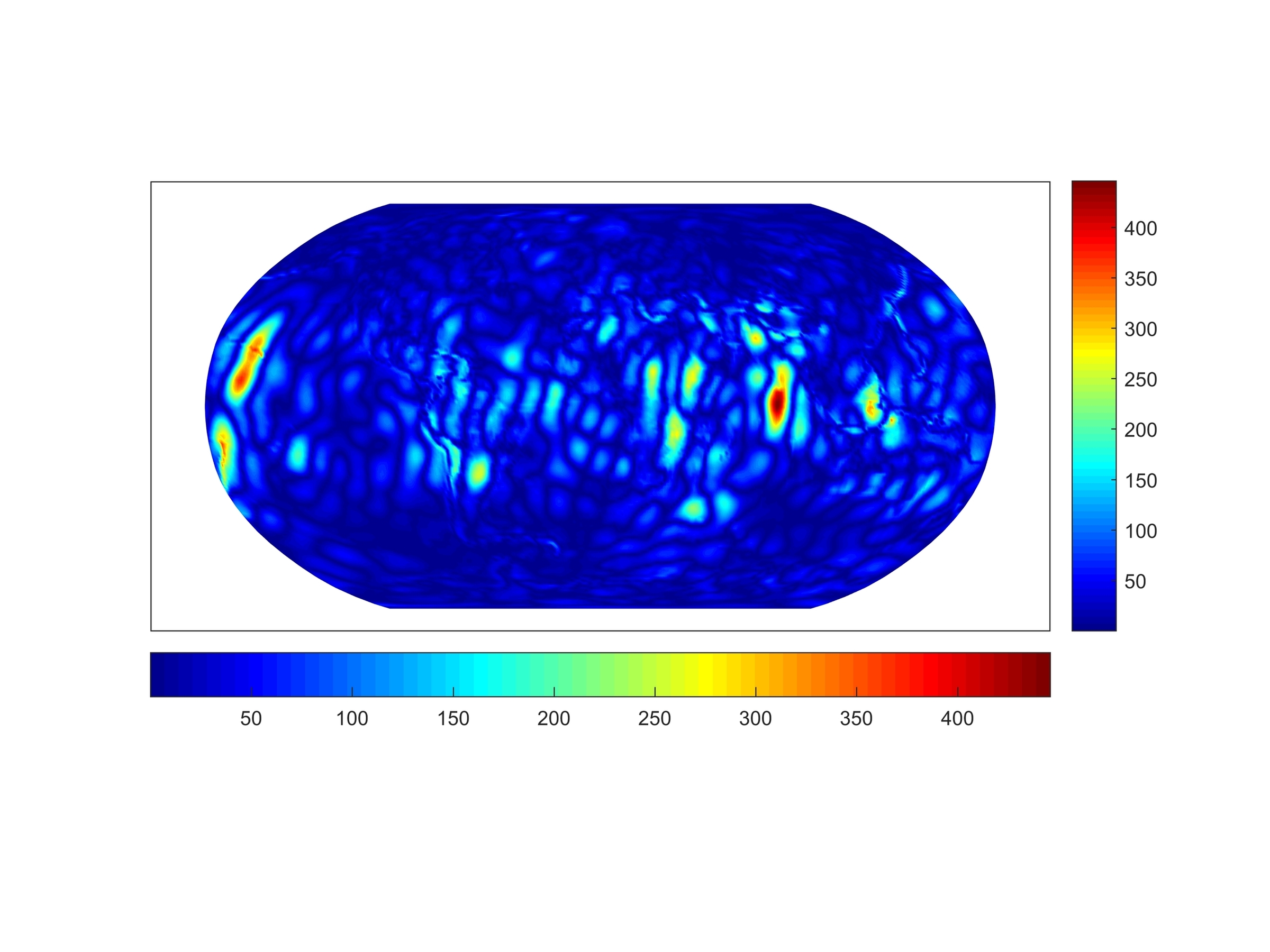}
\end{minipage}
\caption{The approximation from the RFMP for the parameter chosen by the MGCV (left) and the difference to the EGM2008 up to degree $100$ (right). Values in $\mathrm{m^2/s^2}$. The inefficiency amounts to $2.29$.}
\label{fig:rfmp_43}
\end{figure}

For the test case (500km, 5\%, coloured noise, scattered grid) with $\alpha=0.54$ for the AR(1)-process, \cref{fig:rfmp_29} shows the approximation which we obtain by the RFMP for the optimal regularization parameter $\lambda_{29}$ and the difference to the EGM2008 up to degree $100$.
\cref{fig:rfmp_22} shows the approximation belonging  to the regularization parameter $\lambda_{22}$ which is chosen by the GCV.
In \cref{fig:rfmp_43}, we can see the approximation belonging  to the parameter $\lambda_{43}$ which is chosen by the MGCV. We can see that the MGCV chooses the regularization parameter too small and with this choice we obtain a solution which is  underregularized. North-South oriented anomalies occur in the reconstruction which appear to be artefacts due to the noise along the simulated satellite tracks. In contrast, the approximation of the potential for the GCV-based parameter is only slightly worse than the result for the optimal parameter.


\begin{figure}[!ht]
\centering
\begin{minipage}[t]{0.48\textwidth}
\includegraphics[trim=98mm 150mm 142mm 120mm, clip, width=\textwidth]{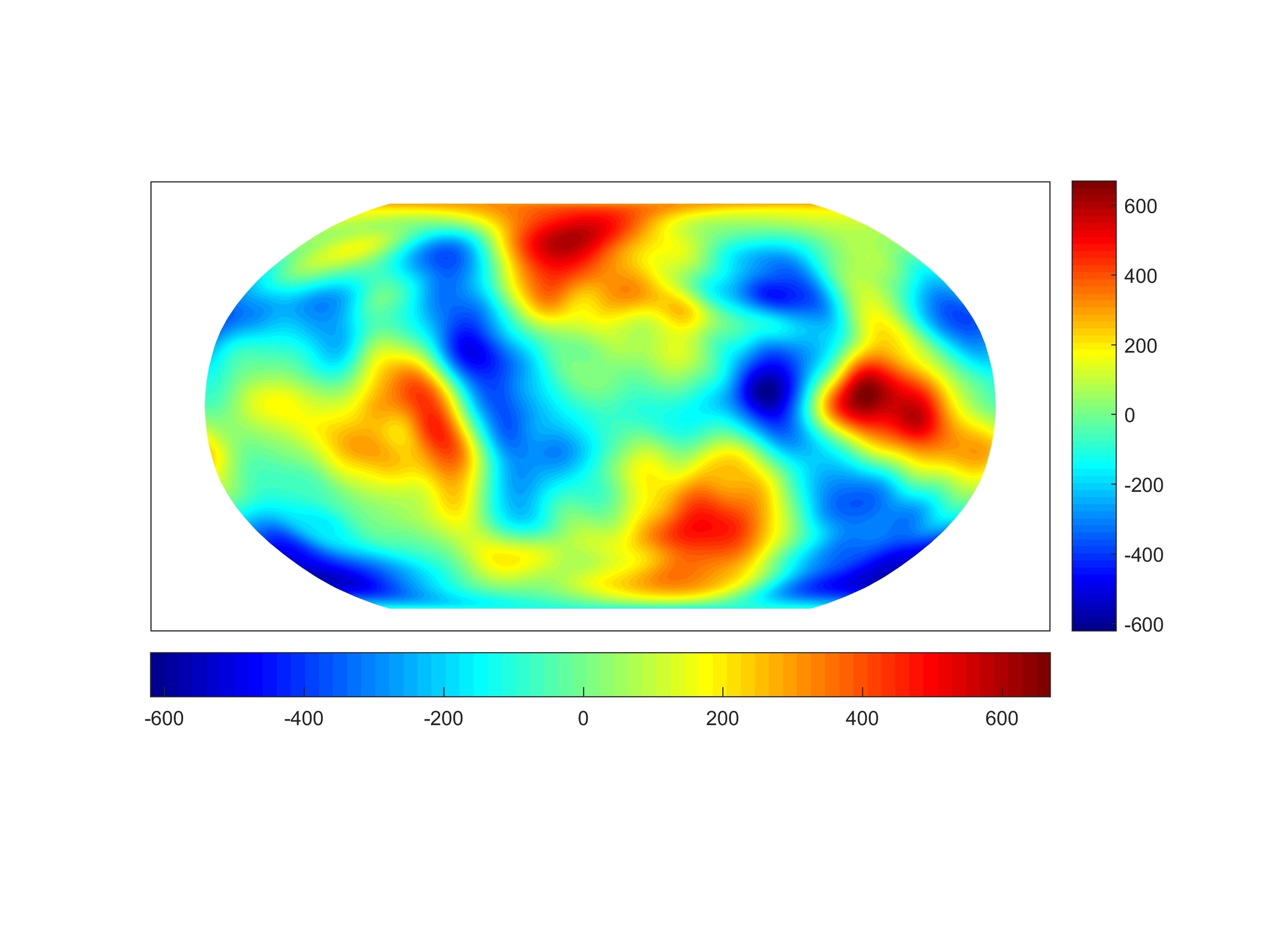}
\end{minipage}
\hspace*{0.2cm}
\begin{minipage}[t]{0.48\textwidth}
\includegraphics[trim=98mm 150mm 142mm 120mm, clip, width=\textwidth]{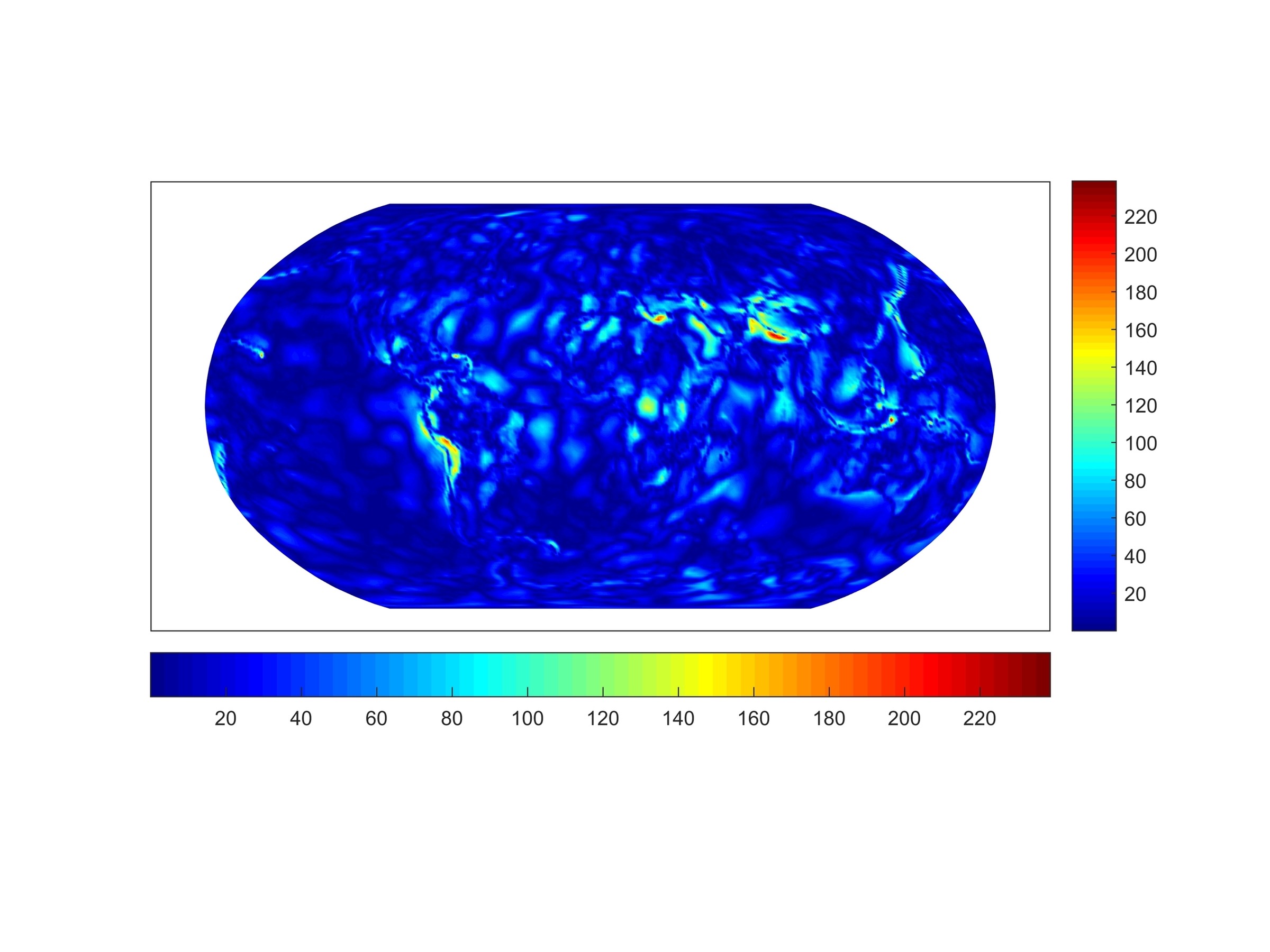}
\end{minipage}
\caption{The approximation from the ROFMP for the best parameter (left) and the difference to the EGM2008 up to degree $100$ (right). Values in $\mathrm{m^2/s^2}$.}
\label{fig:rofmp_29}
\end{figure}

\begin{figure}[!ht]
\centering
\begin{minipage}[t]{0.48\textwidth}
\includegraphics[trim=98mm 150mm 142mm 120mm, clip, width=\textwidth]{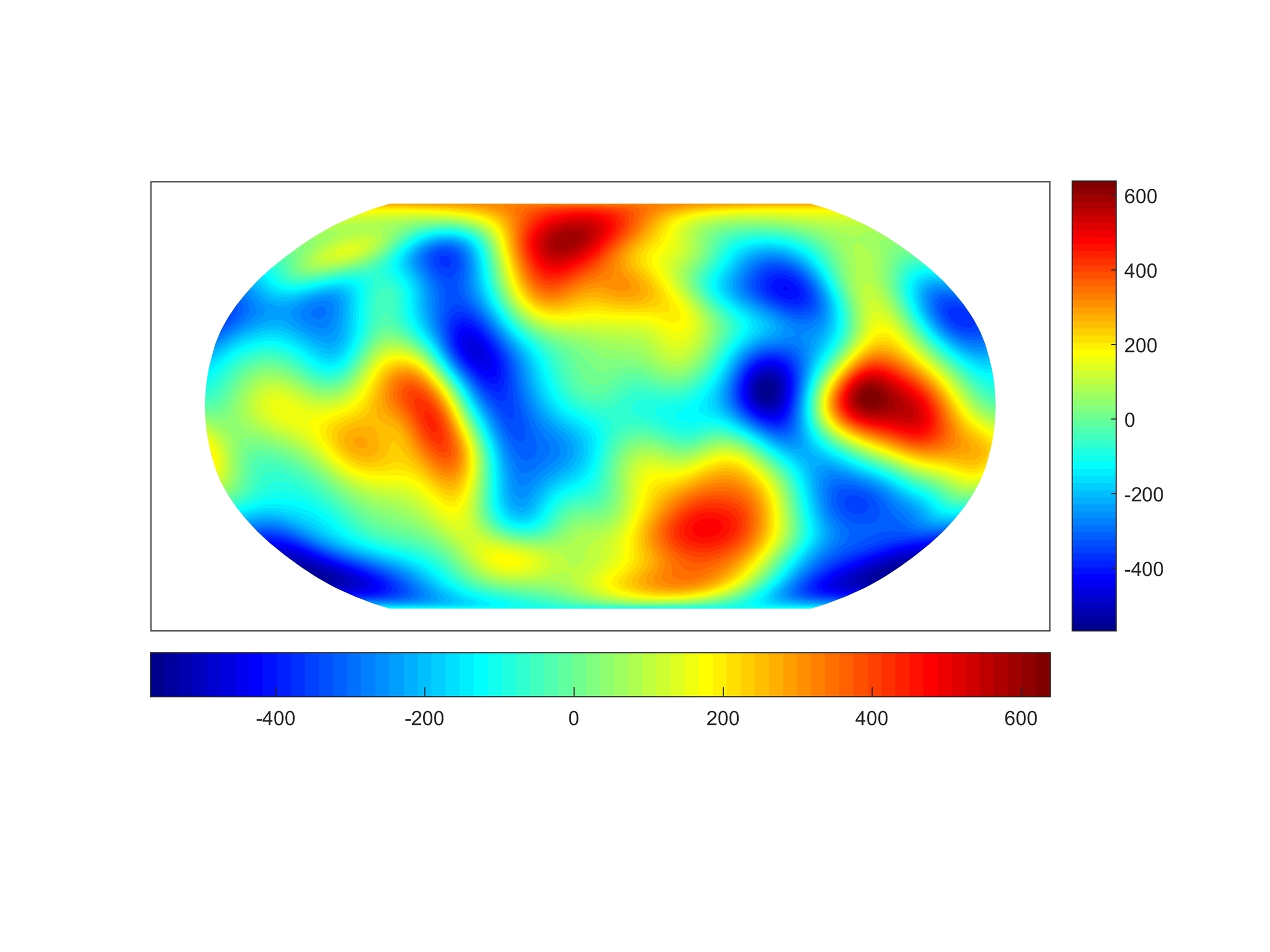}
\end{minipage}
\hspace*{0.2cm}
\begin{minipage}[t]{0.48\textwidth}
\includegraphics[trim=98mm 150mm 142mm 120mm, clip, width=\textwidth]{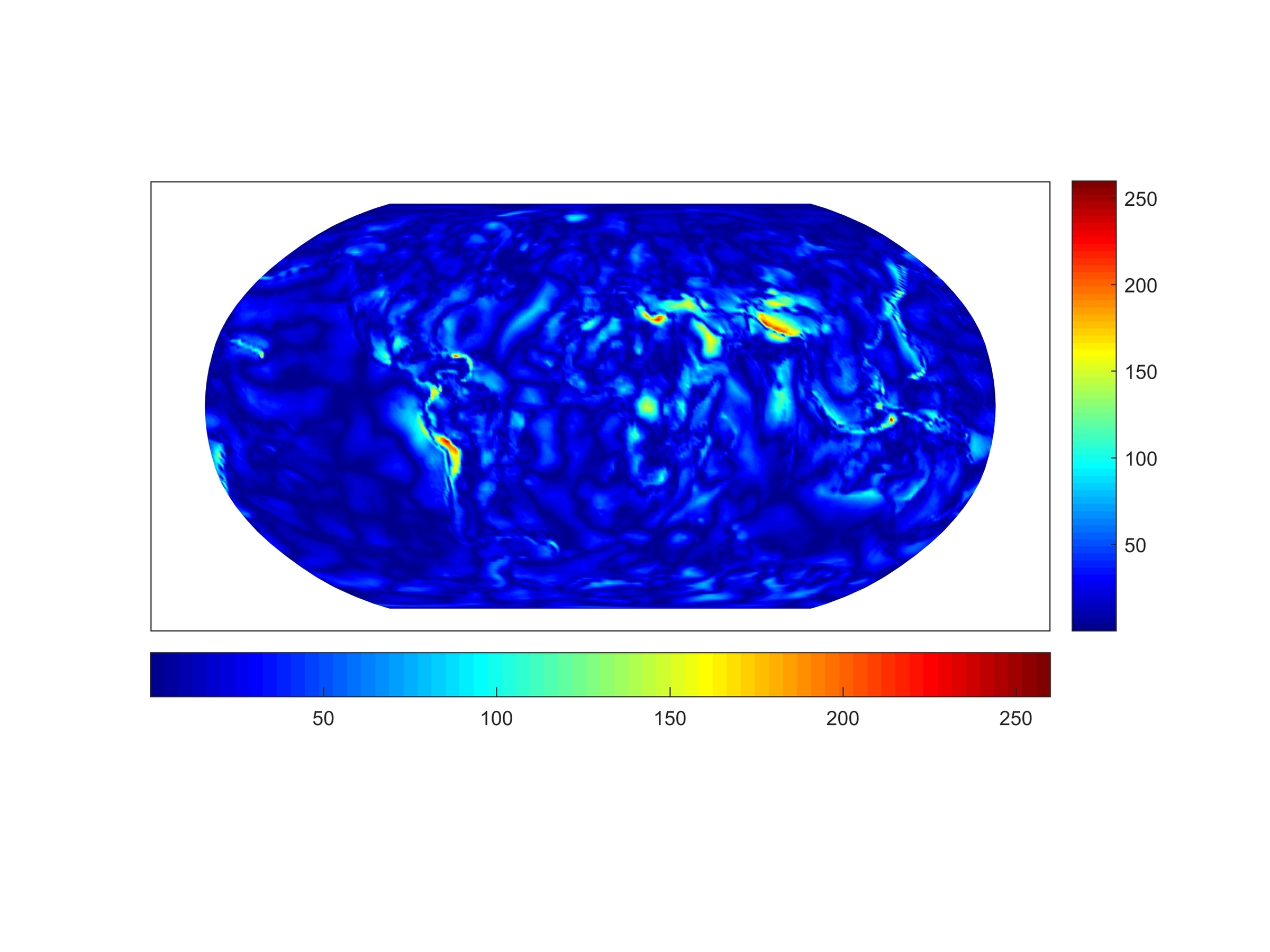}
\end{minipage}
\caption{The approximation from the ROFMP for the parameter chosen by the GCV (left) and the difference to the EGM2008 up to degree $100$ (right). Values in $\mathrm{m^2/s^2}$. The inefficiency amounts to $1.14$.}
\label{fig:rofmp_22}
\end{figure}

\begin{figure}[!ht]
\centering
\begin{minipage}[t]{0.48\textwidth}
\includegraphics[trim=98mm 150mm 142mm 120mm, clip, width=\textwidth]{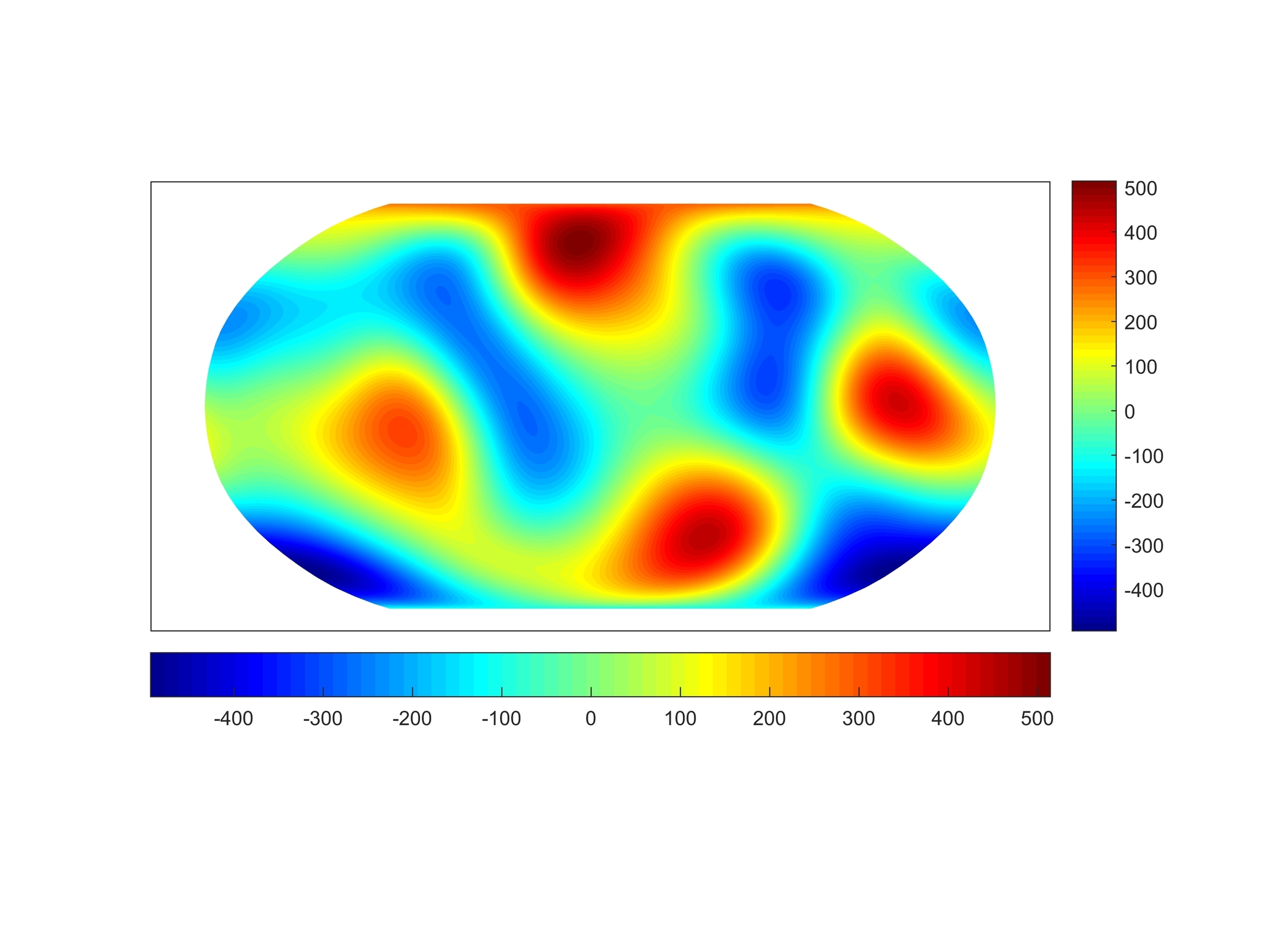}
\end{minipage}
\hspace*{0.2cm}
\begin{minipage}[t]{0.48\textwidth}
\includegraphics[trim=98mm 150mm 142mm 120mm, clip, width=\textwidth]{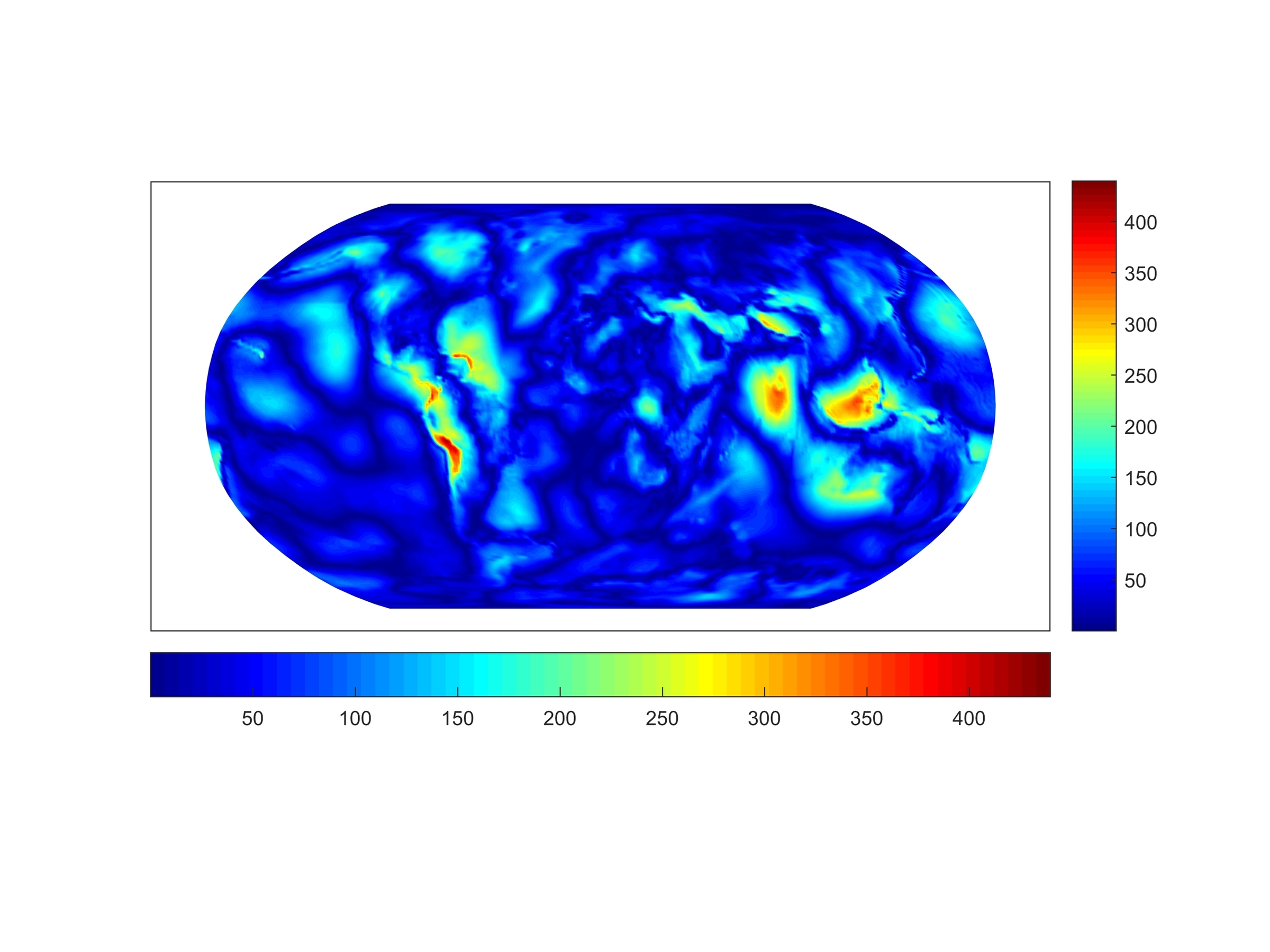}
\end{minipage}
\caption{The approximation from the ROFMP for the parameter chosen by the GML (left) and the difference to the EGM2008 up to degree $100$ (right). Values in $\mathrm{m^2/s^2}$. The inefficiency amounts to $3.03$.}
\label{fig:rofmp_9}
\end{figure}

Furthermore, we show the same test case as above but with the approximation from the ROFMP with $\alpha=0.56$ in the AR(1)-process. \cref{fig:rofmp_29} shows the approximation for the optimal parameter $\lambda_{29}$ and the difference to EGM2008.
In \cref{fig:rofmp_22}, we see the approximation which belongs to the parameter $\lambda_{22}$ which is chosen by the GCV. 
\cref{fig:rofmp_9} shows the approximation with the regularization parameter $\lambda_{9}$ which is chosen by the GML.

Here, the GML chooses a regularization parameter which is too large that means our approximation is overregularized. We get less information and details about the gravitational potential. Essential details such as signals due to the Andes or the region around Indonesia occur in the difference plot -- much stronglier than for the other examples. Again the parameter choice of the GCV yields a good approximation for the gravitational potential. 

\begin{figure}[!ht]
\centering
\includegraphics[trim=25mm 75mm 15mm 90mm, clip, width=0.8\textwidth]{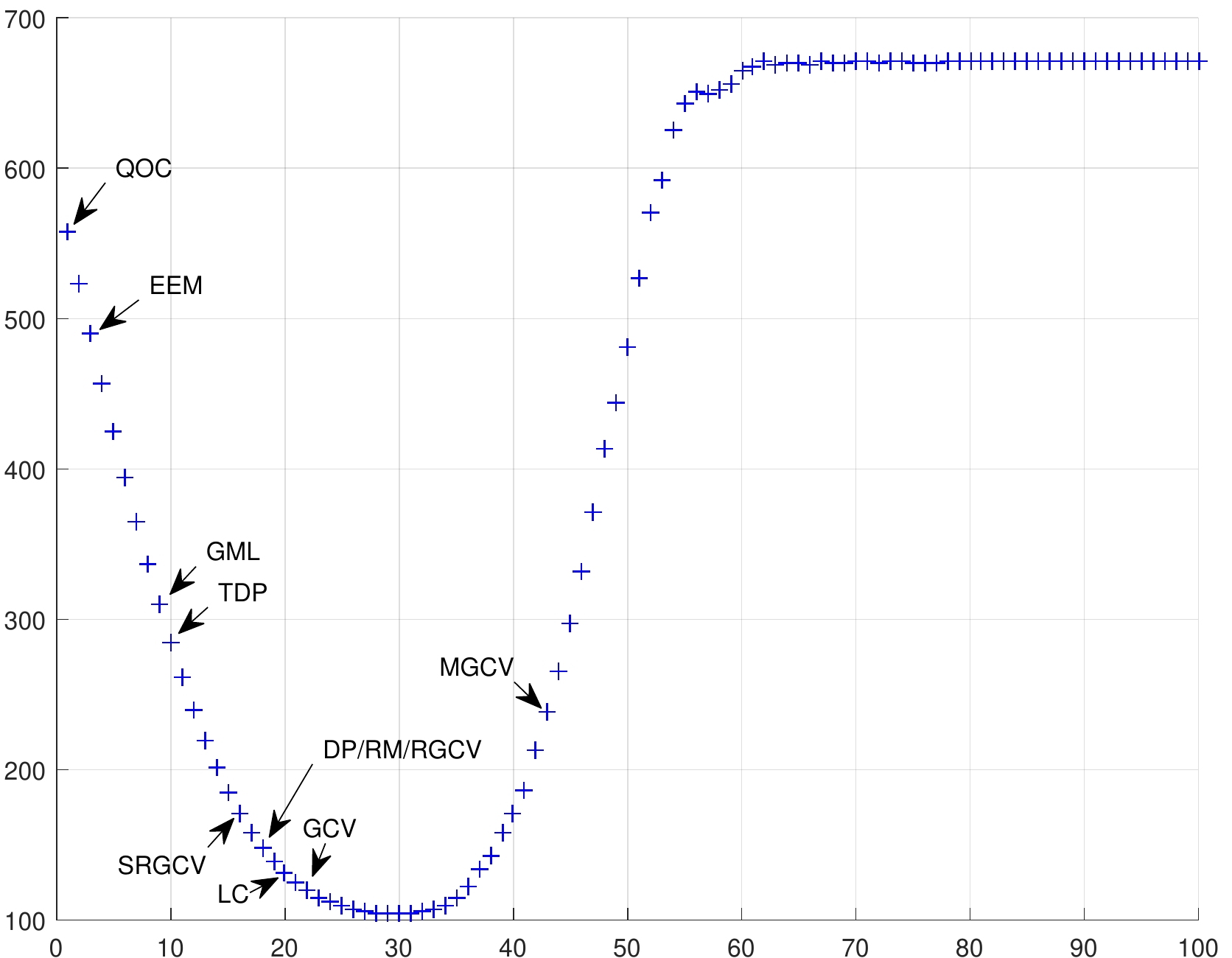}
\caption{The horizontal axis states the index $k$ of the regularization parameter and the vertical axis shows $\|x-x_{k}^\eps\|_{\LO}$ for the RFMP.}
\label{fig:diff1}
\end{figure}

\begin{figure}[!ht]
\centering
\includegraphics[trim=23mm 75mm 15mm 90mm, clip, width=0.8\textwidth]{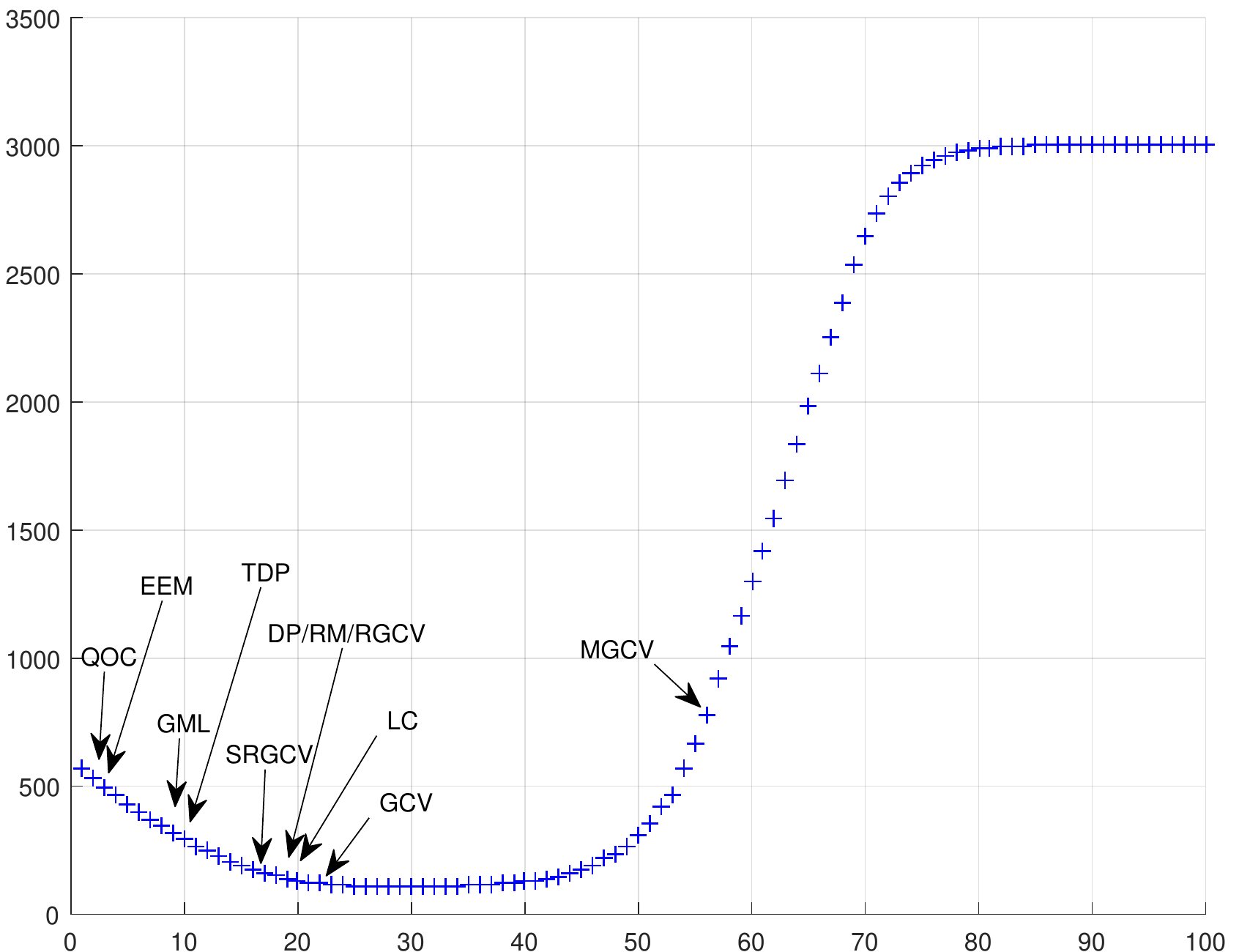}
\caption{The horizontal axis states the index $k$ of the regularization parameter and the vertical axis shows $\|x-x_{k}^\eps \|_{\LO}$ for the ROFMP.}
\label{fig:diff2}
\end{figure}

Finally, \cref{fig:diff1,fig:diff2} show the difference $\| x-x_k^\eps \|_{\LO}$ between the original solution (i.e. EGM2008 up to degree $100$) and the approximation $x_k^\eps$ obtained for the different regularization parameters which were chosen by the considered strategies. The horizontal axis states the index $k$ of the regularization parameter $\lambda_{k}$. The plots refer to the same scenario as Figures \ref{fig:rfmp_29} to \ref{fig:rofmp_9}. The arrows show the parameters which are chosen by the methods. The diagrams confirm our observations that the GCV and the LC yield parameters which are closest to the (theoretical) optimal parameter. We obtain almost equally good results for the DP, the RM and the RGCV.

%Conclusions
\section{Conclusion and outlook}
We tested parameter choice methods for the regularized (orthogonal) functional matching pursuit (RFMP/ROFMP). 
For the evaluation of the parameter choice methods, we constructed eleven different test cases with different satellite heights, data grids, noise types and noise-to-signal ratios (see \cref{tab:testcases}) for the RFMP and ROFMP. For each test case, we generated $32$ noisy data sets. Altogether we ran each algorithm for $352$ data sets and for each data set for $100$ different regularization parameters, that means each algorithm was applied 35200 times. 

Our study shows that the GCV, the LC, the RM, the RGCV and the SRGCV yield the best results in all test cases. The DP provides good to acceptable results. The performance of the QOC seriously depends on the data grid, that means a less regularly distributed grid does not lead to good results. In our experiments, the QOC had good results with the Reuter grid. The MGCV also obtains both good and rather poor results in dependency on the grid and kind of noise we used. Here, the irregularly distributed scattered grid and the coloured noise did not yield good results. At last, the TDP, the EEM and the GML did not always lead to good results in our test cases. 

We want to remark that in average our results were better than in \cite{bauer_gutting_lukas} and \cite{bauer_lukas} for all methods. Some possible reasons for that can be: the coloured noise in our test cases was different and maybe easier to handle for the methods than in the two papers, because we only had an AR(1)-process.
There is a further difference to the other cases in relation to the problem itself. Here we had a data grid given which corresponds to a spatial discretization of the problem. Furthermore, the RFMP and the ROFMP are iterative methods and use stopping criteria which are also some kind of regularization. Since we stop the algorithm at a certain point we do not obtain the approximation of the potential in the limit. 
For these reasons, the outcomes of our experiments and of those in \cite{bauer_gutting_lukas, bauer_lukas} are not really comparable.

The purpose of this paper is to provide a first guideline for the parameter choice for the RFMP and the ROFMP. Certainly, further experiments should be designed in the future.
Maybe, the distribution of our regularization parameters $\lambda_{k}$ could be improved such that the relevant parameters themselves are not too wide apart. Perhaps, the interval from $1$ to $10^{-14}$ should be chosen smaller such that the parameters are closer together. 

Future changes in the implementation could also be the use of other stopping criteria for the RFMP. Furthermore, an enhancement could be the extension of the dictionary to localized trial functions. In addition, the generation of the coloured noise can, for example, use an $AR(k)$-process for $k > 1$ or completely different types of noise can be considered. 
Finally, we can test other tuning parameters for the methods as far as these are required.
Besides, it is possible that the performance of the investigated parameter choice methods in the RFMP/ROFMP depends on the considered inverse problem.

\bibliographystyle{abbrv}
\bibliography{biblio}
\end{document}